\documentclass[a4paper,11pt]{amsart}

\usepackage{fullpage,amsxtra,amssymb,amsthm,amsmath,amscd,xypic,url}
\usepackage[latin1]{inputenc}
\numberwithin{equation}{section}

\renewcommand{\leq}{\leqslant}
\renewcommand{\geq}{\geqslant}

\def\stacksum#1#2{{\stackrel{{\scriptstyle #1}}
{{\scriptstyle #2}}}}

\newcommand{\Cc}{\mathbf{C}}

\newcommand{\Aa}{\mathbf{A}}
\newcommand{\Bb}{\mathbf{B}}
\newcommand{\Zz}{\mathbf{Z}}

\newcommand{\Rr}{\mathbf{R}}
\newcommand{\Gg}{\mathbf{G}}

\newcommand{\Qq}{\mathbf{Q}}
\newcommand{\Fp}{\mathbf{F}}
\newcommand{\Tt}{\mathbf{T}}
\newcommand{\G}{\mathbf{G}}

\newcommand{\Lb}{\mathcal{L}}
\newcommand{\proba}{\mathbf{P}}
\newcommand{\expect}{\mathbf{E}}
\newcommand{\variance}{\mathbf{V}}
\newcommand{\charfun}{\mathbf{1}}
\newcommand{\bon}{\mathcal{B}}

\newcommand{\yc}{y^{\sharp}}

\newcommand{\mods}[1]{\,(\mathrm{mod}\,{#1})}

\newcommand{\ra}{\rightarrow}
\newcommand{\lra}{\longrightarrow}

\newcommand{\fleche}[1]{\stackrel{#1}{\lra}}

\newcommand{\barre}[1]{\overline{{#1}}}

\DeclareMathOperator{\spec}{Spec}

\DeclareMathOperator{\rank}{rank}
\DeclareMathOperator{\res}{Res}
\DeclareMathOperator{\reg}{reg}

\DeclareMathOperator{\Imag}{Im}

\DeclareMathOperator{\frob}{\mathrm{Fr}}
\DeclareMathOperator{\Gal}{Gal}
\DeclareMathOperator{\Ind}{Ind}

\DeclareMathOperator{\Tr}{Tr}
\DeclareMathOperator{\Hom}{Hom}
\DeclareMathOperator{\End}{End}

\DeclareMathOperator{\swan}{Swan}

\newcommand{\eps}{\varepsilon}

\newcommand{\demi}{{\textstyle{\frac{1}{2}}}}

\newcommand{\sheaf}[1]{\mathcal{{#1}}}

\DeclareMathSymbol{\gena}{\mathord}{letters}{"3C}
\DeclareMathSymbol{\genb}{\mathord}{letters}{"3E}

\def\sumb{\mathop{\sum \Bigl.^{\flat}}\limits}
\def\dblsum{\mathop{\sum \sum}\limits}

\def\sums{\mathop{\sum \Bigl.^{*}}\limits}

\theoremstyle{plain}
\newtheorem{theorem}{Theorem}[section]

\newtheorem{lemma}[theorem]{Lemma}
\newtheorem{corollary}[theorem]{Corollary}

\newtheorem{proposition}[theorem]{Proposition}

\theoremstyle{remark}
\newtheorem{remark}[theorem]{Remark}

\theoremstyle{definition}

\newtheorem{definition}[theorem]{Definition}

\newtheorem{example}[theorem]{Example}

\newcommand{\ggeom}{G^g}
\newcommand{\sieve}{\Psi}
\newcommand{\siftable}{\Upsilon}

\begin{document}

\title{The principle of the large sieve}
\author{E. Kowalski}
\address{Universit\'e Bordeaux I - A2X\\
351, cours de la Lib\'eration\\
33405 Talence Cedex\\
France}
\email{emmanuel.kowalski@math.u-bordeaux1.fr}
\dedicatory{Pour les soixantes ans de J-M. Deshouillers}

\date{\today}
\begin{abstract}
  We describe a very general abstract form of sieve based on a
  large sieve inequality which generalizes both the classical sieve
  inequality of Montgomery (and its higher-dimensional variants), and
  our recent sieve for Frobenius over function fields. The general
  framework suggests new applications. We give some first results on the
  number of prime divisors of ``most'' elements of an elliptic
  divisibility sequence, and we develop in some detail
  ``probabilistic'' 
  sieves for random walks on arithmetic
  groups, e.g., estimating the probability of finding a reducible
  characteristic polynomial at some step of  a random walk on
  $SL(n,\mathbf{Z})$.   In addition to the sieve principle, the
  applications depend on bounds for a large sieve constant. To prove
  such bounds involves a variety of deep results, including Property
  $(\tau)$ or expanding properties of Cayley graphs, and the Riemann
  Hypothesis over finite fields. 
\end{abstract}

\maketitle
\tableofcontents

\section{Introduction}

Classical sieve  theory is concerned with the problem of the
asymptotic evaluation of averages of  arithmetic
functions over integers constrained by congruence restrictions modulo
a set of primes. Often the function in question is the characteristic
function of some interesting sequence and the congruence restrictions
are chosen so that those integers remaining after the sieving process
are, for instance, primes or ``almost'' primes. 
\par
If the congruence conditions are phrased as stating that the only
integers $n$ which are allowed are those with reduction modulo a prime
$p$ not in a certain set $\Omega_p$, then a familiar dichotomy arises:
if $\Omega_p$ contains few residue classes (typically, a bounded
number as $p$ increases), the setting is that of a ``small'' sieve.
The simplest such case is the detection of primes with
$\Omega_p=\{0\}$. If, on the other hand, the size of $\Omega_p$
increases, the situation is that of a ``large'' sieve. The first such
sieve was devised by Linnik to investigate the question of Vinogradov
of the size of the smallest quadratic non-residue modulo a prime.
\par
There have already been a number of works extending ``small'' sieves
to more general si\-tuations, where the objects being sifted are not
necessarily integers. Among these, one might quote the vector sieve of
Br\"udern and Fouvry~\cite{brudern-fouvry}, with applications to
Lagrange's theorem with 
almost prime variables, the ``crible étrange'' of Fouvry and
Michel~\cite{fouvry-michel}, with applications to sign changes of
Kloosterman sums, and 
Poonen's striking sieve
procedure for finding smooth hypersurfaces of large degree over finite
fields~\cite{poonen}.
\par
Similarly, the large sieve has been extended in some ways, in
particular (quite early on) to deal with sieves in $\Zz^d$, $d\geq 1$,
or in number fields (see
e.g.~\cite{gallagher}). Interesting applications have been found,
e.g. Duke's theorem on elliptic curves over $\Qq$ with ``maximal''
$p$-torsion fields for all $p$~\cite{duke}. 
All these were much of the same flavor however, and in particular
depended only on the character theory of finite abelian groups as far
as the underlying  harmonic analysis was concerned\footnote{\ There is,
  of course, an enormously important body of work concerning 
  inequalities traditionally called ``large sieve inequalities'' for
  coefficients of automorphic forms of 
  various types which have been developed by Iwaniec,
  Deshouillers-Iwaniec, Duke, Duke-Kowalski, Venkatesh and
  others (a short survey is in~\cite[\S 7.7]{ant}). However, those
  generalize the large sieve inequality for 
  Dirichlet characters, and have usually no relation (except
  terminological) with the traditional sieve principle.}.
\par
In~\cite{sieve}, we have introduced a new large sieve type inequality
for the average distribution of Frobenius conjugacy classes in the
monodromy groups of a family ($\sheaf{F}_{\ell})$ of $\Fp_{\ell}$-adic
sheaves on a variety over a finite field. Although the spirit of the
large sieve is clearly recognizable, the setting is very different,
and the harmonic analysis involves both non-abelian finite groups, and
the deep results of Deligne on the Riemann Hypothesis over finite
fields. Our first application of this new sieve was
related to the ``generic'' arithmetic behavior of the numerator of the
zeta function of a smooth projective curve in a family with large
monodromy, improving significantly a result of
Chavdarov~\cite{chavdarov}. 
\par
Motivated by this first paper, the present one is interested with
foundational issues related to the large sieve. We are able to
describe a very general abstract framework which we call ``the
principle of the large sieve'', with a pun
on~\cite{montgomery-survey}. This 
leads to a sieve 
statement that may in particular be specialized to either the
classical forms of the large sieve, or to a strengthening
of~\cite{sieve}. Roughly speaking, we 
deal with a set $X$ that can be mapped to finite sets
$X_{\ell}$ (for instance, integers can be reduced modulo primes)
and we show how an estimate for the number of those $x\in
X$ which have ``reductions'' outside $\Omega_{\ell}\subset X_{\ell}$
for all or some $\ell$ may be reduced to a bilinear form
estimate of a certain kind. The form of the sieve statement we obtain is 
similar to Montgomery's formulation of the large sieve (see 
e.g.~\cite{montgomery-survey},~\cite{bombieri},~\cite[7.4]{ant}). It
should be mentioned that our ``axioms'' for the sieve may admit other
variations. In fact, Zywina~\cite{zywina} 
has developed a somewhat similar framework, and some of the
flexibility we allow was first suggested by his
presentation.
\par
There remains the
problem of estimating the 
bilinear form.  The classical idea of duality and exponential sums is
one tool in this direction, and we describe it also somewhat
abstractly. We then find a convincing relation with the classical
sieve axioms, related to equidistribution in the finite sets
$X_{\ell}$. 
\par
The bilinear form inequality also seemingly depends on the choice of an
orthonormal basis of certain finite-dimensional Hilbert spaces. 
It turns out that in many applications, the sieve setting is
related to the existence of a group $G$ such that $X_{\ell}$ is the
set of conjugacy classes in a
finite quotient $G_{\ell}$ of $G$ and the reduction $X\ra X_{\ell}$ factors
through $G$. In that case, the bilinear form inequality can be stated
with a distinguished basis arising from the representation theory (or
harmonic analysis) of the finite groups $G_{\ell}$. 
\par
This abstract sieving framework has many
incarnations. As we already stated, we can  recover the 
classical large sieve and the ``sieve for Frobenius''
of~\cite{sieve}, but furthermore, we are led to a number of situations
which are either new (to the author's knowledge), or have
received attention only recently, although not in the same form in
general. One of these concerns (small) sieves in
arithmetic groups and is the subject of ongoing work of Bourgain, Gamburd
and Sarnak~\cite{bgs}, and some of the problems it is suited for have
been raised and partly solved by Rivin~\cite{rivin}, who also
emphasized possible applications to some groups which are ``close in
spirit'' to arithmetic groups, such as mapping class groups of
surfaces or automorphism groups of free groups. Indeed, the large
sieve strengthens significantly the results of Rivin (see
Corollary~\ref{cor-anosov}).
\par
\medskip
Our main interest in writing this paper is the exploration of
the general setting. Consequently, the paper is fairly open-ended and
has a distinctly chatty style. We hope to
come back to some of the new examples 
with more applications in the future. Still, to give a feeling for the
type of results that become available, we finish this introduction
with a few sample statements (the last one could in fact have been
derived in~\cite{sieve}, with a slightly worse bound).

\begin{theorem}\label{th-proba}
Let $(S_n)$ be a simple random walk on $\Zz$, i.e.,
$$
S_n=X_1+\cdots +X_n
$$
where $(X_k)$ is a sequence of independent random variables with
$\proba(X_k=\pm 1)=\demi$ for all $k$.
\par
Let $\eps>0$ be given, $\eps\leq 1/4$. For any odd $q\geq
1$, any $a$ coprime with $q$, we have 
$$
\proba(S_n\text{ is prime and } \equiv a\mods{q})\ll \frac{1}{\varphi(q)}
\frac{1}{\log n}
$$
if $n\geq 1$, $q\leq n^{1/4-\eps}$, the implied constant depending
only on $\eps$.
\end{theorem}


\begin{theorem}\label{th-sln}
Let $n\geq 2$ be an integer, let $G=SL(n,\Zz)$ and let
$S=S^{-1}\subset G$ be a finite generating set of $G$, e.g., the
finite set of elementary matrices with $\pm 1$ entries off the
diagonal.
Let $(X_k)$ be the simple left-invariant random walk on $\Gamma$, i.e.,
a sequence of $\Gamma$-valued random variables such that $X_0=1$ and
$$
X_{k+1}=X_k\xi_{k+1}\text{ for } k\geq 0,
$$
where $(\xi_k)$ is a sequence of $S$-valued independent random
variables with $\proba(\xi_k=s)=\frac{1}{|S|}$ for all $s\in S$. 
Then, almost surely, there are only finitely many $k$ for which the
characteristic polynomial $\det(X_k-T)\in \Zz[T]$ is reducible, or in
other words, the set of matrices with reducible characteristic
polynomials in $SL(n,\Zz)$ is transient for the random walk.
\end{theorem}

In fact (see Theorem~\ref{th-sln-bis}), we will derive this by showing
that the probability that $\det(X_k-T)$ be reducible decays 
exponentially fast with $k$ (in the case $n\geq 3$ at least).  An
analogue of this result (with some extra conditions) has the
geometric/topological consequence that the set of non-pseudo-Anosov
elements is transient for random walks on mapping class
groups of closed orientable surfaces, answering a question of
Maher~\cite[Question 1.3]{maher} (see Corollary~\ref{cor-anosov} 
for details; this application was suggested by Rivin's
paper~\cite{rivin}).

\begin{theorem}\label{th-non-conj}
Let $n\geq 3$ be an integer, let $G=SL(n,\Zz)$, and let
$S=S^{-1}\subset G$ be a finite symmetric generating set. 
Then there exists $\beta>0$ such that for any $N\geq 1$, we have 
$$
|\{w\in S^N\,\mid\, 
\text{one entry of the matrix $g_w$ is a square}\}|\ll |S|^{N(1-\beta)},
$$ 
where $g_w=s_1\cdots s_N$ for $w=(s_1,\ldots,s_N)\in S^N$, and $\beta$
and the implied constant depend only on $n$ and $S$. 
\par
Equivalently, for the random walk $(X_k)$ on $G$ defined as in the
statement of the previous theorem, we have
$$
\proba(\text{one entry of the matrix $X_k$ is a square})
\ll \exp(-\delta k)
$$
for $k\geq 1$ and some constant $\delta>0$, where $\delta$ and the
implied constant depend only on $n$ and $S$.
\end{theorem}


\begin{theorem}\label{th-ell}
Let $E/\Qq$ be an elliptic curve with rank $r\geq 1$ given by a
Weierstrass equation
$$
y^2+a_1xy+a_3y=x^3+a_2x^2+a_4x+a_6,\quad\text{ where } a_i\in\Zz.
$$
For $x\in
E(\Qq)$, let $\omega_E(x)$ be the number of primes, without
multiplicity, dividing the denominator of the coordinates of $x$, with
$\omega_E(0)=+\infty$. Let $h(x)$ denote the canonical height on $E$.
\par
Then for any
fixed real number $\kappa$ with $0<\kappa<1$, we have
$$
|\{x\in E(\Qq)\,\mid\, h(x)\leq T\text{ and }
\omega_E(x)<\kappa \log\log T
\}|\ll T^{r/2}(\log\log T)^{-1},
$$
for $T\geq 3$, where the implied constant depends only on $E$ and
$\kappa$. 
\end{theorem}

\begin{theorem}\label{th-sq-jac}
Let $q$ be a power of a prime number $p\geq 5$, $g\geq 1$ an integer
and let $f\in \Fp_q[T]$
be a squarefree polynomial of degree $2g$. For $t$ not a zero of
$f$, let $C_t$ denote the smooth projective model of the
hyperelliptic curve
$$
y^2=f(x)(x-t),
$$
and let $J_t$ denote its Jacobian variety.
Then we have
\begin{gather*}
|\{
t\in \Fp_q\,\mid\, f(t)\not=0\text{ and 
$|C_t(\Fp_q)|$ is a square} \}|\ll
q^{1-\gamma}(\log q),\\
|\{
t\in \Fp_q\,\mid\, f(t)\not=0\text{ and 
$|J_t(\Fp_q)|$ is a square} \}|\ll
q^{1-\gamma}(\log q)
\end{gather*}
where $\gamma=(4g^2+2g+4)^{-1}$, and the implied constants are absolute.
\end{theorem}

It is well-known that the strong form of the large sieve is as 
efficient (qualitatively) as the best small sieves, as far as upper
bound sieves are concerned. To put this in context, we will briefly
recall the principles of small sieves (in the same abstract context)
in an Appendix, and we will give a sample application
(Theorem~\ref{th-lower-sieve}) related to Theorem~\ref{th-sq-jac}. 
\par
\medskip
\par
The plan of this paper is as follows. In the first sections, the
abstract sieve setting is described, and the abstract large sieve
inequality is derived; this is a pleasant and rather
straightforward algebraic exercise. In
Sections~\ref{sec-spec-group} and~\ref{sec-spec-coset}, we specialize
the general  
setting to two cases (``group sieve'' and ``coset sieve'') related to
group theory, using the representation theory of finite groups.
This leads to the natural problem of finding precise estimates for
the degree and the sum of degrees of
irreducible representations of some finite groups of Lie type, which
we consider in some cases in Section~\ref{sec-reprs}. For
this we use Deligne-Lusztig characters, and arguments shown to the
author by J. Michel; this section may be omitted in a first reading.
\par
Turning to examples of sieves, already in Section~\ref{sec-examples}
we show how many classically-known uses of the large sieve are special
cases of the setting of Section~\ref{sec-spec-group}. In the same
section, we also indicate the relation with the inclusion-exclusion
technique in probability  and combinatorics, which shows in particular
that the general sieve bound is sharp (see Example~\ref{ex-excl-incl}).
\par
New
(or emerging) situations are considered next, in four sections which
are quite independent of one another (all of them involve either group
or coset sieves). ``Probabilistic'' sieves 
are discussed briefly in Section~\ref{sec-proba-sieve}, leading to
Theorem~\ref{th-proba}. Sieving in 
arithmetic groups is described in Section~\ref{sec-sieve-arith}, where
Theorem~\ref{th-sln} is proved. The crucial point (as in the work of
Bourgain, Gamburd and Sarnak) is  the expanding properties of Cayley
graphs  of $SL(n,\Zz/d\Zz)$, phrased in terms of Property
$(\tau)$. Then comes  an  amusing  
``elliptic sieve'' which is related to the number of prime divisors of
the denominators of rational points on an elliptic
curve, leading to Theorem~\ref{th-ell}. In turn, this is linked to
the analysis of the prime 
factorization of elements of so-called ``elliptic divisibility
sequences'', and we find  that ``most'' elements have many prime
factors. This  
complements recent heuristics and results of Silverman, 
Everest, Ward and others concerning the paucity of primes and prime
powers in such sequences. 
Finally, in
Section~\ref{sec-ls-ff}, we extend the sieve result of~\cite{sieve} 
concerning the distribution of geometric Frobenius conjugacy classes
in finite monodromy groups over finite fields, and derive some new
applications. 
To conclude, Appendix A briefly indicates the link with small
sieve situations, for the purpose of comparison and reference, with
a sample application, and Appendix B contains the proofs of some
``local'' density computations in matrix groups over finite
fields. Those estimates have been used previously, but we defer the
proof to not distract from the main thrust of the arguments underlying
the principle of the sieve. Note that the techniques underlying those
computations 
are in fact quite advanced and of independent interest, and involve
work of Chavdarov~\cite{chavdarov} and non-trivial estimates for
exponential sums over finite fields.

\medskip
\textbf{Notation.} As usual, $|X|$ denotes the cardinality of a
set; however if $X$ is a measure space with measure $\mu$, we
sometimes write $|X|$ instead of $\mu(X)$. 
\par
For a group $G$, $G^{\sharp}$
denotes the set of its conjugacy 
classes, and for a conjugacy-invariant subset $X\subset G$,
$X^{\sharp}\subset G^{\sharp}$ is the corresponding set of conjugacy
classes. The conjugacy class of $g\in G$ is denoted $g^{\sharp}$.
\par
By $f\ll g$ for $x\in X$, or $f=O(g)$ for $x\in X$, where $X$ is an
arbitrary set on which $f$ is defined, we mean synonymously that there
exists a constant $C\geq 0$ such that $|f(x)|\leq Cg(x)$ for all $x\in
X$. The ``implied constant'' is any admissible value of $C$. It may
depend on the set $X$ which is always specified or clear in context.
The notation $ f\asymp g$ means $f\ll g$ and $g\ll f$.
On the other hand $f(x)=o(g(x))$ as $x\ra x_0$
is a topological statement meaning that $f(x)/g(x)\ra 0$ as $x\ra
x_0$.
\par
For $n\geq 1$ an integer, $\omega(n)$ is the number of primes dividing
$n$, without counting multiplicity. For $z\in\Cc$, we denote
$e(z)=\exp(2i\pi z)$.
\par
In probabilistic contexts, $\proba(A)$ is the probability of an event,
$\expect(X)$ is the expectation of a random variable $X$,
$\variance(X)$ its variance, and
$\charfun_A$ is the characteristic function of an event $A$.
\par
\medskip
\textbf{Acknowledgments.} D. Zywina has developed~\cite{zywina} an
abstract setup of the 
large sieve similar to the conjugacy sieve described in
Section~\ref{sec-spec-group}. His remarks have 
been very helpful both for the purpose of straightening out the
assumptions used, and as motivation for the search of new ``unusual''
applications. One of his nice tricks (the use of general sieve
support) is also used here.  
The probabilistic setting was suggested in part by
Rivin's preprint~\cite{rivin}, who also mentioned to me the work of
Bourgain, Sarnak and Gamburd.
I also wish to thank P. Sarnak for sending me a copy
of his email~\cite{sarnak-email} to his coauthors. Finally, I thank
J. Michel for  
providing the ideas of the proof of Proposition~\ref{pr-degrees} and
explaining some basic properties of representations of finite groups of
Lie type, and P. Duchon and M-L. Chabanol for help, advice and
references concerning probability theory and graph theory.
\par

\section{The principle of the large sieve}
\label{sec-prelim}

We will start by describing a very general type of sieve. The goal is
to reach an analogue of the large sieve inequality, in the sense of a
reduction of a sieve bound to a bilinear form estimate.
\par
We start by introducing the notation and terminology. The \emph{sieve
  setting}\label{pg-setting} is a triple $\sieve=(Y,\Lambda,(\rho_{\ell}))$
consisting of
\begin{itemize}
\item A  set  $Y$;
\item An index set $\Lambda$;
\item For all $\ell\in\Lambda$, a surjective map $\rho_{\ell}\,:\,
  Y\ra Y_{\ell}$ where $Y_{\ell}$ is a finite set.
\end{itemize}
\par
In combinatorial terms, this might be thought as a family of colorings
of the set $Y$. In applications, $\Lambda$ will  often be a
subset of primes (or prime ideals in 
some number field), but
as first pointed out by Zywina, this is not necessary for the formal
part of setting up the sieve, and although the generality is not
really abstractly greater, it is convenient to allow arbitrary
$\Lambda$.  
\par
Then, a \emph{siftable set}\label{pg-siftable} associated to
$\sieve=(Y,\Lambda,(\rho_{\ell}))$ is a triple
$\siftable=(X,\mu,F)$ consisting of
\begin{itemize}
\item A measure space $(X,\mu)$ with $\mu(X)<+\infty$;
\item A map $F\,:\, X\ra Y$ such that the composites $X\ra Y\ra
  Y_{\ell}$ are measurable, i.e., the sets $\{x\in X\,\mid\,
  \rho_{\ell}(F_x)=y\}$ are measurable for all $\ell$ and all $y\in
  Y_{\ell}$. 
\end{itemize}
\par
The simplest case is when $X$ is a finite set and $\mu$ is counting
measure. We call this the \emph{counting case}. Even when this is not
the case,  for notational
convenience, we will usually write $|B|=\mu(B)$
for the measure of a measurable set $B\subset X$.
\par
The last piece of data is a finite subset $\Lb^*$ of $\Lambda$, called
the \emph{prime sieve support}\label{pg-prime-sieve-support}, and a
family $\Omega=(\Omega_{\ell})$ of \emph{sieving
  sets\footnote{\ Sometimes, $\Omega$ will also denote a probability
    space, but no confusion should arise.}} of $Y_{\ell}$, defined for $\ell\in \Lb^*$.
\par
With this final data $(\sieve,\siftable,\Lb^*,\Omega)$, we can define
the sieve problem. 

\begin{definition}\label{def-sieve}
Let $\sieve=(Y,\Lambda,(\rho_{\ell}))$ be a sieve setting,
$\siftable=(X,\mu,F)$ a siftable set, $\Lb^*$ a prime sieve support
and   $\Omega$ a family of sieving sets. Then the  \emph{sifted
  sets}\label{pg-sifted} are
\begin{gather*}
S(Y,\Omega;\Lb^*)=\{y \in Y\,\mid\, \rho_{\ell}(y)\notin
\Omega_{\ell}\text{ for all } \ell\in \Lb^*\},
\\
S(X,\Omega;\Lb^*)=\{x\in X\,\mid\, \rho_{\ell}(F_x)\notin
\Omega_{\ell}\text{ for all } \ell\in \Lb^*\}.
\end{gather*}
The latter is also $F^{-1}(S(Y,\Omega;\Lb^*))$ and is a measurable
subset of $X$.
\end{definition}

The problem we will consider is to find estimates for the measure
$|S(X,\Omega;\Lb^*)|$ of the sifted set. Here we think that the sieve
setting is fixed, while there usually will be an infinite sequence of
siftable sets with size  $|X|$ going to infinity; this size will be
the main variable in the estimates. 

\begin{example}\label{ex-classical-new}
The classical sieve arises as follows: the sieve setting is
$$
\sieve=(\Zz,\{\text{primes}\},\Zz\ra \Zz/\ell\Zz)
$$
and the siftable sets are $X=\{n\,\mid\, M<n\leq M+N\}$ with counting
measure and $F_x=x$ for $x\in X$. Then the sifted sets become the
classical sets of integers in an interval with reductions modulo
primes in $\Lb^*$ lying outside a subset $\Omega_{\ell}\subset
\Zz/\ell\Zz$ of residue classes.
\par
In most cases, $(X,\mu)$ will be a finite set with counting
measure, and often $X\subset Y$ with $F_x=x$ for $x\in X$. See 
Section~\ref{sec-ls-ff} for a conspicuous example where $F$ is not the
identity, Section~\ref{sec-proba-sieve} for interesting situations
where the measure space $(X,\mu)$ is a probability space, and $F$ a random
variable, and  Section~\ref{sec-sieve-arith} for another
example.
\end{example}

We will now indicate one type of inequality that reduces the sieve
problem to the estimation of a \emph{large sieve constant}
$\Delta$. The latter is a more analytic problem, and can be attacked
in a number of ways. This large sieve constant depends on most of the
data involved, but is independent of the sieving sets.
\par
First we need some more notation. Given a sieve setting $\sieve$, we
let $S(\Lambda)$ denote the set of finite subsets $m\subset
\Lambda$. Since $S(\Lambda)$ may be identified with the set of
squarefree integers $m\geq 1$ in the classical case where $\Lambda$ is
the set of primes, to simplify notation we write $\ell\mid m$ for
$\ell\in m$ when $\ell\in\Lambda$ and $m\in S(\Lambda)$, and similarly
for $n\mid m$ instead of $n\subset m$ if $n$, $m\in S(\Lambda)$.
\par
A \emph{sieve support} $\Lb$ associated to a prime sieve support
$\Lb^*$ is any finite subset of $S(\Lambda)$ such that 
\begin{equation}\label{eq-sieve-support}
\ell\in m,\ m\in \Lb\text{ implies } \ell\in\Lb^*,\text{ and }
\{\ell\}\in \Lb\text{  if }\ell\in \Lb^*.
\end{equation}
\par
This implies that $\Lb$ determines $\Lb^*$ (as the set of elements of
singletons in $\Lb$).
If $\Lambda$ is a set of primes, $\Lb$ ``is'' a set of squarefree integers
only divisible by primes in $\Lb^*$ and containing
$\Lb^*$ (including possibly $m=1$, not divisible by any prime). 
\par
For $m\in S(\Lambda)$, let
$$
Y_m=\prod_{\ell\mid m}{Y_{\ell}}
$$
and let $\rho_m\,:\, Y\ra Y_m$ be the obvious product map. (In other
words, we look at all ``refined''  colorings of $Y$ obtained by
looking at all possible finite tuples of colorings). If $m=\emptyset$,
$Y_m$ is a set with a single element, and $\rho_m$ is a constant map. 
\par
We will consider functions on the various sets $Y_m$, and it will be
important to endow the space of complex-valued functions on $Y_m$ with
appropriate and consistent inner products. For this purpose, we assume
given for $\ell\in\Lambda$ a density 
$$
\nu_{\ell}\,:\, Y_{\ell}\ra [0,1]
$$
(often denoted
simply $\nu$ when no ambiguity is possible) such that the inner
product on functions $f\,:\, Y_{\ell}\ra \Cc$ is given by
$$
\langle f,g\rangle=\sum_{y\in Y_{\ell}}
{\nu_{\ell}(y)f(y)\overline{g(y)}}.
$$
\par
We assume that $\nu(y)>0$ for all $y\in Y_{\ell}$, in order that this
hermitian form be positive definite (it will be clear that $\nu(y)\geq
0$ would suffice, but the stronger assumption is no problem for
applications), and that $\nu$ 
is a probability density, i.e., we have
\begin{equation}\label{eq-proba-density}
\sum_{y\in Y_{\ell}}{\nu_{\ell}(y)}=1.
\end{equation}
\par
Using the product structure we define corresponding inner products
and measures on the spaces of functions $Y_m\ra
\Cc$. Property~(\ref{eq-proba-density}) still holds.
We will interpret $\nu$ as a measure on $Y_{\ell}$ or $Y_{m}$, so we will 
write for instance
$$
\nu(\Omega_{\ell})=\sum_{y\in\Omega_{\ell}}{\nu(y)},\quad\quad
\text{ for }\Omega_{\ell}\subset Y_{\ell}.
$$
\par
We denote by $L^2(Y_m)$ the space of complex-valued functions on $Y_m$
with the inner product thus defined.
\par
The simplest example is when
$\nu(y)=1/|Y_m|$, but see Sections~\ref{sec-spec-group}
and~\ref{sec-spec-coset} for
important natural cases  where $\nu$ is not uniform.
It will be clear in the remarks and sections following the statement
of the sieve inequality that, in general, the apparent choice of
$\nu_{\ell}$ is illusory (only one choice will lead to good results). 
\par
Note that
$\rho_m$ is not necessarily surjective, but it turns out to be true, and
a crucial fact, in most applications of the sieve, so we make a
definition (the terminology will be clearer in later applications).

\begin{definition}\label{def-chinese}
A sieve setting $\sieve=(Y,\Lambda,(\rho_{\ell}))$ is \emph{linearly
disjoint} if the map $\rho_m\,:\, Y\ra Y_m$ is onto for all $m\in
S(\Lambda)$. 
\end{definition}

Here is now the first sieve inequality.

\begin{proposition}\label{pr-comb-sieve}
Let $\sieve$, $\siftable$, $\Lb^*$ be as above. For any sieve support
$\Lb$ associated to $\Lb^*$, i.e,  any finite
subset of $S(\Lambda)$ satisfying~\emph{(\ref{eq-sieve-support})},
let $\Delta=\Delta(X,\Lb)$ denote the \emph{large sieve
constant}, which is by definition the smallest non-negative real
number such that 
\begin{equation}\label{eq-comb-ls}
\sum_{m\in\Lb}{\sum_{\varphi\in\bon_m^*}{
\Bigl|
\int_{X}{\alpha(x)\varphi(\rho_m(F_x))d\mu(x)}
\Bigr|^2
}}\leq \Delta\int_X{|\alpha(x)|^2d\mu(x)}
\end{equation}
for any square integrable function $\alpha\,:\, X\ra \Cc$,
where $\varphi$ in the outer sum ranges over $\bon_m^*$, where
$\bon_{\ell}^*=\bon_{\ell}-\{1\}$,
$\bon_{\ell}$ is an orthonormal basis, containing the constant
function $1$, of the space $L^2(Y_{\ell})$,
and for all $m$ we let
$$
\bon_m=\prod_{\ell\mid m}{\bon_{\ell}},\quad\quad 
\bon_m^*=\prod_{\ell\mid m}{\bon_{\ell}^*},
$$
the function on $Y_m$ corresponding to $(\varphi_{\ell})$ being given
by
$$
(y_{\ell})\mapsto \prod_{\ell\mid m}{\varphi_{\ell}(y_{\ell})},
$$
and for $m=\emptyset$, we have $\bon_m=\bon_m^*=\{1\}$.
\par
Then for arbitrary sieving sets $\Omega=(\Omega_{\ell})$, we have
$$
|S(X,\Omega;\Lb^*)|\leq \Delta H^{-1}
$$
where
\begin{equation}\label{eq-h}
H=\sum_{m\in\Lb}{\prod_{\ell\mid m}{\frac{\nu(\Omega_{\ell})}
{\nu(Y_{\ell}-\Omega_{\ell})}}}=
\sum_{m\in\Lb}{\prod_{\ell\mid m}{\frac{\nu(\Omega_{\ell})}
{1-\nu(\Omega_{\ell})}}}.
\end{equation}
\end{proposition}

\begin{remark}
The large sieve constant as defined above is independent of the choices
of basis $\bon_{\ell}$ (containing the constant function
$1$). Here is a more intrinsic definition which shows this, and 
provides a first hint of the link with classical (small) sieve axioms.
It's not clear how much this intrinsic definition can be useful in
practice, which explains why we kept a concrete version in the
statement of Proposition~\ref{pr-comb-sieve}.
\par
By definition, the inequality~(\ref{eq-comb-ls}) means that $\Delta$
is the square of the norm of the linear operator
$$
T\ 
\left\{
\begin{matrix}
L^2(X,\mu)&\longrightarrow
&\displaystyle{\bigoplus_{m\in\Lb}{L(L_0^2(Y_m),\Cc)}}
\\
\alpha&\mapsto&\displaystyle{ \Bigl(f\mapsto 
\int_X{\alpha(x)f(\rho_m(F_x))d\mu(x)}
\Bigr)_m}
\end{matrix}
\right.
$$
where the direct sum over $m$ is orthogonal and $L(L^2_0(Y_m),\Cc)$ is
the space of linear functionals on 
$$
L^2_0(Y_m)=\bigotimes_{\ell\mid m}{L^2_0(Y_{\ell})},\quad \text{ where }\quad
L^2_0(Y_{\ell})=\{f\in L^2(Y_{\ell})\,\mid\, \langle f,1\rangle=
\sum_y{\nu(y)f(y)}=0\}
$$
(the space $L^2_0(Y_m)$ may be thought of as the ``primitive''
subspace of the functions on $Y_m$), with the norm
$$
\|f^*\|=\max_{f\not=0}\frac{|\langle f^*,f\rangle|}{\|f\|}.
$$
\par
Since we are dealing with Hilbert spaces, $L(L^2_0(Y_m),\Cc)$ is
canonically isometric to $L^2_0(Y_m)$, and $\Delta$ is the square of
the norm of the operator 
$$
T_1
\left\{
\begin{matrix}
L^2(X,\mu)&\longrightarrow &
\displaystyle{\bigoplus_{m\in\Lb}{L_0^2(Y_m)}}\\
\alpha&\mapsto& T_1(\alpha)
\end{matrix}
\right.
$$
where $T_1(\alpha)$ is the vector such that $\langle
f,T_1(\alpha)\rangle =T(\alpha)(f)$ for $f\in L^2_0(Y_m)$,
$m\in\Lb$. This vector is easy to identify: we have
$$
\int_X{\alpha(x)f(\rho_m(F_x))d\mu(x)}=
\sum_{y\in
  Y_{\ell}}{f(y)\Bigl(\int_{\{\rho_m(F_x)=y\}}{\alpha(x)d\mu(x)}
\Bigr)}
$$
which means that $T_1(\alpha)$ is the complex-conjugate of the
projection to $L^2_0(Y_m)$ of the function
$$
y\mapsto \frac{1}{\nu_m(y)}
\int_{\{\rho_m(F_x)=y\}}{\alpha(x)d\mu(x)}
$$
on $Y_{m}$. For $m=\{\ell\}$, this projection is
obtained by subtracting the contribution of the  constant function,
i.e., subtracting the average over $y$: it is 
\begin{align*}
y\mapsto  &\frac{1}{\nu(y)}
\int_{\{\rho_m(F_x)=y\}}{\alpha(x)d\mu(x)}-
\sum_{y}{\int_{\{\rho_m(F_x)=y\}}{\alpha(x)d\mu(x)}}
\\
&=\frac{1}{\nu(y)}\int_{\{\rho_m(F_x)=y\}}{\alpha(x)d\mu(x)}-
\int_{X}{\alpha(x)d\mu(x)}.
\end{align*}
\par
In the case of counting measure and a uniform density $\nu$, this
becomes the quantity
$$
\sum_{\rho_m(F_x)=y}{\alpha(x)}-\frac{1}{|Y_{\ell}|}
\sum_{x}{\alpha(x)}
$$
after multiplying by $\nu(y)$, which is a typical ``error term''
appearing in sieve axioms. 
\end{remark}

To prove Proposition~\ref{pr-comb-sieve}, 
we start with two lemmas. For $m\in S(\Lambda)$, $y\in Y_m$, an
element $\varphi$ of the basis $\bon_m$, and a square-integrable
function $\alpha\in L^2(X,\mu)$, we denote  
\begin{equation}\label{eq-defs}
S(m,y)=\int_{\{\rho_m(F_x)=y\}}{\alpha(x)d\mu(x)},\quad\text{ and }
\quad
S(\varphi)=\int_X{\alpha(x)\varphi(\rho_m(F_x))d\mu(x)},
\end{equation}
where the integral is defined because $\mu(X)<+\infty$ by assumption.
The first lemma is the following: 

\begin{lemma}\label{lm-1b}
We have for all $\ell\in\Lambda$ the relation
$$
\sum_{\varphi\in\bon_{\ell}^*}{|S(\varphi)|^2}=
\sum_{y\in  Y_{\ell}}
{\frac{|S(\ell,y)|^2}{\nu(y)}}
-\Bigl|\int_X{\alpha(x)d\mu(x)}\Bigr|^2.
$$
\end{lemma}

\begin{proof}
Expanding the square by Fubini's Theorem, the left-hand side is
$$
\int_X\int_X{\alpha(x)\overline{\alpha(y)}
\sum_{\varphi\in\bon_{\ell}^*}{\varphi(\rho_{\ell}(F_x))\overline{
\varphi(\rho_{\ell}(F_y))}
d\mu(x)d\mu(y)}}.
$$
\par
Since $(\varphi)_{\varphi\in\bon_{\ell}}$ is an orthonormal basis of
the space of functions on $Y_{\ell}$, expanding the delta function
$z\mapsto \delta(y,z)$ in the basis gives
$$
\sum_{\varphi\in\bon_{\ell}}{\varphi(y)\overline{\varphi(z)}}
=\frac{1}{\nu(y)}\delta(y,z).  
$$
\par
Taking on the right-hand side the contribution of the constant
function $1$, we get in particular
$$
\sum_{\varphi\in\bon_{\ell}^*}
{\varphi(\rho_{\ell}(F_x))
\overline{\varphi(\rho_{\ell}(F_y))}}=\frac{1}{\nu(F_x)}
\delta(F_x,F_y)-1.
$$
Inserting this in the first relation, we obtain
\begin{align*}
\sum_{\varphi\in\bon_{\ell}^*}{|S(\varphi)|^2}&=
\int\int_{\{F_x=F_y\}}{\frac{
\alpha(x)\overline{\alpha(y)}}{\nu(F_x)}d\mu(x)d\mu(y)}-
\int_X\int_{X}{\alpha(x)\overline{\alpha(y)}d\mu(x)d\mu(y)}
\\
&=\sum_{z\in Y_{\ell}}{\frac{1}{\nu(z)}
\int\int_{\{F_x=z=
F_y\}}{\alpha(x)\overline{\alpha(y)}d\mu(x)d\mu(y)}}
-\Bigl|\int_{X}{\alpha(x)d\mu(x)}\Bigr|^2
\\
&=\sum_{y\in
  Y_{\ell}}{\frac{|S(\ell,y)|^2}{\nu(y)}}
-\Bigl|\int_{X}{\alpha(x)d\mu(x)}\Bigr|^2,
\end{align*}
as desired.
\end{proof}

Here is the next lemma.

\begin{lemma}\label{lm-2b}
Let $(\sieve,\siftable, \Omega,\Lb^*)$ be as above, and let $\Lb$ be
any sieve support associated to $\Lb^*$. For any square-integrable
function $x\mapsto \alpha(x)$ on $X$  supported
on the sifted set $S(X,\Omega;\Lb^*)\subset X$, and for any $m\in \Lb$,
we have
$$
\sum_{\varphi\in\bon_m^*}{|S(\varphi)|^2}\geq
\Bigl|
\int_{X}{\alpha(x)d\mu(x)}
\Bigr|^2
\prod_{\ell\mid m}{
\frac{\nu(\Omega_{\ell})}{\nu(Y_{\ell}-\Omega_{\ell})}},
$$
where $S(\varphi)$ is given by~\emph{(\ref{eq-defs})}.
\end{lemma}

\begin{proof}
Since this does not change the sifted set, we may 
replace $\Lb$ if necessary by the full power set of $\Lb^*$. 
Then, as in the classical case (see 
e.g.~\cite[Lemma 7.15]{ant}), the proof 
proceeds by induction on the number of elements in $m$. If
$m=\emptyset$, the inequality is trivial (there is equality, in
fact). If $m=\{\ell\}$ with $\ell\in\Lambda$ (in the arithmetic case,
$m$ is a prime), then $\ell\in\Lb^*$ by~(\ref{eq-sieve-support}). Using
Cauchy's inequality and the definition of the sifted set with the
assumption on $\alpha(x)$ to restrict the 
support of integration to elements where $\rho_{\ell}(F_x)\notin
\Omega_{\ell}$, we obtain: 
\begin{align*}
\Bigl|
\int_{X}{\alpha(x)d\mu(x)}
\Bigr|^2=\Bigl|
\sum_{\stacksum{y\in Y_{\ell}}{y\notin \Omega_{\ell}}}
{S(\ell,y)}
\Bigr|^2
&\leq \Bigl(
\sum_{y\notin \Omega_{\ell}}
{\nu(y)}
\Bigr)
\Bigl(
\sum_{y\in Y_{\ell}}
{\frac{|S(\ell,y)|^2}{\nu(y)}}
\Bigr)\\
&= 
\nu(Y_{\ell}-\Omega_{\ell})
\sum_{y\in Y_{\ell}}
{\frac{|S(\ell,y)|^2}{\nu(y)}}
\\
&=\nu(Y_{\ell}-\Omega_{\ell})
\Bigl\{
\sum_{\varphi\in\bon_{\ell}^*}{|S(\varphi)|^2}
+\Bigl|
\int_{X}{\alpha(x)d\mu(x)}
\Bigr|^2
\Bigr\}
\end{align*}
(by Lemma~\ref{lm-1b}), hence the result by moving $|\int
\alpha(x)d\mu|^2$ on the left-hand side, since $\nu(Y_{\ell})=1$.
\par
The induction step is now immediate, relying on the fact that the
function $\alpha$ is arbitrary  and the
sets $\bon_m^*$ are ``multiplicative'': for $m\in\Lb$, not a singleton,
write $m=m_1m_2=m_1\cup m_2$ with $m_1$ and $m_2$ non-empty. Then  we
have\footnote{\ Here we use the enlargement of $\Lb$ at the beginning to
  ensure that $m_i\in\Lb$.}  
$$
\sum_{\varphi\in\bon_{m_1m_2}^*}{|S(\varphi)|^2}=
\sum_{\varphi_1\in\bon_{m_1}^*}{\sum_{\varphi_2\in\bon_{m_2}^*}{
|S(\varphi_1\otimes\varphi_2)|^2
}}
$$
where $\varphi_1\otimes\varphi_2$ is the function $(y,z)\mapsto
\varphi_1(y)\varphi_2(z)$. For fixed $\varphi_1$, we can express the
inner sum as 
$$
S(\varphi_1\otimes\varphi_2)=
\int_X{\beta(x)\varphi_2(\rho_{m_2}(F_x))d\mu(x)}
$$
with $\beta(x)=\alpha(x)\varphi_1(\rho_{m_1}(F_x)))$, which is also
supported on $S(X,\Omega;\Lb^*)$. By the induction hypothesis applied
first to $m_2$, then to $m_1$,
we obtain
\begin{align*}
\sum_{\varphi\in\bon_{m_1m_2}^*}{|S(\varphi)|^2}&\geq 
\prod_{\ell\mid m_2}
{\frac{\nu(\Omega_{\ell})}{\nu(Y_{\ell}-\Omega_{\ell})}}
\sum_{\varphi_1\in\bon_{m_1}^*}{
\Bigl|
\int_{X}{\beta(x)d\mu(x)}
\Bigr|^2
}\\
&=\prod_{\ell\mid m_2}
{\frac{\nu(\Omega_{\ell})}{\nu(Y_{\ell}-\Omega_{\ell})}}
\sum_{\varphi_1\in\bon_{m_1}^*}{|S(\varphi_1)|^2}\geq
\prod_{\ell\mid m_1m_2}
{\frac{\nu(\Omega_{\ell})}{\nu(Y_{\ell}-\Omega_{\ell})}}
\Bigl|\int_{X}{\alpha(x)d\mu(x)}\Bigr|^2. 
\end{align*}
\end{proof}

Now the proof of Proposition~\ref{pr-comb-sieve} is easy.

\begin{proof}[Proof of Proposition~\ref{pr-comb-sieve}]
Take $\alpha(x)$ to be the characteristic function of $S(X,\Omega;\Lb^*)$
and sum over $m\in \Lb$ the inequality of Lemma~\ref{lm-2b}; since
$$
\int_{X}{\alpha(x)d\mu(x)}=\int_{X}{\alpha(x)^2d\mu(x)}=|S(X,\Omega;\Lb^*)|,
$$
it follows that 
$$
|S(X,\Omega;\Lb^*)|^2\sum_{m\in \Lb}{
\prod_{\ell\mid m}
{\frac{\nu(\Omega_{\ell})}{\nu(Y_{\ell}-\Omega_{\ell})}}}
\leq \sum_{m\in \Lb}
{\sum_{\varphi\in \bon_m^*}{|S(\varphi)|^2}}\leq \Delta |S(X,\Omega;\Lb^*)|,
$$
hence the result.
\end{proof}

\begin{example}
In the classical case, with $Y=\Zz$ and $Y_{\ell}=\Zz/\ell\Zz$, we can
identity $Y_m$ with $\Zz/m\Zz$ by the Chinese Remainder Theorem. With
$\nu(y)=1/\ell$ for all $\ell$ and all $y$, the usual basis of functions on
$Y_m$ is that of additive characters
$$
x\mapsto e\Bigl(\frac{ax}{m}\Bigr)
$$
for $a\in \Zz/m\Zz$. It is easy to check that such a character belongs
to $\bon_m^*$ if and only if $a$ and $m$ are coprime.
\end{example}

At this point a ``large sieve inequality'' will be an
estimate for the quantity $\Delta$. There are various techniques
available for this purpose; see~\cite[Ch. VII]{ant} for a survey of
some of them. 
\par
The simplest  technique is to use the familiar
duality principle for bilinear forms or linear operators. Since
$\Delta$ is the square of the norm of a linear operator, it is the
square of the norm of its adjoint. Hence we have:

\begin{lemma}\label{lm-dual}
Let $\sieve=(Y,\Lambda,(\rho_{\ell}))$ be a sieve setting,
$(X,\mu,F)$ a siftable set, $\Lb$ a sieve support associated to $\Lb^*$. 
Fix orthonormal basis $\bon_{\ell}$ and define $\bon_m$ as above. 
Then 
the large sieve constant $\Delta(X,\Lb)$ is the smallest number
$\Delta$ such that
\begin{equation}\label{eq-dual-sieve}
\int_{X}{\Bigl|
\sum_{m\in \Lb}{\sum_{\varphi\in\bon_m^*}
{
\beta(m,\varphi)\varphi(\rho_m(F_x))}}
\Bigr|^2d\mu(x)}\leq \Delta\sum_{m}{\sum_{\varphi}{|\beta(m,\varphi)|^2}}
\end{equation}
for all vectors of complex numbers $(\beta(m,\varphi))$.
\end{lemma}

The point is that this leads to another bound for $\Delta$ in terms of
bounds for the ``dual'' sums
$W(\varphi,\varphi')$ obtained by 
expanding the square in this inequality, i.e.\label{pg-exp}
$$
W(\varphi,\varphi')=
\int_{X}{\varphi(\rho_m(F_x))
\overline{\varphi'(\rho_n(F_x))}d\mu(x)},
$$
where $\varphi\in\bon_m$ and $\varphi'\in\bon_n$ for some $m$ and $n$ in
$S(\Lambda)$. 
Precisely, we have:

\begin{proposition}\label{pr-w-to-sieve}
Let $\sieve=(Y,\Lambda,(\rho_{\ell}))$ be a sieve setting,
$\siftable=(X,\mu,F)$ a siftable set, $\Lb^*$ a prime sieve support
and $\Lb$ an associated sieve support. Then the
large sieve constant satisfies 
\begin{equation}\label{eq-w-to-sieve}
\Delta\leq \max_{n\in \Lb}\max_{\varphi\in\bon_n^*}{
\sum_{m\in \Lb}{\sum_{\varphi'\in\bon_m^*}{
|W(\varphi,\varphi')|}}}.
\end{equation}
\end{proposition}

\begin{proof}
Expanding the left-hand side of~(\ref{eq-dual-sieve}), we have
$$
\int_{X}{\Bigl|
\sum_{m\in \Lb}{\sum_{\varphi\in\bon_m^*}
{
\beta(m,\varphi)\varphi(\rho_m(F_x))
}}
\Bigr|^2d\mu(x)}=
\dblsum_{m,n}{\dblsum_{\varphi,\varphi'}
{\beta(m,\varphi)\overline{\beta(n,\varphi')}
W(\varphi,\varphi')}}
$$
and applying $|uv|\leq \demi(|u|^2+|v^2|)$ the result follows as
usual. 
\end{proof}

The point is that sieve results are now reduced to
individual \emph{uniform} estimates for the ``sums''
$W(\varphi,\varphi')$. Note that, here, the choice of the
orthonormal basis may well be very important in estimating
$W(\varphi,\varphi')$ and therefore $\Delta$.
\par
However, at least formally, we can proceed in full generality as
follows, where the idea is that in applications $\rho_m(F_x)$ should range
fairly equitably (with respect to the density $\nu_{m}$) over  the
elements of $Y_m$, so the sum  
$W(\varphi,\varphi')$ should be estimated by exploiting the
``periodicity'' of
$\varphi(\rho_m(F_x))\overline{\varphi'(\rho_n(F_x))}$. To do 
this, we introduce  further notation.
\par
Let $m$, $n$ be two elements of $S(\Lambda)$, $\varphi\in\bon_m$,
$\varphi'\in\bon_n$. Let $d=m\cap n$ be the intersection (g.c.d. in the case
of integers) of $m$ and $n$, and write
$m=m'd=m'\cup d$, $n=n'd=n'\cup d$ (disjoint unions). According to the
multiplicative definition of $\bon_m$ and $\bon_n$, we can write
$$
\varphi=\varphi_{m'}\otimes \varphi_d,\quad
\varphi'=\varphi'_{n'}\otimes \varphi'_d 
$$
for some unique basis elements $\varphi_{m'}\in\bon_{m'}$, $\varphi_d$,
$\varphi'_d\in \bon_d$ and $\varphi'_{n'}\in\bon_{n'}$.
\par
Let $[m,n]=mn=m\cup n$ be the ``l.c.m'' of $m$ and $n$. We have the
decomposition
$$
Y_{[m,n]}=Y_{m'}\times Y_d\times Y_{n'},
$$
the (not necessarily surjective) map $\rho_{[m,n]}\,:\, Y\ra Y_{[m,n]}$
and the function
\begin{equation}\label{eq-varphi-varphi}
[\varphi,\overline{\varphi'}]=\varphi_{m'}\otimes (\varphi_d
\overline{\varphi'_d})
\otimes \overline{\varphi'_{n'}}\,:\, (y_1,y_d,y_2)\mapsto
\varphi_{m'}(y_1)\varphi_d(y_d)\overline{\varphi'_d(y_d)}
\overline{\varphi'_{n'}(y_2)},
\end{equation}
(which is not usually a basis element in $\bon_{[m,n]}$). 
\par
The motivation for all this is the following tautology:

\begin{lemma}\label{lm-3}
Let $m$, $n$, $\varphi$ and $\varphi'$ be as before. We have
$$
[\varphi,\overline{\varphi'}](\rho_{[m,n]}(y))=\varphi(\rho_m(y))
\overline{\varphi'(\rho_n(y))}
$$
for all $y\in Y$, hence
$$
W(\varphi,\varphi')=
\int_{X}{[\varphi,\overline{\varphi'}](\rho_{[m,n]}(F_x))d\mu(x)}.
$$
\end{lemma}

Now we can hope to split the integral according to the value
of $y=\rho_{[m,n]}(F_x)$ in $Y_{[m,n]}$, and evaluate it by summing the
main term in an equidistribution statement.
\par
More precisely, for $d\in S(\Lambda)$ and $y\in Y_d$, we \emph{define}
$r_d(X;y)$ as the
``error term'' in the expected equidistribution statement:
\begin{equation}
\label{eq-rd}
|\{\rho_d(F_x)=y\}|
=\int_{\{\rho_d(F_x)=y\}}{d\mu(x)}=
\nu_d(y)|X|+r_d(X;y).
\end{equation}

Then we can write
$W(\varphi,\varphi')$  as described before:
\begin{align}
W(\varphi,\varphi')&=
\int_{X}{[\varphi,\overline{\varphi'}](\rho_{[m,n]}(F_x))d\mu(x)}
\nonumber\\
&=\sum_{y\in Y_{[m,n]}}
{
[\varphi,\overline{\varphi'}}](y)
\int_{\{\rho_{[m,n]}(F_x)=y\}}
{d\mu(x)}\nonumber\\
&=m([\varphi,\overline{\varphi'}])|X|
+O\Bigl(
\sum_{y\in Y_{[m,n]}}
{\|[\varphi,\overline{\varphi'}]\|_{\infty}|r_{[m,n]}(X;y)|}
\Bigr)\label{eq-periodique}
\end{align}
after inserting~(\ref{eq-rd}), where the implied constant is of
modulus $\leq 1$ and
$$
m([\varphi,\overline{\varphi'}])=
\sum_{y\in Y_{[m,n]}}{\nu_{[m,n]}(y)[\varphi,\overline{\varphi'}](y)}=
\langle [\varphi,\overline{\varphi'}],1\rangle,
$$ 
the inner product in $L^2(Y_{[m,n]})$.  One would then hope that
$m([\varphi,\varphi'])$ is the delta-symbol
$\delta((m,\varphi),(n,\varphi'))$ which would select the diagonal in
the main term of the sums $W(\varphi,\varphi')$.
In 
Sections~\ref{sec-spec-group} and~\ref{sec-spec-coset}, we will see
how to evaluate this quantity for the special case of group and coset
sieves. But first, a short digression...

\section{The ``dual'' sieve}\label{sec-dual}

The equivalent definition of the large sieve constant by means of the
duality principle (i.e, Lemma~\ref{lm-dual}) is quite useful 
in itself. For instance, it yields the 
following type of sieve inequality, which in the classical case goes
back to Rényi.

\begin{proposition}\label{pr-weaker-ls}
Let $(Y,\Lambda,(\rho_{\ell}))$  be a sieve setting, $(X,\mu,F)$ a
siftable set and $\Lb^*$ a prime sieve support. Let $\Delta$ be the 
large sieve constant for $\Lb=\Lb^*$.\footnote{\ Precisely, $\Lb$ is the
  set of singletons $\{\ell\}$ for $\ell\in\Lb^*$.} Then for any
sifting sets $(\Omega_{\ell})$, we have
\begin{equation}
\label{eq-weaker-ls}
\int_{X}{\Bigl(P(x,\Lb)-P(\Lb)\Bigr)^2d\mu(x)} \leq \Delta Q(\Lb)
\end{equation}
where
\begin{equation}\label{eq-pxl-pl}
P(x,\Lb)=\sum_{\stacksum{\ell\in\Lb}{\rho_{\ell}(F_x)\in
    \Omega_{\ell}}}
{1},\quad\quad
P(\Lb)=\sum_{\ell\in\Lb}{\nu(\Omega_{\ell})},\quad\quad
Q(\Lb)=\sum_{\ell\in\Lb}{\nu(\Omega_{\ell})(1-\nu(\Omega_{\ell}))}.
\end{equation} 
\end{proposition}

\begin{proof}
By expanding the characteristic 
function $\chi(\Omega_{\ell})$ of $\Omega_{\ell}\subset Y_{\ell}$ in
the orthonormal basis $\bon_{\ell}$, we obtain
$$
P(x,\Lb)=P(\Lb)+\sum_{\ell\in\Lb}{\sum_{\varphi\in\bon_{\ell}^*}{
\beta(\ell,\varphi)\varphi(\rho_{\ell}(F_x))}},
$$
where 
$$
\beta(\ell,\varphi)=\sum_{y\in \Omega_{\ell}}{
\nu_{\ell}(y)\overline{\varphi}(y)},
$$
and we used the fact that
$\bon_{\ell}^*=\bon_{\ell}-\{1\}$ for $\ell\in\Lambda$. Thus we get
\begin{align*}
\int_{X}{\Bigl(P(x,\Lb)-P(\Lb)\Bigr)^2d\mu(x)}
&=\int_{X}{\Bigl|
\sum_{\ell\in\Lb}{\sum_{\varphi\in\bon_{\ell}^*}{
\beta(\ell,\varphi)\varphi(\rho_{\ell}(F_x))}}
\Bigr|^2d\mu(x)}\nonumber\\
&\leq \Delta \sum_{\ell\in\Lb}{\sum_{\varphi\in\bon_{\ell}^*}
{|\beta(\ell,\varphi)|^2}}
\end{align*}
by applying~(\ref{eq-dual-sieve}). Since we have
$$
\sum_{\varphi\in\bon_{\ell}^*}{|\beta(\ell,\varphi)|^2}
=\sum_{\varphi\in\bon_{\ell}}{|\beta(\ell,\varphi)|^2}
-|\beta(\ell,1)|^2
=\|\chi(\Omega_{\ell})\|^2-\nu(\Omega_{\ell})^2=\nu(\Omega_{\ell})
(1-\nu(\Omega_{\ell})),
$$
this implies the result.
\end{proof}

In particular, since $P(x,\Lb)=0$ for $x\in
S(X;\Omega,\Lb^*)$ and $Q(\Lb)\leq P(\Lb)$, we get (by positivity
again) the estimate 
$$
|S(X;\Omega,\Lb^*)|\leq \Delta P(\Lb)^{-1},
$$
which is the analogue of the inequalities used e.g. by Gallagher
in~\cite[Th. A]{gallagher}, and by the author in~\cite{sieve}. This
inequality also follows from Proposition~\ref{pr-comb-sieve} if  we
take $\Lb$ containing only singletons (in the arithmetic case, this
means using only the primes), since we get the estimate
$$
|S(X,\Omega;\Lb^*)|\leq \Delta H^{-1}
\quad\quad\text{with}\quad\quad
H= \sum_{\ell\in \Lb}{\frac{\nu(\Omega_{\ell})}
{\nu(Y_{\ell}-\Omega_{\ell})}}\geq 
\sum_{\ell\in \Lb}{\nu(\Omega_{\ell})}=P(\Lb)
$$
(in fact, by Cauchy's inequality, we have $P(\Lb)^2\leq HQ(\Lb)$).
\par
This type of result is also related to
Tur\'an's method in probabilistic number theory. In counting primes
with the classical setting, or more generally in ``small sieve''
situations, it may seem quite weak (it 
only implies $\pi(X)\ll X(\log\log X)^{-1}$). However, it is really a
different type of statement, which has additional flexibility: for
instance, it still implies that for $X\geq 3$ we have
$$ 
|\{n\leq X\,\mid\, \omega(n)<\kappa\log\log X\}|
\ll \frac{X}{\log\log X}
$$
for any $\kappa\in ]0,1[$, the implied constant depending only on
$\kappa$. This estimate is now \emph{of the right order of magnitude},
and this shows in particular that one can not hope to
improve~(\ref{eq-weaker-ls}) by using information related to all
``squarefree'' numbers;  in other words,
Proposition~\ref{pr-comb-sieve} \emph{can not} be extended ``as is'' to
an upper bound for the variance on the left of~(\ref{eq-weaker-ls}).
\par
These remarks indicate that Proposition~\ref{pr-weaker-ls} has its own
interest in cases where the ``stronger'' form of the large sieve is
in fact not adapted to the type of situation considered.
In  Section~\ref{sec-elliptic-amusing}, we will describe an
amusing use of the inequality~(\ref{eq-weaker-ls}), where the ``pure
sieve'' bound would indeed be essentially trivial.

\section{Group and conjugacy sieves}\label{sec-spec-group} 

We now come to the description of a more specific type of sieve setting,
related to a group structure on $Y$. Together with the coset sieves of
Section~\ref{sec-spec-coset}, this exhausts most examples of
applications we know at the moment.
\par
A \emph{group sieve} corresponds to a sieve setting
$\sieve=(G,\Lambda,(\rho_{\ell}))$ where $G$ is a group and the maps
$\rho_{\ell}\,:\, G\ra G_{\ell}$ are homomorphisms onto finite
groups. A \emph{conjugacy sieve}, similarly, is a sieve setting
$\sieve=(G,\Lambda,(\rho_{\ell}))$ where $\rho_{\ell}\,:\, G\ra
G_{\ell}^{\sharp}$ is a surjective map from $G$ to the finite set of
conjugacy classes $G_{\ell}^{\sharp}$ of a finite group $G_{\ell}$,
that factors as 
$$
G\ra G_{\ell}\ra G_{\ell}^{\sharp}
$$
where $G\ra G_{\ell}$ is a surjective homomorphism.
Obviously, if $G$ is abelian, group and conjugacy sieves are
identical, and any group sieve induces a conjugacy sieve.
\par
The group structure suggests a natural choice of orthonormal basis
$\bon_{\ell}$ for functions on $G_{\ell}$ or $G_{\ell}^{\sharp}$, as
well as natural densities $\nu_{\ell}$. We start with the simpler
conjugacy sieve.
\par
From the classical representation theory of finite groups (see,
e.g.,~\cite{serre-rep}), we know that for any $\ell\in\Lambda$, the
functions 
$$
y\mapsto \Tr \pi(y),
$$
on $G_{\ell}$, where $\pi$ runs over the set $\Pi_{\ell}$ of
(isomorphism classes of) 
irreducible linear representations $\pi\,:\, G_{\ell}\ra GL(V_{\pi})$,
form an 
orthonormal basis of the space $\mathcal{C}(G_{\ell})$ of functions on
$G_{\ell}$ invariant under conjugation, with the inner product
$$
\langle f,g\rangle=\frac{1}{|G_{\ell}|}\sum_{y\in G_{\ell}}
{f(y)\overline{g(y)}}.
$$
\par
Translating this statement to functions on the set $G_{\ell}^{\sharp}$
of conjugacy classes, this
means that the functions 
$$
\varphi(\yc)=\Tr\pi(\yc)
$$
on  $G_{\ell}^{\sharp}$  form an
orthonormal basis $\bon_{\ell}$ of $L^2(G_{\ell}^{\sharp})$ with the
inner product 
$$
\langle f,g\rangle=\frac{1}{|G_{\ell}|}\sum_{\yc\in G}
{|\yc|f(\yc)\overline{g(\yc)}}.
$$
\par
Moreover, the trivial representation $1$ of $G_{\ell}$ has for
character the constant function $1$, so we can use the basis
$\bon_{\ell}=(\Tr\pi(\yc))_{\pi}$ for computing the large sieve
constant if the density
$$
\nu_{\ell}(\yc)=\frac{|\yc|}{|G_{\ell}|}
$$
is used. Note that this is the image on $G_{\ell}^{\sharp}$ of the
uniform density on $G_{\ell}$.
\par
Note also that in the abelian case, the representations are
one-dimensional, and the basis thus described is the basis of
characters of $G_{\ell}$, with  the uniform density, i.e., that of
group homomorphisms $G_{\ell}\ra \Cc^{\times}$ with
$$
\langle f,g\rangle=\frac{1}{|G_{\ell}|}\sum_{y\in G_{\ell}}
{f(y)\overline{g(y)}}.
$$
\par
Coming back to a general group sieve, the basis and densities extended
to the sets 
$$
G_m^{\sharp}=\prod_{\ell\mid m}{G_{\ell}^{\sharp}}
$$
for $m\in S(\Lambda)$ have a similar interpretation. Indeed,
$G_m^{\sharp}$ identifies clearly with the set of conjugacy classes of
the finite group  $G_m=\prod G_{\ell}$. The density $\nu_m$ is
therefore still given by
$$
\nu_m(\yc)=\frac{|\yc|}{|G_m|}.
$$
\par
Also, it it well-known that the
irreducible representations of $G_m$ are of the form
$$
\pi\,:\, g\mapsto \boxtimes_{\ell\mid m}{\pi_{\ell}(g)}
$$
for some uniquely defined irreducible  representations $\pi_{\ell}$ of
$G_{\ell}$, where $\boxtimes$ is the external tensor product defined
by 
$$
g=(g_{\ell})\mapsto \bigotimes_{\ell\mid m}{\rho_{\ell}(g_{\ell})}.
$$
In other words, the set $\Pi_m$ of irreducible linear representations
of $G_m$ is identified canonically with $\prod{\Pi_{\ell}}$. Moreover,
the character of  a representation of $G_m$ of this form is simply 
$$
\Tr\pi(g)=\prod_{\ell\mid m}{\Tr\pi_{\ell}(g_{\ell})},
$$
so that the basis $\bon_m$ obtained from $\bon_{\ell}$ is none other
than the basis of functions $\yc\mapsto \Tr\pi(\yc)$ for $\pi$
ranging over $\Pi_m$.
\par
Given a siftable set $(X,\mu,F)$ associated to a conjugacy sieve
$(G,\Lambda,(\rho_{\ell}))$, the sums $W(\varphi,\varphi')$ become 
\begin{equation}\label{eq-sumexp}
W(\pi,\tau)=\int_{X}{\Tr\pi(\rho_m(F_x))
\overline{\Tr\tau(\rho_n(F_x))}d\mu(x)}
\end{equation}
for irreducible representations $\pi$ and $\tau$ of $G_m$ and $G_n$
respectively, 
which can usually be interpreted as \emph{exponential sums} (or integrals)
over $X$, since the character values, as traces of matrices of finite
order, are sums of finitely many roots of unity.
\par
We summarize briefly by reproducing the general sieve results in this
context, phrasing things as related to the conjugacy sieve induced
from a group sieve (which seems most natural for applications). In
such a situation, the sieving sets $\Omega_{\ell}$ are naturally given
as conjugacy-invariant subsets of $G_{\ell}$, and are identified with
subsets $\Omega_{\ell}^{\sharp}$ of $G_{\ell}^{\sharp}$. Note that we
have then  
$\nu_m(\Omega_{\ell}^{\sharp})=|\Omega_{\ell}|/|G_{\ell}|$.

\begin{proposition}\label{pr-conjug-sieve}
Let $(G,\Lambda,(\rho_{\ell}))$ be a group sieve setting,
$(X,\mu,F)$ an associated siftable set. For any prime sieve support
$\Lb^*$ and an associated sieve support $\Lb$
satisfying~\emph{(\ref{eq-sieve-support})}, and for any conjugacy
invariant sifting sets $(\Omega_{\ell})$, we have
$$
|S(X,\Omega;\Lb^*)|\leq \Delta H^{-1}
$$
where $\Delta$ is the smallest non-negative real number such that
$$
\sum_{m\in\Lb}{\sum_{\pi\in\Pi_m^*}{\Bigl|
\int_X{\alpha(x)\Tr\pi(\rho_m(F_x))d\mu(x)}
\Bigr|^2}}\leq \Delta\int_X{|\alpha(x)|^2d\mu(x)}
$$
for all square-integrable function $\alpha\in L^2(X,\mu)$, where $\pi$
ranges over the set $\Pi_m^*$ of \emph{primitive} irreducible linear
representations of $G_m$, i.e., those such that 
no component $\pi_{\ell}$ for $\ell\mid m$ is trivial, and where
$$
H=\sum_{m\in\Lb}{\prod_{\ell\mid m}{\frac{|\Omega_{\ell}|}
{|G_{\ell}|-|\Omega_{\ell}|}}}.
$$
\par
Moreover we have
$$
\Delta\leq
\max_{m\in\Lb}\max_{\pi\in\Pi_m^*}\sum_{n\in\Lb}\sum_{\tau\in\Pi_n^*}{ 
|W(\pi,\tau)|},
$$
where
$$
W(\pi,\tau)=\int_{X}{\Tr\pi(\rho_m(F_x))
\overline{\Tr\tau(\rho_n(F_x))}d\mu(x)}.
$$
\end{proposition}


The general sieve setting can also be applied to problems where the
sieving sets are not conjugacy-invariant, using the basis of
\emph{matrix coefficients} of irreducible representations.
Let $(G,\Lambda,(\rho_{\ell}))$ be a group sieve setting. For each
$\ell$ and each irreducible representation $\pi\in\Pi_{\ell}$, choose an
orthonormal basis $(e_{\pi,i})$ of the space $V_{\pi}$ of the
representation (with respect to a $G_{\ell}$-invariant inner
product $\langle\cdot,\cdot\rangle_{\pi}$). Then (see, e.g.,~\cite[\S
I.5]{knapp}, which treats compact 
groups), the family $\bon_{\ell}$ of functions of the type
$$
\varphi_{\pi,e,f}\,:\, x\mapsto \sqrt{\dim \pi}
\langle \pi(x)e, f\rangle_{\pi},\quad\quad e=e_{\pi,1},\ldots,
e_{\pi,\ell},\quad
f=e_{\pi,1},\ldots, e_{\pi,\ell}
$$
is an orthonormal basis of $L^2(G_{\ell})$ for the inner product
$$
\langle f,g\rangle=\frac{1}{|G_{\ell}|}\sum_{x\in G_{\ell}}{f(x)
\overline{g(x)}},
$$
i.e., corresponding to the density $\nu_{\ell}(x)=1/|G_{\ell}|$ for
all $x\in G_{\ell}$.
Moreover, for $\pi=1$, and an arbitrary choice of $e\in\Cc$ with
$|e|=1$, the function $\varphi_{1,e,e}\in \bon_{\ell}$ is the
constant function $1$.
\par
If we extend the basis $\bon_{\ell}$ to orthonormal basis $\bon_m$ of
$L^2(G_m)$ for all $m\in S(\Lambda)$, by multiplicativity, the
functions in $\bon_m$ are of the type
$$
\varphi_{\pi,e,f}\,:\,
(x_{\ell})\mapsto \sqrt{\dim \pi}\prod_{\ell\mid m}{\langle
\pi_{\ell}(x_{\ell})e_{\ell},
f_{\ell}
\rangle_{\pi_{\ell}}}
$$
where $e=\otimes e_{\ell}$ and $f=\otimes f_{\ell}$ run
over elements of the orthonormal basis
$$
\Bigl(\bigotimes_{\ell\mid m}{e_{\pi_{\ell},i_{\ell}}}\Bigr),\quad\quad
1\leq i_{\ell}\leq \dim \pi_{\ell},
$$
constructed from the chosen bases $(e_{\pi,i})$ of the
components, the inner product on the space of $\boxtimes\pi_{\ell}$ being
the natural $G_{m}$-invariant one.
\par
The sums $W(\varphi,\varphi')$ occurring in
Proposition~\ref{pr-w-to-sieve} to estimate the large sieve constant
are given by
\begin{equation}\label{eq-matrix-coeffs}
W(\varphi_{\pi,e,f},\varphi_{\tau,e',f'})=
\sqrt{(\dim\pi)(\dim\tau)}
\int_X{
\langle \pi(\rho_m(F_x))e,f\rangle_{\pi}
\overline{
\langle \tau(\rho_n(F_x))e',f'\rangle_{\tau}
}
d\mu(x)}.
\end{equation}
\par
If we apply Lemma~\ref{lm-3} to  elements $\varphi_{\pi,e,f}$,
$\varphi_{\tau,e',f'}$ of the basis $\bon_m$ and $\bon_n$ of
$L^2(G_m)$, 
the function $[\varphi_{\pi},\overline{\varphi_{\tau}}]$ which is
integrated can be written as a matrix coefficient of
the representation
\begin{equation}\label{eq-lcm-reps}
[\pi,\bar{\tau}]=\pi_{m'}\boxtimes
(\pi_d\otimes\overline{\tau_d})\boxtimes \overline{\tau_{n'}}
\end{equation}
of $G_{[m,n]}$, where we write $\pi=\pi_{m'}\boxtimes \pi_d$,
$\tau=\tau_{n'}\boxtimes\tau_d$, with the obvious meaning of the
components $\pi_{m'}$, $\pi_d$, $\tau_d$, $\tau_{n'}$, and the bar
indicates taking the contragredient representation. 
\par
Indeed, we have
$$
[\varphi_{\pi,e,f},\overline{\varphi_{\tau,e',f'}}](x_{\ell})=
\sqrt{(\dim\pi)(\dim\tau)}\langle
[\pi,\bar{\tau}](\rho_{[m,n]}(F_x))\tilde{e},\tilde{f}
\rangle_{[\pi,\bar{\tau}]}
$$
for $(x_{\ell})\in G_{[m,n]}$, with $\tilde{e}=e\otimes e'$,
$\tilde{f}=f\otimes f'$.
\par
Concretely, this means that in order to deal with the sums
$W(\varphi,\varphi')$ 
to estimate the large sieve constant using the basis $\bon_m$ of
matrix coefficients, it suffices to be able to estimate all integrals
of the type
\begin{equation}\label{eq-int-varpi}
\int_{X}{\langle \varpi(F_x)e,f\rangle_{\varpi} d\mu(x)}.
\end{equation}
where $\varpi$   is a representation of $G$ that factors through a
finite product of groups $G_{\ell}$, and $e$, $f$ are vectors in the
space of the representation $\varpi$ (the inner product being
$G$-invariant). 
See the proof of Theorem~\ref{th-non-conj} for an application of this.

\begin{remark}
Another potentially useful sieve setting associated to a group sieve
setting $(G,\Lambda,\rho_{\ell})$ is obtained by replacing
$\rho_{\ell}$ with the projections $G\ra G_{\ell}\ra
G_{\ell}/K_{\ell}=Y_{\ell}$ for $\ell\in\Lambda$, where $K_{\ell}$ is
an arbitrary subgroup of $G_{\ell}$. Considering the density on
$Y_{\ell}$ which is image of the uniform density on $G_{\ell}$, an
orthonormal basis $\bon_{\ell}$ of $L^2(Y_{\ell})$ is then obtained by
taking the functions
$$
\varphi_{\pi,e,f}\,:\, gK_{\ell}\mapsto \langle \pi(g)e,f\rangle
$$
where $\pi$ runs over irreducible representations of $G_{\ell}$,  $e$
runs over an orthonormal basis of the $K_{\ell}$-invariant subspace in
the space $V_{\pi}$ of $\pi$, and $f$ over a full orthonormal basis of
$V_{\pi}$.
\par
Indeed, the restriction on $e$ ensures that such functions are well-defined
on $G_{\ell}/K_{\ell}$ (i.e., the matrix coefficient is
$K_{\ell}$-invariant), and since those are matrix coefficients, there
only remains to check that they span $L^2(Y_{\ell})$. However, the
total number of functions is
$$
\sum_{\pi}{(\dim\pi)
\langle \res^{G_{\ell}}_{K_{\ell}}\pi,1\rangle_{K_{\ell}}}=
\sum_{\pi}{(\dim\pi)
\langle \pi,\Ind_{K_{\ell}}^{G_{\ell}}1\rangle_{G_{\ell}}}=
\dim \Ind_{K_{\ell}}^{G_{\ell}}1=|Y_{\ell}|
$$
and since they are independent, the result follows.
\par
Because this basis is a sub-basis of the previous one, any estimate
for the large sieve constant for the group sieve will give one for
this sieve setting.
\end{remark}


\section{Elementary and classical examples}\label{sec-examples}

We first describe how some classical uses of the large sieve are special
cases of the group sieve setting of the previous section, and conclude
this section with a ``new'' example of the general case which is
particularly easy to analyze (and of little practical use), and hence
somewhat enlightening.

\begin{example}
As already mentioned, the classical
large sieve arises from the group sieve setting
$$
\sieve=(\Zz,\{\text{primes}\},\Zz\ra \Zz/\ell\Zz)
$$
where the condition for an additive character $x\mapsto e(ax/m)$ of
$G_{m}=(\Zz/m\Zz)$ to be primitive is equivalent with the classical
condition that $(a,m)=1$.
\par
In the most typical case, the siftable sets are
$$
X=\{n\geq 1\,\mid\, N\leq n<N+M\}
$$
with $F_x=x$, 
and the abstract sieving problem becomes the ``original'' one of finding
integers in $X$ which lie outside certain residue classes modulo some
primes $\ell$.
\par
More generally, take 
$$
\sieve=(\Zz^r, \{\text{primes}\},\Zz^r\ra
(\Zz/\ell\Zz)^r)
$$
(the reduction maps) and $X=\{(a_1,\ldots, a_r)\,\mid\, N_i\leq
a_i<N_i+M_i\}$, with $F$ the identity map again. Then what results
is the higher-dimensional large sieve (see e.g.~\cite{gallagher}).
\par
For completeness, we recall the estimates available for the
large sieve constant in those two situations, when we take $\Lb^*$ to
be the set of primes $\leq L$, and $\Lb$ to
be the set of squarefree integers $\leq L$, for some $L\geq 1$. We
write $S(X;\Omega,L)$ instead of $S(X;\Omega,\Lb^*)$.
\begin{theorem}
With notation as above, we have $\Delta\leq N-1+L^2$ for $r=1$ and
$\Delta\leq  (\sqrt{N}+L)^{2r}$ for all $r\geq 1$. In particular, for
any sieve problem, we have
\begin{gather*}
|S(X;\Omega,L)|\leq (N-1+L^2)H^{-1},\quad\quad\text{ if } r=1,\\
|S(X;\Omega,L)|\leq (\sqrt{N}+L)^{2r}H^{-1},\quad\quad\text{ if } r\geq 1,
\end{gather*}
where
$$
H=\sumb_{m\leq L}{\prod_{\ell\mid m}
{\frac{|\Omega_{\ell}|}{\ell^r-|\Omega_{\ell}|}}},
$$
the notation $\sumb$ indicating a sum restricted to squarefree numbers.
\end{theorem}

In the one-variable case, this is due essentially to Montgomery, and
to Selberg with the constant $N-1+L^2$, see e.g.~\cite[\S 7.5]{ant};
the  higher-dimensional case as stated is due to Huxley,
see~\cite{huxley}. Note that modern treatments deduce such estimates
from an analytic inequality which is more general than the ones
we used, namely, the inequality
$$
\sum_{r}{\Bigl|\sum_{M<n\leq M+N}{a_n e(n\xi_r)}\Bigr|^2}\leq
(N-1+\delta^{-1})\sum_{n}{|a_n|^2}
$$
for arbitrary sets $(\xi_r)$ of elements in $\Rr/\Zz$ which are
$\delta$-spaced, i.e., the distance $d(\xi_r,\xi_s)$ in $\Rr/\Zz$ is
at least $\delta$ if $r\not=s$ 
(this was first considered by Bombieri and Davenport; see,
e.g.,~\cite[Th. 7.7]{ant}). This amounts, roughly, to considering sums 
$$
\sum_{M<n\leq M+N}{e((\xi_r-\xi_s)n)}=W(\pi_r,\pi_s)
$$
where $\pi_r\,:\, n\mapsto e(n\xi_r)$ and $\pi_s$ are representations
of $G=\Zz$ which \emph{do not} factor through a finite index
subgroup. This suggests trying to prove similar inequalities for
general groups sieves, i.e., essentially, consider
integrals~(\ref{eq-int-varpi}) for arbitrary (unitary) representations
$\varpi$ of $G$.  
\par
Note that for $r=1$, the equidistribution assumption~(\ref{eq-rd})
becomes
$$
\sum_{\stacksum{N\leq n<N+M}{n\equiv y\mods{d}}}{1}=
\frac{M}{d}+r_d(X;y),
$$
which holds with $|r_d(X;y)|\leq 1$ for any  $y\in \Zz/d\Zz$. 
From~(\ref{eq-ex-statement}) we obtain the estimate $\Delta\leq
N+L^4$, which is by no means ridiculous. (See 
Section~\ref{sec-spec-coset} for the computation of the quantity
$m([\varphi,\varphi'])$ for group sieves, or do the exercise).
\par
Classical sieve theory is founded on such assumptions
as~(\ref{eq-rd}), usually stated merely for $y=0$, and on further
assumptions  concerning the resulting level of distribution, i.e.,
bounds for 
$r_d(X;0)$ on average over $d$ in a range as large as possible
(compared with the size of $X$). More general bounds for $r_d(X;y)$ do
occur however.
\par
Note that, even if this is classical, the general framework clearly
shows that to sieve an arbitrary set of integers
$X\subset \{n\,\mid\, n\geq 1\}\subset \Zz$, it suffices (at least up
to a point!) to have estimates for exponential sums
$$
\sum_{x\in X}{e\Bigl(\frac{ax}{m}-\frac{bx}{n}\Bigr)}
$$
with $n$, $m$ squarefree and $(a,m)=(b,n)=1$. It suffices, in
particular, to have equidistribution of $X$ in (all) arithmetic
progressions. This means for instance that some measure of large sieve
is usually doable for any sequence for which the classical ``small''
sieves work. This is of particular interest if $X$ is ``sparse'', in
the sense that e.g. $X\subset \{n\,\mid\, N<n\leq 2N\}$ for some $N$
with $|X|/N$ going to zero.
\par
It would also be interesting, as a problem in itself, to investigate
the values of the large sieve constant when using other sieve support
than squarefree integers up to $L$, for instance when the sieve
support is the support of a combinatorial (small) sieve.
\end{example}

\begin{example}
Can the multiplicative large sieve inequality for Dirichlet characters
be related to our general setting? Indeed, in at least two ways. 
First, let $q\geq 2$ be
given, let $G$ be the multiplicative subgroup of
$\Qq^{\times}$ generated by primes $p>q$, and take 
$$
\sieve=(G,\{\text{primes $\ell>q$}\},G\ra
(\Zz/\ell\Zz)^{\times}=G_{\ell}).
$$
In that context, we can take 
$$
X=\{n\leq N\,\mid\, p\mid n\Rightarrow p>q\},
$$
and $F_x=x$, and if $\Lb^*$ is the set of primes $\leq L\leq q$, and
$\Lb$ is the set of squarefree numbers $\leq L$, the sifted sets
become 
$$
S(X;\Omega,L)=\{ n\leq N\,\mid\, (p\mid n\Rightarrow p>q)\text{ and }
n\mods{\ell}\notin \Omega_{\ell}\text{ for } \ell\leq L\leq q\},
$$
where $\Omega_{\ell}\subset (\Zz/\ell\Zz)^{\times}$. A simple check
shows that the inequality defining the large sieve constant $\Delta$
becomes
\begin{equation}\label{eq-mult-ls}
\sumb_{m\leq L}{\sums_{\chi\mods{m}}{\Bigl|
\sum_{n\in X}{a_n \chi(n) \Bigr|^2}}}\leq \Delta
\sum_{n\in X}{|a_n|^2}
\end{equation}
for any complex numbers $a_n$, where $\chi$ runs over primitive
characters modulo $m$, and hence 
$\Delta\leq N-1+L^2$ by the multiplicative large sieve inequality (see
e.g.~\cite[Th. 7.13]{ant}).  
\par
Alternately, if we allow the density $\nu_{\ell}$ to have zeros, we may
take the classical sieve setting $Y=\Zz$, $Y\ra Y_{\ell}=\Zz/\ell\Zz$,
$X=\{N<n\leq M+N\}$, $F_x=x$, with density
$$
\nu_{\ell}(y)=\frac{1}{\ell-1}\quad\text{if } y\not=0,\quad\quad
\nu_{\ell}(0)=0,
$$
and then check that since the final statements do not involve the
inverse of $\nu_{\ell}(y)$, although the proofs involved division by
$\nu_{\ell}(y)$, it remains true that for $\Omega_{\ell}\subset
\Zz_{\ell}-\{0\}$, we have
$$
|S(X,\Omega;L)|\leq \Delta H^{-1}
$$
where 
$$
H=\sumb_{m\leq L}{\prod_{\ell\mid m}
{\frac{|\Omega_{\ell}|}{\ell-1-|\Omega_{\ell}|}}}
$$
and $\Delta$ is the multiplicative large sieve constant defined
by~(\ref{eq-mult-ls}) (e.g., use a positive perturbation
$\nu_{\ell,\eps}(y)>0$ of the density so that $H_{\eps}\ra H$ and
$\Delta_{\eps}\ra \Delta$, as $\eps\ra 0$). Again, we have
$\Delta\leq N-1+L^2$.

\end{example}

\begin{example}
Serre~\cite{serre-brauer} has used a variant of the higher-dimensional
large sieve where 
$$
\sieve=(\Zz^r,\{\text{primes}\},\Zz^r\ra (\Zz/\ell^2\Zz)^r)
$$
and
$$
X=\{(x_1,\ldots ,x_r)\in \Zz^r\,\mid\, |x_i|\leq
N\}
$$
with $F_x=x$. With suitable sifting sets, this provides estimates for
the number of trivial specializations of elements of $2$-torsion in
the Brauer group of  $\Qq(T_1,\ldots,T_r)$. 
\end{example}

\begin{example}
Here is a new example, which is a number field analogue of the
situation of~\cite{sieve} (described also in the
Section~\ref{sec-ls-ff}). It is related to Serre's discussion
in~\cite{serre-higher} of a higher-dimensional Chebotarev density
theorem over number fields (see also~\cite{pink} for an independent
treatment with more details). Let $Y/\Zz$ be a separated scheme of
finite type, and let $Y_{\ell}\ra Y$ be a family of
étale Galois
coverings,\footnote{\ Or better with ``controlled'' ramification, if not
  étale, since this is likely to be needed for some natural
  applications.} 
corresponding to surjective maps $G=\pi_1(Y,\bar{\eta})\ra
G_{\ell}$. The sieve setting is $(G,\{\text{primes}\},G\ra G_{\ell})$.
Now let $|Y|$ denote the set of
closed points of $Y$, which means those where the residue field $k(y)$
is finite, and let
$$
X=\{y\in |Y|\,\mid\, |k(y)|\leq T\}
$$
for some $T\geq 2$, which is finite. For $y\in X$, denote by $F_x\in
G$ the corresponding geometric Frobenius automorphism (or conjugacy
class rather) to obtain a
siftable set $(X,\text{counting measure},F)$ associated with the
conjugacy sieve.
It should be possible
to obtain a large sieve inequality in this context, at least assuming
GRH and the Artin conjecture. 
\par
Note that if $Y$ is the spectrum of the ring of integers in some
number field (or even $Y=\spec(\Zz)$ itself), this becomes the sieve
for Frobenius considered by Zywina~\cite{zywina}, with applications
(under GRH) to
the Lang-Trotter conjecture, and to Koblitz's conjecture for
elliptic curves over number fields.
\end{example}

\begin{example}\label{ex-excl-incl}
The next example illustrates the general sieve setting, showing  that
it includes (and extends) the inclusion-exclusion familiar in
combinatorics and  probability theory, and also that the large sieve
inequality is sharp in this general context (i.e., there may be
equality $|S(X,\Omega;\Lb^*)|=\Delta H^{-1}$).
\par
Let $(\Omega,\Sigma,\proba)$ be a probability space and $A_{\ell}\subset
\Sigma$, for $\ell\in\Lambda$, a countable family  of events. Consider
the event 
$$
A=\{\omega\in\Omega\,\mid\, \omega\notin A_{\ell}\text{ for any
}\ell\in\Lambda\}.
$$
\par
For $m\in S(\Lambda)$, denote 
$$
A_{m}=\bigcap_{\ell\in m}{A_{\ell}},\quad\quad A_{\emptyset}=\Omega.
$$
\par
If $\Lambda$ is finite, which we now assume, the inclusion-exclusion
formula is 
$$
\proba(A)=\sum_{m\in S(\Lambda)}{(-1)^{|m|}\proba(A_m)},
$$
and in particular, if the events are independent (as a whole), we have
$$
\proba(A_m)=\prod_{\ell\in m}{\proba(A_{\ell})},\quad\quad\text{ and }
\quad\quad 
\proba(A)=\prod_{\ell\in\Lambda}{(1-\proba(A_{\ell}))}.
$$
\par
Take the sieve setting $(\Omega,\Lambda,\charfun_{A_{\ell}})$, where
$\charfun_{B}$ is the characteristic function of an event $B$, with
$Y_{\ell}=\{0,1\}$ for all $\ell$, and the siftable set
$(\Omega,\proba,\mathrm{Id})$. 
Choose the density $\nu_{\ell}=\charfun_{A_{\ell}}(P)$, i.e., put
$$
\nu_{\ell}(1)=\proba(A_{\ell}),\quad\quad \nu_{\ell}(0)=1-\proba(A_{\ell}).
$$
With sieving sets $\Omega_{\ell}=\{1\}$ for $\ell\in\Lambda$, we
have precisely $S(X,\Omega;\Lambda)=A$.
\par
The large sieve inequality yields
$$
\proba(A)\leq \Delta H^{-1}
$$
where 
$$
H=\sum_{m\in\Lb}{\prod_{\ell\in m}{
\frac{\proba(A_{\ell})}{1-\proba(A_{\ell})}}},
$$
and $\Delta$ is the large sieve constant for the sieve support $\Lb$,
which may be any collection of subsets of $\Lambda$ such that
$\{\ell\}\in \Lb$ for all $\ell\in\Lambda$.
\par
Coming to the large sieve constant, note that $L^2_0(Y_{\ell})$ is
one-dimensional for all $\ell$, hence so is $L^2_0(Y_m)$ for all $m$
(including $m=\emptyset$). The basis function $\varphi_{\ell}$ for
$L^2_0(Y_{\ell})$ (up to multiplication by a complex number with
modulus $1$) is given  by
$$
\varphi_{\ell}(y)=\frac{y-p_{\ell}}{\sqrt{p_{\ell}(1-p_{\ell})}}
$$
where $p_{\ell}=\proba(A_{\ell})$ for simplicity, so that
$$
\varphi_{\ell}(\charfun_{A_{\ell}})=
\frac{\charfun_{A_{\ell}}-\proba(A_{\ell})}
{\sqrt{\variance(\charfun_{A_{\ell}})}},
$$
and in particular 
$$
\expect(\varphi_{\ell}(\charfun_{A_{\ell}}))=\langle
\varphi_{\ell},1\rangle=0
,\quad\quad
\expect(\varphi_{\ell}(\charfun_{A_{\ell}})^2)=\|\varphi_{\ell}\|^2=1.
$$
\par
Hence, for $\ell$, $\ell'\in \Lambda$,
$W(\varphi_{\ell},\varphi_{\ell'})$ is given by
$$
W(\varphi_{\ell},\varphi_{\ell'})=
\expect(\varphi_{\ell}(\charfun_{A_{\ell}})\varphi_{\ell'}
(\charfun_{A_{\ell'}}))
$$
and it is (by definition) the correlation coefficient of the
random variables $\charfun_{A_{\ell}}$ and $\charfun_{A_{\ell'}}$,
explicitly
$$
W(\varphi_{\ell},\varphi_{\ell'})=
\begin{cases}
1&\text{if } \ell=\ell',\\
\displaystyle{\frac{
\proba(A_{\ell}\cap A_{\ell'})-\proba(A_{\ell})\proba(A_{\ell'})}
{\sqrt{p_{\ell}(1-p_{\ell})p_{\ell'}(1-p_{\ell'})}}
}
&\text{ otherwise.}
\end{cases}
$$
\par
If (and only if) the $(A_{\ell})$ form a family of pairwise
independent events, 
we see that 
$W(\varphi_{\ell},\varphi_{\ell'})=\delta(\ell,\ell')$. More
generally, in all cases, for any $m$, $n\subset \Lambda$, we have
$$
W(\varphi_m,\varphi_n)=
\expect\Bigl(\prod_{\ell\in m}{\varphi_{\ell}(\charfun_{A_{\ell}})}
\prod_{\ell\in n}{\varphi_{\ell}(\charfun_{A_{\ell}})}\Bigr)
$$
which is a multiple normalized centered moment of the
$\charfun_{A_{\ell}}$.  
\par
If the $(A_{\ell})$ are globally independent, we obtain
\begin{align*}
W(\varphi_m,\varphi_n)&=
\prod_{\stacksum{\ell\in m\cup n}{\ell\notin m\cap n}}
\frac{\expect(\charfun_{A_{\ell}}-p_{\ell})}
{\sqrt{\variance(\charfun_{A_{\ell}})}}
\prod_{\ell\in m\cap n}{
\frac{\expect((\charfun_{A_{\ell}}-p_{\ell})^2)}
{
{\sqrt{\variance(\charfun_{A_{\ell}})}}
}}
\\
&=\delta(m,n)
\end{align*}
(since the third factor vanishes if the product is not empty, i.e., if
$m\not=n$, and the third term is
$1$ by orthonormality of $\varphi_{\ell}$).
It follows by~(\ref{eq-w-to-sieve}) that 
$\Delta\leq 1$, and in fact there must be equality.
Moreover, in this situation, if $\Lb$ contains all subsets of
$\Lambda$, we have 
$$
H=\prod_{\ell\in\Lambda}{\Bigl(1+\frac{p_{\ell}}{1-p_{\ell}}\Bigr)}
=\prod_{\ell\in\Lambda}{
\frac{1}{1-p_{\ell}}
},
$$
so that we find
$$
\Delta H^{-1}\leq \prod_{\ell\in\Lambda}{(1-\proba(A_{\ell}))}=\proba(A),
$$
i.e., the large sieve inequality is an equality here.
\par
Similarly, the inequality~(\ref{eq-weaker-ls}) becomes an equality if
the events are pairwise independent, and reflects the formula for the
variance of a sum of (pairwise) independent random variables.
\par
In the general case of possibly dependent events, on the other hand,
we have a quantitative inequality for $\proba(A)$ which may be of some
interest (and may be already known!). In fact, we have several
possibilities depending on the choice of sieve support. It would be
interesting to determine if those inequalities are of some use in
probability theory.


To conclude this example, note that  \emph{any} sieve, once the prime
sieve support $\Lb^*$ and the sieving sets
$(\Omega_{\ell})$ are chosen, may be considered as  a similar ``binary''
sieve with $Y_{\ell}=\{0,1\}$ for all $\ell$, by replacing the sieve
setting $(Y,\Lambda,(\rho_{\ell}))$  with
$(Y,\Lb^*, \charfun_{\Omega_{\ell}})$.
\end{example}

\begin{example}
There are a few examples of the use of simple sieve methods in
combinatorics. An example is  a paper of Liu and
Murty~\cite{liu-murty} (mentioned to us by A. Granville), which
explores (with some interesting 
combinatorial applications) a simple form of the dual sieve. Their
sieve setting amounts to taking $\sieve=(A,B,\charfun_b)$ 
where $A$ and $B$ are finite sets, and for each $b\in B$, we have a map
$\charfun_b\,:\, A\ra \{0,1\}$ (in loc. cit., the authors see $(A,B)$
as a bipartite graph, and $\charfun_b(a)=1$ if and only if there is an
edge from $a$ to $b$); the siftable set is $A$ with identity map and
counting measure, and the density is determined by
$\nu_b(1)=|\charfun_b^{-1}(1)|/|A|$. In other words, this is also a
special case of the previous example, and Theorem~1 and Corollary~1 of
loc. cit. can also be trivially deduced from this (though they are
simple enough to be better considered separately).
\end{example}

\section{Coset sieves}\label{sec-spec-coset}

Our next subject is a generalization of group sieves, which is the
setting in which the Frobenius sieve over finite fields
of~\cite{sieve} and Section~\ref{sec-ls-ff} operates. 
\par
As in Section~\ref{sec-spec-group}, we start with a group $G$ and a
family of surjective homomorphisms $G\ra G_{\ell}$, for $\ell\in\Lambda$,
onto finite groups. 
However, we also assume that there is a normal subgroup $\ggeom$ of
$G$ such that the quotient $G/\ggeom$ is abelian, and hence
we obtain a commutative diagram with exact rows
\begin{equation}\label{eq-cd}
\begin{CD}
1 @>>> \ggeom @>>> G
@>d>> G/\ggeom @>>>  1\\
@. @VVV   @V\rho_{\ell} VV  @Vp VV \\
1 @>>> \ggeom_{\ell} @>>>  G_{\ell}
@>d>> \Gamma_{\ell} @>>>  1,
\end{CD}
\end{equation}
where the downward arrows are surjective
and the quotient groups $\Gamma_{\ell}$ thus defined are finite abelian
groups. 
\par
After extending the definition of $G_m$ to elements of $S(\Lambda)$ by
multiplicativity, we can also define
$$
\ggeom_m=\prod_{\ell\mid m}{\ggeom_{\ell}}
$$
and we still can write commutative diagrams with exact rows
\begin{equation}\label{eq-cd-sqf}
\begin{CD}
1 @>>> \ggeom @>>> G
@>d>> G/\ggeom @>>>  1\\
@. @VVV   @V\rho_m VV  @Vp VV \\
1 @>>> \ggeom_{m} @>>>  G_{m}
@>d>> \Gamma_{m} @>>>  1,
\end{CD}
\end{equation}
(but the downward arrows are no longer necessarily surjective).
\par
The sieve setting for a coset sieve is then $(Y,\Lambda,(\rho_{\ell}))$
where $Y$ is the set of $G$-conjugacy classes in $d^{-1}(\alpha)$ for
some fixed $\alpha\in G/\ggeom$. Since $\ggeom$ is normal in $G$, this
set is indeed invariant under conjugation by the whole of $G$ (this is
an important point). We let $\rho_{\ell}$ be the induced map 
$$
Y\ra Y_{\ell}=\{\yc\in G_{\ell}^{\sharp}\,\mid\,
d(\yc)=p(\alpha)\}
\subset G_{\ell}^{\sharp}.
$$
\par
The natural density to consider (which arises in the sieve for
Frobenius) is still
$$
\nu_{\ell}(\yc)=\frac{|\yc|}{|\ggeom_{\ell}|},
\quad\text{  and hence } 
\quad
\nu_m(\yc)=\frac{|\yc|}{|\ggeom_m|}
$$
for a conjugacy class $\yc$. Note that this means that for any
conjugacy-invariant subset $\Omega_{\ell}\subset G_{\ell}$, union of
a set  $\Omega_{\ell}^{\sharp}$ of conjugacy classes such that
$\Omega_{\ell}^{\sharp}\subset d^{-1}(p(\alpha))=Y_{\ell}$, we have
$$
\nu(\Omega_{\ell}^{\sharp})=\frac{|\Omega_{\ell}|}{|\ggeom_{\ell}|}.
$$
\par
We turn to the question of finding a suitable orthonormal basis of
$L^2(Y_{\ell},\nu_{\ell})$. This is provided by the following general
lemma, which applies equally to $H=G_{\ell}$ and to $H=G_m$ for $m\in
S(\Lambda)$. 

\begin{lemma}\label{lm-group-theory}
Let $H$ be a finite group, $H^g$ a subgroup with abelian quotient
$\Gamma=H/H^g$. Let $\alpha\in\Gamma$ and $Y$ the set of conjugacy
classes of $G$ with image $\alpha$ in $\Gamma$.
\par
For an irreducible linear representation $\pi$ of $H$, let
$\varphi_{\pi}$ be the function
$$
\varphi_{\pi}\,:\,\yc\mapsto \Tr\pi(\yc)
$$
on $H^{\sharp}$.
\par
\emph{(1)} For $\pi$, $\tau$ irreducible linear representations of
$H$, we have
\begin{equation}\label{eq-scalar}
\langle \varphi_{\pi},\varphi_{\tau}\rangle =
\begin{cases}
0,&\text{if either $\varphi_{\pi}\mid H^g\not=\varphi_{\tau}\mid H^g$ or
  $\varphi_{\pi}\mid Y=0$,}\\ 
\overline{\psi(\alpha)}|\hat{\Gamma}^{\pi}|,&\text{ where
  $\psi\in\hat{\Gamma}$ satisfies 
  $\pi\otimes\psi\simeq\tau$, otherwise},
\end{cases}
\end{equation}
where $\hat{\Gamma}$ is the group of characters of $\Gamma$,
$\hat{\Gamma}^{\pi}=\{\psi\in\hat{\Gamma}\,\mid\, \pi\simeq
\pi\otimes\psi\}$, 
and the inner product is
$$
\langle f,g\rangle=\frac{1}{|H^g|}\sum_{\yc\in Y}{
|\yc|f(\yc)\overline{g(\yc)}}.
$$
\par
\emph{(2)} Let $\bon$ be the family of functions
$$
\yc\mapsto \frac{1}{\sqrt{|\hat{\Gamma}^{\pi}|}}\varphi_{\pi}(\yc),
$$
restricted to $Y$, where $\pi$ ranges over the subset of a set
of representatives for the equivalence relation 
$$
\pi\sim \tau\text{ if and only if }\pi\mid H^g\simeq \tau\mid  H^g,
$$
consisting of those representatives such that $\varphi_{\pi}\mid
Y\not=0$. Then $\bon$ is an orthonormal basis of $L^2(Y)$ for the
above inner product.
\end{lemma}

In the second case of~(\ref{eq-scalar}), the  
existence of the character $\psi$ will follow from the proof below.

\begin{proof}
We have
\begin{align*}
\langle\varphi_{\pi},\varphi_{\tau}\rangle&=\frac{1}{|H^g|}
\sum_{\stacksum{y\in H}{d(y)=\alpha}}{\Tr\pi(y)\overline{
\Tr\tau(y)}}\\
&=\frac{1}{|H^g|}\frac{1}{|\Gamma|}
\sum_{y\in H}{\Bigl(\sum_{\psi\in\hat{\Gamma}}{
\overline{\psi(\alpha)}\psi(y)\Bigr)\Tr\pi(y)\overline{
\Tr\tau(y)}}}\\
&=\sum_{\psi\in\hat{\Gamma}}{\overline{\psi(\alpha)}
\langle \pi\otimes\psi,\tau\rangle_H}
=\sum_{\psi\in\hat{\Gamma}}{\overline{\psi(\alpha)}
\delta(\pi\otimes\psi,\tau)},
\end{align*}
by orthogonality of characters of irreducible representations in
$L^2(H)$. 
\par
First of all, this is certainly zero unless there  exists at least
one $\psi$ such that $\pi\otimes\psi\simeq \tau$. In such a case we have
$\pi\mid H^g\simeq \tau\mid H^g$ since $H^g\subset \ker(\psi)$, so we
have shown that the condition $\pi\mid H^g\not\simeq \tau\mid H^g$
implies that the inner product is zero.
\par
Assume now that $\pi\mid H^g\simeq \tau\mid H^g$; then repeating the
above with $\alpha=1$ (i.e., $Y=H^g$), it follows from $\langle
\pi,\tau\rangle_{H^g}\not=0$ that there exists one $\psi$ at least
such that $\pi\otimes\psi=\tau$. 
\par
Fixing one such character $\psi_0$,
the characters $\psi'$ for which $\pi\otimes\psi'\simeq \tau$ are
given by $\psi'=\psi\psi_0$ where $\psi\in\hat{\Gamma}^{\pi}$. Then we
find
$$
\langle\varphi_{\pi},\varphi_{\tau}\rangle=
\sum_{\psi\in\hat{\Gamma}}{\overline{\psi(\alpha)}
\delta(\pi\otimes\psi,\pi\otimes\psi_0)}=
\overline{\psi_0(\alpha)}\sum_{\psi\in\hat{\Gamma}^{\pi}}{
\overline{\psi(\alpha)}}.
$$
\par
For any $\psi\in\hat{\Gamma}^{\pi}$ and $\yc\in Y$, we have the
character relation 
$$
\Tr\pi(\yc)=\psi(\yc)\Tr\pi(\yc)=\psi(\alpha)\Tr\pi(\yc),
$$
hence either $\psi(\alpha)=1$ for all $\psi$, or $\Tr\pi(\yc)=0$ for
all $\yc$,  i.e., $\varphi_{\pi}$ restricted to $Y$ vanishes. In this
last case, we have trivially $\varphi_{\tau}=0$ also on $Y$, and the
inner product vanishes.
\par
So we are led to the last case where $\pi\mid H^g=\tau\mid H^g$ but
$\psi(\alpha)=1$ for all $\psi\in\hat{\Gamma}^{\pi}$. Then the inner
product formula is clear from the above.
\par
Now to prove (2) from (1), notice first that the family $\bon$ is a
generating set of $L^2(Y)$ (indeed, all $\varphi_{\pi}$ generate
$L^2(H^{\sharp})$, but those $\pi$ for which $\varphi_{\pi}=0$ on $Y$
are clearly not needed, and if $\pi\sim\tau$, we have
$\varphi_{\tau}=\psi(\alpha) \varphi_{\pi}$ on $Y$, where $\psi$
satisfies $\tau\simeq\psi\otimes \pi$, so one element of each
equivalence class suffices for functions on $Y$). Then the fact that
we have an orthonormal basis follows from the inner product formula,
observing that if $\tau\simeq\pi\otimes\psi$, we have in fact
$\pi=\tau$ by definition of the equivalence relation, so $\psi=1$
in~(\ref{eq-scalar}). 
\end{proof}

\begin{example}\label{ex-biz}
In this lemma we emphasize that distinct representations of $H$
may give the same restriction on $H^g$, in which case they correspond to a
single element of the basis, and that it is possible that a
$\varphi_{\pi}$ vanish on $Y$, in which case their representative is
discarded from the basis.
\par
Take for instance $G=D_{n}$, a dihedral group of order $2n$. There is
an exact sequence 
$$
1\ra \Zz/n\Zz\ra G\fleche{d} \Zz/2\Zz\ra 1
$$
and if $Y=d^{-1}(1)\subset G$ and $\pi$ is any representation of $G$
of degree $2$, we have $\Tr\pi(x)=0$ for all $x\in Y$ (see
e.g.~\cite[5.3]{serre-rep}). 
\par
In particular, note that even though both cosets of $\Zz/n\Zz$ 
in $G$ have four elements, the sets of conjugacy classes in each do
not have the same cardinality (there are $5$ conjugacy classes, $3$ in
$\ker d$ and $2$ in the other coset). In other words, in a coset
sieve, the spaces $Y_m$ usually depend on the value of $\alpha$ (they
are usually not even isomorphic).
\end{example}

If we apply Lemma~\ref{lm-group-theory} to the groups $G_m$ and their
subgroups $\ggeom_m$, 
we clearly obtain orthonormal bases of $L^2(Y_m)$ containing the
constant function $1$, for the density $\nu_m$ above. Although it was
not phrased in this manner, this is what was used in~\cite{sieve}
(with minor differences, e.g., the upper bound $\kappa$ for the order
of $\hat{\Gamma}_m^{\pi}$ that occurs in
loc. cit., and can be removed -- as also noticed independently by
Zywina in a private email).
\par
As before we summarize the sieve statement:
\begin{proposition}\label{pr-coset-sieve}
Let $G$ be a group, $\ggeom$ a normal subgroup with abelian quotient,
$\rho_{\ell}\,:\, G\ra G_{\ell}$ a family of surjective
homomorphisms onto finite groups. Let 
$(Y,\Lambda,(\rho_{\ell}))$ be the coset sieve setting associated with
some $\alpha\in G/\ggeom$, and 
$(X,\mu,F)$ an associated siftable set. 
\par
For $m\in S(\Lambda)$, let  $\Pi_m$ be a set of 
representatives of the set of irreducible representations of $G_m$
modulo equality restricted to $\ggeom_m$, containing the constant
function $1$. Moreover, let $\Pi_m^*$ be the subset of primitive
representations, i.e., those such that when $\pi$ is decomposed as
$\boxtimes_{\ell\mid m}{\pi_{\ell}}$, no component $\pi_{\ell}$ is
trivial, and $\Tr\pi_{\ell}$ is not identically zero  on
$Y_{\ell}$. 
\par
Then, for any prime sieve support
$\Lb^*$ and associated sieve support $\Lb$, i.e., such
that~\emph{(\ref{eq-sieve-support})} holds, and for any conjugacy
invariant sifting sets $(\Omega_{\ell})$ with $\Omega_{\ell}\subset
Y_{\ell}$ for $\ell\in\Lb^*$, we have
$$
|S(X,\Omega;\Lb^*)|\leq \Delta H^{-1}
$$
where $\Delta$ is the smallest non-negative real number such that
$$
\sum_{m\in\Lb}{\sum_{\pi\in\Pi_m^*}{\Bigl|
\int_X{\alpha(x)\Tr\pi(\rho_m(F_x))d\mu(x)}
\Bigr|^2}}\leq \Delta\int_X{|\alpha(x)|^2d\mu(x)}
$$
for all square-integrable function $\alpha\in L^2(X,\mu)$, and where 
$$
H=\sum_{m\in\Lb}{\prod_{\ell\mid m}{\frac{|\Omega_{\ell}|}
{|\ggeom_{\ell}|-|\Omega_{\ell}|}}}.
$$
\par
Moreover we have
\begin{equation}\label{eq-max-cosets}
\Delta\leq
\max_{m\in\Lb}\max_{\pi\in\Pi_m^*}\sum_{n\in\Lb}\sum_{\tau\in\Pi_n^*}
{|W(\pi,\tau)|},
\end{equation}
where
$$
W(\pi,\tau)=\frac{1}{\sqrt{|\hat{\Gamma}_m^{\pi}|
|\hat{\Gamma}_n^{\tau}|}}\int_{X}{\Tr\pi(\rho_m(F_x))
\overline{\Tr\tau(\rho_n(F_x))}d\mu(x)}.
$$
\end{proposition}

We now consider what happens of the equidistribution approach in this
context. (Some of this also applies to group conjugacy sieves, where
$\ggeom=G$). 
\par
If we apply Lemma~\ref{lm-3} to the elements $\varphi_{\pi}$,
$\varphi_{\tau}$ of the basis $\bon_m$ and $\bon_n$ of $L^2(Y_m)$, we
see that the function $[\varphi_{\pi},\overline{\varphi_{\tau}}]$ defined
in~(\ref{eq-varphi-varphi}) is the character of the representation 
$$
[\pi,\bar{\tau}]=\pi_{m'}\boxtimes
(\pi_d\otimes\overline{\tau_d})\boxtimes \overline{\tau_{n'}}
$$
of $G_{[m,n]}$, already defined in~(\ref{eq-lcm-reps}). Hence we have
\begin{equation}\label{eq-wpitau}
W(\pi,\tau)=\frac{1}{\sqrt{|\hat{\Gamma}_m^{\pi}|
|\hat{\Gamma}_n^{\tau}|}}
\int_{X}{\Tr([\pi,\bar{\tau}]\rho_{[m,n]}(F_x))d\mu(x)}.
\end{equation}
\par
In applications, this means that to estimate the integrals
$W(\pi,\tau)$ it suffices (and may be more 
convenient) to be able to deal with integrals of the form
$$
\int_{X}{\Tr(\varpi(F_x))d\mu(x)}
$$
where $\varpi$  is a representation of $G$ that factors through a
finite product of groups $G_{\ell}$ (see
Section~\ref{sec-sieve-arith} for an instance of this).
\par
If we try to approach those integrals using the equidistribution
method, then the analogue of~(\ref{eq-rd}) is the identity
\begin{equation}\label{eq-rdb}
|\{\rho_d(F_x)=y^{\sharp}\}|
=\int_{\{\rho_d(F_x)=y^{\sharp}\}}{d\mu(x)}=
\frac{|y^{\sharp}|}{|\ggeom_d|}|X|+r_d(X;y^{\sharp}),
\end{equation}
defining $r_d(X;\yc)$ for $\yc\in Y_d$. Then~(\ref{eq-periodique}) becomes
$$
W(\pi,\tau)=\frac{|X|}{\sqrt{|\hat{\Gamma}_m^{\pi}|
|\hat{\Gamma}_n^{\tau}|}}
m([\pi,\bar{\tau}])
+O\Bigl(
\frac{1}{\sqrt{|\hat{\Gamma}_m^{\pi}|
|\hat{\Gamma}_n^{\tau}|}}
\sum_{\yc\in Y_{[m,n]}}
{\dim [\pi,\bar{\tau}]|r_{[m,n]}(X;y^{\sharp})|}
\Bigr)
$$
where, comparing with Lemma~\ref{lm-group-theory} with
$H=G_{[m,n]}$, we have 
\begin{equation}\label{eq-mpitau}
m([\pi,\bar{\tau}])=
\langle \varphi_{\pi},\varphi_{\tau}\rangle
\end{equation}
where the inner product is in $L^2(Y_{[m,n]})$ and both $\pi$ and
$\tau$ are extended to (irreducible) representations of $G_{[m,n]}$ by
taking trivial components at those $\ell\in [m,n]$ not in $m$ or $n$
respectively. Hence by~(\ref{eq-scalar}), we have
$m([\pi,\bar{\tau}])=0$  unless $\pi$ and $\tau$ thus extended are
isomorphic restricted to $\ggeom_{[m,n]}$, which clearly can occur
only if $m=n$ and then if $\pi=\tau$ by orthogonality of $\bon_n^*$.
When $(m,\pi)=(n,\tau)$, the inner product is equal to
$|\hat{\Gamma}_m^{\pi}|$ by~(\ref{eq-scalar}).
\par
Using this and~(\ref{eq-periodique}), we get
$$
W(\pi,\tau)=\delta(\pi,\tau)|X|+O
\Bigl(
\sum_{\yc\in Y_{[m,n]}}
{\dim[\pi,\bar{\tau}]|r_{[m,n]}(X;y^{\sharp})|}
\Bigr),
$$
where the implied constant is $\leq 1$. 
Hence  for any sieve support $\Lb$, the large sieve bound of
Proposition~\ref{pr-w-to-sieve} holds with
\begin{equation}\label{eq-delta-equi}
\Delta\leq |X|+R(X,\Lb)
\end{equation}
where
\begin{equation}\label{eq-ex-statement}
R(X;\Lb)=\max_{n\in \Lb}\max_{\pi\in\Pi_n^*}\Bigl\{
\sum_{m\in \Lb}\sum_{\tau\in\Pi_m^*}{
\sum_{\yc\in Y_{[m,n]}}
{\dim [\pi,\bar{\tau}]|r_{[m,n]}(X;y^{\sharp})|}\Bigr\}}.
\end{equation}

\par
For later reference, we also note the following fact:
\begin{lemma}\label{lm-ortho}
Let $m$, $n$ in $S(\Lambda)$, $\pi\in\Pi_m^*$, $\tau\in\Pi_n^*$.
The multiplicity of the
trivial representation in the restriction of $[\pi,\bar{\tau}]$ to
$\ggeom_{[m,n]}$ is equal to zero if
$(m,\pi)\not=(n,\tau)$, and is equal to
$|\hat{\Gamma}_m^{\pi}|$ if $(m,\pi)=(n,\tau)$. 
\end{lemma}

\begin{proof}
This multiplicity is by definition $\langle [\pi,\bar{\tau}],1\rangle$
computed in $L^2(\ggeom_{[m,n]})$, i.e., it is 
$\langle \varphi_{\pi},\varphi_{\tau}\rangle$ in $L^2(Y_{[m,n]})$ in
the case $\alpha=1\in G/\ggeom$ (with the same convention on extending
$\pi$ and $\tau$ to $G_{[m,n]}$ as before). So the result is a
consequence of Lemma~\ref{lm-group-theory}. 
\end{proof}

\section{Degrees and sums of degrees of representations of finite
  groups}
\label{sec-reprs}

This section is essentially independent from the rest of the
paper, and is devoted to proving some inequalities which are likely to
be useful in estimating quantities such as~(\ref{eq-max-cosets}) or
$R(X,\Lb)$ in~(\ref{eq-ex-statement}). Indeed, we will use them later
on in Section~\ref{sec-sieve-arith} and Section~\ref{sec-ls-ff}.
\par
In practice, the bound for the individual exponential sums
$W(\pi,\tau)$ is likely to involve the order of the groups
$G$ and the degrees of its representations, and their combination
in~(\ref{eq-max-cosets}) will involve sums of the degrees. For
instance, in the next sections, we will need to bound
\begin{gather*}
\max_{m,\pi}{\Bigl\{(\dim\pi)\sum_{n}{|G_{[m,n]}|
\sum_{\tau\in\Pi_n^*}{
(\dim \tau)}}}\Bigr\},\\
\max_{m,\pi}{\Bigl\{(\dim\pi)\sum_{n}{
\sum_{\tau\in\Pi_n^*}{(\dim \tau)}}}\Bigr\}.
\end{gather*}
\par
In applications, the groups $G_{\ell}$ are often (essentially) classical
linear groups over $\Fp_{\ell}$, but they are not entirely known (it
may only be known that they have bounded index in $GL(n,\Fp_{\ell})$
as $\ell$ varies, for instance¸ see~\cite{sieve} and
Section~\ref{sec-ls-ff}). Our results are biased to this case. 
\par
For a finite group $G$ and $p\in [1,+\infty]$, we denote
$$
A_p(G)=\Bigl(\sum_{\rho}{\dim(\rho)^p}\Bigr)^{1/p},\quad\text{ if
$p\not=+\infty$, } 
\quad
A_{\infty}(G)=\max\{\dim(\rho)\}
$$
where $\rho$ runs over irreducible linear representations of $G$ (in
characteristic zero). For example, we have $A_2(G)=\sqrt{|G|}$ for all
$G$ and if $G$ is abelian, then 
$A_p(G)=|G|^{1/p}$ for all $p$. Moreover
$$
\lim_{p\ra +\infty}{A_p(G)}=A_{\infty}(G).
$$
\par
We are primarily 
interested in $A_1(G)$ and $A_{\infty}(G)$, but $A_{5/2}(G)$ will also
occur in the proof of Theorem~\ref{th-non-conj}, and other cases may turn  
out to be useful in other sieve settings. 
We start with an easy monotonicity lemma.

\begin{lemma}\label{lm-incr}
Let $G$ be a finite group and $H\subset G$ a  subgroup, $p\in
[1,+\infty]$. We have 
$$
A_p(H)\leq A_p(G).
$$
\end{lemma}

\begin{proof}
For any  irreducible representation $\rho$ of $H$, choose
(arbitrarily) an irreducible representation $\pi(\rho)$ of $G$ 
that  occurs with positive
multiplicity in the induced representation
$\Ind_{H}^G\rho$. 
\par
Let  $\pi$ be a representation of $G$ in the image of
$\rho\mapsto\pi(\rho)$. For any $\rho$ where $\pi(\rho)=\pi$, we have 
$$
\langle
\rho,\res_{H}^G\pi\rangle_{H}=\langle \Ind_{H}^G\rho,\pi\rangle_G>0,
$$
by Frobenius reciprocity, i.e., all  $\rho$ with $\pi(\rho)=\pi$
occur in the restriction of $\pi$ to $H$. Hence for $p\not=+\infty$ we
obtain 
$$
\sum_{\stacksum{\rho}{\pi(\rho)=\pi}}
{\dim(\rho)^p}
\leq \Bigl(\sum_{\stacksum{\rho}{\pi(\rho)=\pi}}{\dim (\rho)}\Bigr)^p
\leq \dim(\pi)^p,
$$
and summing over all possible $\pi(\rho)$ gives the inequality
$$
A_p(H)^p\leq A_p(G)^p
$$
by positivity.  This settles the case $p\not=+\infty$, and the other
case only requires noticing that $\dim(\rho)\leq \dim(\pi(\rho))\leq
A_{\infty}(G)$. 
\end{proof}


We come to the main result of this section.  The terminology, which
may not be familiar to all readers, is explained by examples after the
proof. We hope that there will be no confusion between $p$ and the
characteristic of the finite field $\Fp_q$ which occurs...

\begin{proposition}\label{pr-bounds-p}
Let $\G/\Fp_{q}$ be a split connected reductive linear algebraic
group of dimension $d$ and rank $r$ over a finite field, with
connected center. Let $W$ be its Weyl group and $G=\G(\Fp_{q})$ the
finite group of rational points of $\G$.
\par
\emph{(1)} 
For any subgroup $H\subset G$ and $p\in [1,+\infty]$, we have
$$
A_p(H)\leq (q+1)^{(d-r)/2+r/p}\Bigl(1+\frac{2r|W|}{q-1}\Bigr)^{1/p},
$$
with the convention $r/p=0$ if $p=+\infty$, in particular the second
factor is $=1$ for $p=+\infty$.
\par
\emph{(2)} If $\G$ is a product of groups of type $A$ or $C$, i.e., of
linear and symplectic groups, then
$$
A_p(H)\leq (q+1)^{(d-r)/2+r/p}.
$$
\end{proposition}

The proof is based on a simple interpolation argument from the extreme
cases $p=1$, $p=+\infty$. Indeed by Lemma~\ref{lm-incr} we can clearly
assume $H=G$ and by writing the obvious inequality
$$
A_p(G)^p=\sum_{\rho}{\dim(\rho)^p}\leq A_{\infty}(G)^{p-1} A_1(G),
$$
we see that it suffices to prove the following:

\begin{proposition}\label{pr-degrees}
Let $\G/\Fp_{q}$ be a split connected reductive linear algebraic
group of dimension $d$ with connected center, and
let $G=\G(\Fp_{q})$ be the finite group of its rational
points. Let $r$ be the rank of  $\G$. Then
we have 
\begin{equation}\label{eq-stronger}
  A_{\infty}(G)\leq\frac{|G|_{p'}}{(q-1)^r}\leq
  (q+1)^{(d-r)/2},\quad\text{ and } 
  \quad
  A_1(G)\leq (q+1)^{(d+r)/2}\Bigl(1+\frac{2r|W|}{q-1}\Bigr),
\end{equation}
where $n_{p'}$ denotes the prime-to-$p$ part of a rational number $n$,
$p$ being the characteristic of $\Fp_q$. Moreover, if the principal
series of $G$ is not empty\footnote{\ In particular if $q$ is large
  enough given $\Gg$.}, there is equality 
$$
A_{\infty}(G)=\frac{|G|_{p'}}{(q-1)^r}
$$
and $\dim\rho=A_{\infty}(G)$ if and only if $\rho$ is in the
principal series.
\par
Finally if $\G$ is a product of groups of type $A$ or $C$, then the
factor $(1+2r|W|/(q-1))$ may be removed in the bound for $A_1(G)$.
\end{proposition}

It seems very possible that the factor $(1+2r|W|/(q-1))$ could always
be removed, but we haven't been able to figure this out using
Deligne-Lusztig characters, and in fact for groups of type $A$ or $C$,
we simply quote \emph{exact formulas} for $A_1(G)$ due to Gow,
Klyachko and Vinroot, which are proved in completely different
ways.\footnote{\ The ``right'' upper bound for the case of groups of type $A$
(i.e, $GL(r)$) may
be recovered using the structure of unipotent representations of such
groups.}  The extra factor is not likely to be a problem in many
applications where $q\ra +\infty$, but it may be questionable for
uniformity with respect to the rank. 
\par
The ideas in the proof were suggested and explained by J. Michel. 


\begin{proof}
This is based on properties of the Deligne-Lusztig
generalized characters.
We will mostly refer to~\cite{digne-michel} and \cite{carter}
for all facts which are needed (using notation 
from~\cite{digne-michel}, except for writing simply $G$ for what is
denoted $\G^F$ there). We identify irreducible representations of $G$
(up to isomorphism) with their characters seen as complex-valued
functions on $G$.  
\par
First, for a connected reductive group $\G/\Fp_q$ over a
finite field, Deligne and Lusztig have constructed (see
e.g.~\cite[11.14]{digne-michel})  a family
$R_{\Tt}^{\G}(\theta)$ of generalized representations of
$G=\G(\Fp_q)$ (i.e., linear combinations with integer
coefficients of ``genuine'' representations of $G$), parametrized by
pairs $(\Tt,\theta)$ consisting of a maximal torus $\Tt\subset \Gg$
defined over $\Fp_q$ and a (one-dimensional) character
$\theta$ of the finite abelian group $T=\Tt(\Fp_q)$. The
$R_{\Tt}^{\G}(\theta)$ are not 
all irreducible, but any irreducible character occurs
(with positive or negative multiplicity)  in
the decomposition of at least one such character. Moreover,
$R_{\Tt}^{\G}(\theta)$ only depends (up to isomorphism) on the $G$-conjugacy
class of the pair $(\Tt,\theta)$. 
\par
We quote here a useful classical fact:
for any $\Tt$ we have 
\begin{equation}
  \label{eq-order-torus}
  (q-1)^r\leq |T|\leq (q+1)^r
\end{equation}
(see e.g.~\cite[13.7 (ii)]{digne-michel}), and moreover $|T|=(q-1)^r$
if and only if $\Tt$ is a split torus (i.e., $\Tt\simeq \Gg_m^r$ over
$\Fp_q$). Indeed, we have
$$
|T|=|\det(q^n-w\mid Y_0)|
$$ 
where $w\in W$ is such that $\Tt$ is obtained from a split torus
$\Tt_0$ by ``twisting with $w$'' (see
e.g.~\cite[Prop. 3.3.5]{carter}), and $Y_0\simeq \Zz^r$ is the group
of cocharacters of $\Tt_0$. If $\lambda_1$, \ldots, $\lambda_r$ are the
eigenvalues of $w$ acting on $Y_0$, which are roots of unity, then we
have 
$$
|T|=\prod_{i=1}^r{(q-\lambda_i)},
$$ 
and so $|T|=(q-1)^r$ if and only if each $\lambda_i$ is equal to $1$,
if and only if $w$ acts trivially on $Y_0$, if and only if $w=1$ ($W$
acts faithfully on $Y_0$) and $\Tt$ is split.
\par
As in~\cite[12.12]{digne-michel}, we denote by $\rho\mapsto p(\rho)$ the
orthogonal projection of the space $\mathcal{C}(G)$ of
complex-valued conjugacy-invariant functions on $G$ to the subspace
generated by Deligne-Lusztig characters, where $\mathcal{C}(G)$ is
given the standard inner product
$$
\langle f,g\rangle =\frac{1}{|G|}\sum_{x\in G}{f(x)\overline{g(x)}},
$$
and for a representation $\rho$, we of course denote
$p(\rho)=p(\Tr\rho)$ the projection of its character.
\par
For any representation $\rho$, we have
$\dim(\rho)=\dim(p(\rho))$, where $\dim(f)$, for an arbitrary function
$f\in\mathcal{C}(G)$ is obtained by linearity from the degree of
characters. Indeed, for any $f$ standard character
theory shows that
$$
\dim(f)=\langle f,\reg_G\rangle
$$
where $\reg_G$ is the regular representation of $G$. 
From~\cite[12.14]{digne-michel}, the regular representation is in the
subspace spanned by the Deligne-Lusztig characters, so by definition
of an orthogonal projector we have
$$
\dim(\rho)=\langle \rho,\reg_G\rangle=\langle p(\rho),
\reg_G\rangle=\dim(p(\rho)).
$$
\par
Now because the characters $R_{\Tt}^{\G}(\theta)$ for distinct conjugacy
classes of $(\Tt,\theta)$ are 
orthogonal (see e.g.~\cite[11.15]{digne-michel}), we can write
$$
p(\rho)=\sum_{(\Tt,\theta)}{
\frac{\langle \rho,R_{\Tt}^{\G}(\theta)\rangle}{
\langle R_{\Tt}^{\G}(\theta),R_{\Tt}^{\G}(\theta)\rangle}
R_{\Tt}^{\G}(\theta)}
$$
(sum over all distinct Deligne-Lusztig characters) and so
$$
\dim(p(\rho))=\sum_{(\Tt,\theta)}{
\frac{\langle \rho,R_{\Tt}^{\G}(\theta)\rangle}{
\langle R_{\Tt}^{\G}(\theta),R_{\Tt}^{\G}(\theta)\rangle}
\dim(R_{\Tt}^{\G}(\theta))}.
$$
By~\cite[12.9]{digne-michel} we have 
\begin{equation}\label{eq-dim-rtt}
\dim(R_{\Tt}^{\G}(\theta))=\eps_{\G}\eps_{\Tt}|G|_{p'}|T|^{-1},
\end{equation} 
where
$\eps_{\G}=(-1)^{r}$ and $\eps_{\Tt}=(-1)^{r(\Tt)}$,
$r(\Tt)$ being the 
$\Fp_q$-rank of $\Tt$ (see~\cite[p. 66]{digne-michel} for the
definition). This yields the formula
\begin{equation}
  \label{eq-formula-dim}
  \dim(p(\rho))=|G|_{p'}\sum_{(\Tt,\theta)}{
\frac{1}{|T|}
\frac{\langle \rho,\eps_{\Gg}\eps_{\Tt}R_{\Tt}^{\G}(\theta)\rangle}
{\langle R_{\Tt}^{\G}(\theta),R_{\Tt}^{\G}(\theta)\rangle}}.
\end{equation}
\par
Now we use the fact that pairs $(\Tt,\theta)$ are partitioned in
\emph{geometric conjugacy classes}, defined as follows: two pairs
$(\Tt,\theta)$ and $(\Tt',\theta')$ are geometrically conjugate if and
only if there exists $g\in \Gg(\bar{\Fp}_q)$ such that
$\Tt=g\Tt'g^{-1}$ and for all $n$ such that $g\in \Gg(\Fp_{q^n})$, we
have
$$
\theta(N_{\Fp_{q^n}/\Fp_q}(x))=\theta'(N_{\Fp_{q^n}/\Fp_q}(g^{-1}xg))
\quad\quad\text{ for } x\in\Tt(\Fp_{q^n}),
$$
(see e.g.~\cite[13.2]{digne-michel}). The point is the following
property of geometric conjugacy classes: if 
the generalized characters $R_{\Tt}^{\Gg}(\theta)$ and
$R_{\Tt'}^{\Gg}(\theta')$ have a common irreducible component, then
$(\Tt,\theta)$ and $(\Tt',\theta')$ are geometrically conjugate
(see e.g.~\cite[13.2]{digne-michel}). 
\par
In particular, for a given
$\rho$, if $\langle
\rho,R_{\Tt}^{\Gg}(\theta)\rangle$ is non-zero for some
$(\Tt,\theta)$, then only pairs $(\Tt',\theta')$
geometrically conjugate to $(\Tt,\theta)$ may satisfy $\langle
\rho,R_{\Tt'}^{\Gg}(\theta)\rangle\not=0$. So we have
$$
  \dim(p(\rho))=|G|_{p'}\sum_{(\Tt,\theta)\in\kappa}{
\frac{1}{|T|}
\frac{\langle \rho,\eps_{\Gg}\eps_{\Tt}R_{\Tt}^{\G}(\theta)\rangle}
{\langle R_{\Tt}^{\G}(\theta),R_{\Tt}^{\G}(\theta)\rangle}},
$$
for some geometric conjugacy class $\kappa$, depending on $\rho$. By
Cauchy-Schwarz, we obtain
\begin{equation}
\label{eq-cauchy}
\dim(p(\rho))\leq 
|G|_{p'}
\Bigl(\sum_{(\Tt,\theta)\in\kappa}{
\frac{1}{|T|^2}\frac{1}{\langle R_{\Tt}^{\G}(\theta),
R_{\Tt}^{\G}(\theta)\rangle}}\Bigr)^{1/2}
\Bigl(\sum_{(\Tt,\theta)\in\kappa}{
\frac{|\langle \rho,R_{\Tt}^{\G}(\theta)\rangle|^2}
{\langle R_{\Tt}^{\G}(\theta),R_{\Tt}^{\G}(\theta)\rangle}}
\Bigr)^{1/2}.
\end{equation}
\par
The second term on the right is simply $\langle
p(\rho),p(\rho)\rangle\leq \langle \rho,\rho\rangle=1$. As for the
first term we have
$$
\sum_{(\Tt,\theta)\in\kappa}{
\frac{1}{|T|^2}\frac{1}{\langle R_{\Tt}^{\G}(\theta),
R_{\Tt}^{\G}(\theta)\rangle}}
\leq \frac{1}{(q-1)^{2r}}\sum_{(\Tt,\theta)\in\kappa}{
\frac{1}{\langle R_{\Tt}^{\G}(\theta),
R_{\Tt}^{\G}(\theta)\rangle}}
$$
by~(\ref{eq-order-torus}). Now it is known that for each class
$\kappa$, the assumption 
that $\Gg$ has connected center implies that the generalized
character
$$
\chi(\kappa)=\sum_{(\Tt,\theta)\in\kappa}
{\frac{\eps_{\G}\eps_{\Tt}R_{\Tt}^{\G}(\theta)}
{\langle  R_{\Tt}^{\G}(\theta),R_{\Tt}^{\G}(\theta)\rangle}}
$$
is in fact an irreducible character of $G$ (such characters are called
\emph{regular} characters; see e.g.~\cite[Prop. 8.4.7]{carter}).  This
implies that 
$$
\sum_{(\Tt,\theta)\in\kappa}{
\frac{1}{\langle R_{\Tt}^{\G}(\theta),
R_{\Tt}^{\G}(\theta)\rangle}}
=\langle\chi(\kappa),\chi(\kappa)\rangle=1,
$$
and so we have
\begin{equation}\label{eq-ainf-almost}
\dim p(\rho)\leq \frac{|G|_{p'}}{(q-1)^r}.
\end{equation}
\par
Now observe that we will have equality in this argument if
$\rho$ is itself of the form $\pm R_{\Tt}^{\Gg}(\theta)$, and if
$|T|=(q-1)^r$. Those conditions hold for representations of the
principal series, i.e., characters $R_{\Tt}^{\Gg}(\theta)$ for an
$\Fp_q$-split torus $\Tt$ and a character $\theta$ ``in general
position'' (see e.g.~\cite[Cor. 7.3.5]{carter}). Such characters are
also, more elementarily, induced characters 
$\Ind_B^G(\theta)$, where $B=\Bb(\Fp_q)$ is a Borel subgroup
containing $T$, for some Borel subgroup $\Bb$ defined over $\Fp_q$
containing $\Tt$ (which exist for a split torus $\Tt$) and $\theta$
is extended to $B$ by setting $\theta(u)=1$ for unipotent elements
$u\in B$. For this, see e.g.~\cite[Prop. 2.6]{lusztig}.
\par
Conversely, let $\rho$ be such that
$$
\dim\rho=\frac{|G|_{p'}}{(q-1)^r}
$$
and let $\kappa$ be the associated geometric conjugacy class. From the
above, for any $(\Tt,\theta)$ in
$\kappa$, we have  $|T|=(q-1)^r$, i.e., $\Tt$ is $\Fp_q$-split. Now it
follows from 
Lemma~\ref{lm-geo-conj} (probably well-known) that this implies that
$R_{\Tt}^{\Gg}(\theta)$ is an irreducible
representation, so must be equal to
$\rho$. 
\par
We now come to $A_1(G)$. To deal with the fact that
in~(\ref{eq-formula-dim}), $|T|$ depends on 
$(\Tt,\theta)\in\kappa$, we write
\begin{equation}\label{eq-depend}
 \dim(p(\rho))=\frac{|G|_{p'}}{(q-1)^r}\sum_{\kappa}{
\langle\rho,\chi(\kappa)\rangle}
+|G|_{p'}\sum_{(\Tt,\theta)}{\Bigl(
\frac{1}{|T|}-\frac{1}{(q-1)^r}
\Bigr)
\frac{\eps_{\Gg}\eps_{\Tt}\langle\rho,R_{\Tt}^{\G}(\theta)\rangle}
{\langle
  R_{\Tt}^{\G}(\theta),R_{\Tt}^{\G}(\theta)\rangle}}
\end{equation}
(since by~(\ref{eq-order-torus}), the dependency is rather weak).
\par
Now summing over $\rho$, consider the first term's contribution. Since
$\chi(\kappa)$ is an 
irreducible character, the sum
$$
\sum_{\rho}{\sum_{\kappa}{\langle\rho,\chi(\kappa)\rangle}}
$$
is simply the number of geometric conjugacy classes. This is given by
$q^{r'}|Z|$ by~\cite[14.42]{digne-michel} 
or~\cite[Th. 4.4.6 (ii)]{carter}, where $r'$ is  
the semisimple rank of $\Gg$ and $Z=Z(\Gg)(\Fp_q)$ is the group of
rational points of the center of $\Gg$. For this quantity, note that
the center of $\Gg$ being 
connected implies that $Z(\Gg)$ is the radical of $\Gg$ (see
e.g.~\cite[Pr. 7.3.1]{springer}) so $Z(\Gg)$ is a torus and $r=r'+\dim
Z(\Gg)$. So using again the bounds~(\ref{eq-order-torus}) for the
cardinality of the group of rational points of a torus, we obtain
\begin{equation}
\label{eq-conj-geo}
|Z|q^{r'}\leq (q+1)^r.
\end{equation}
\par
To estimate the sum of the contributions in the second term, say $\sum
t(\rho)$, we write 
$$
\sum_{\rho}{t(\rho)}=
|G|_{p'}\sum_{(\Tt,\theta)}{\Bigl(
\frac{1}{|T|}-\frac{1}{(q-1)^r}
\Bigr)
\frac{\eps_{\Gg}\eps_{\Tt}\langle\sum_{\rho}{\rho},
R_{\Tt}^{\G}(\theta)\rangle}
{\langle
  R_{\Tt}^{\G}(\theta),R_{\Tt}^{\G}(\theta)\rangle}},
$$
and we bound
\begin{equation}\label{eq-to-refine}
\Bigl|\langle\sum_{\rho}{\rho},
R_{\Tt}^{\G}(\theta)\rangle\Bigr|\leq \langle
R_{\Tt}^{\G}(\theta),R_{\Tt}^{\G}(\theta)\rangle
\end{equation}
for any $(\Tt,\theta)$, since we can write
$$
R_{\Tt}^{\G}(\theta)=\sum_{\rho}{a(\rho)\rho}\quad\text{with }
a(\rho)\in\Zz, 
$$
and therefore
\begin{equation}\label{eq-to-refinebis}
\Bigl|\langle\sum_{\rho}{\rho},
R_{\Tt}^{\G}(\theta)\rangle\Bigr|=\Bigl|\sum_{\rho}{a(\rho)}\Bigr|\leq
\sum_{\rho}{|a(\rho)|^2}=\langle
R_{\Tt}^{\G}(\theta),R_{\Tt}^{\G}(\theta)\rangle.
\end{equation}
Thus
$$
\sum_{\rho}{t(\rho)}\leq \frac{|G|_{p'}}{(q-1)^r}\frac{2r}{q-1}
|\{(\Tt,\theta)\}|.
$$
There are at most $|W|$ different choices of $\Tt$ up to
$G$-conjugacy, and for each there are at most $|T|\leq (q+1)^r$
different characters, and so we have
\begin{equation}\label{eq-sum-trho}
\sum_{\rho}{t(\rho)}\leq
\frac{|G|_{p'}}{(q-1)^r}\frac{2r|W|}{q-1}(q+1)^r,
\end{equation}
and
\begin{equation}\label{eq-a1-almost}
\sum_{\rho}{\dim\rho}\leq (q+1)^r\frac{|G|_{p'}}{(q-1)^r}
\Bigl(1+\frac{2r|W|}{q-1}\Bigr).
\end{equation}
\par
To conclude, we use the classical formula
$$
|G|=q^N\prod_{1\leq i\leq r}{(q^{d_i}-1)},
$$
where $N$ is the number of positive roots of $\G$, and the
$d_i$ are the degrees of invariants of the Weyl group (this is
because $\Gg$ is split; see e.g.~\cite[2.4.1 (iv); 2.9,
p. 75]{carter}). So  
$$
|G|_{p'}=\prod_{1\leq i\leq r}{(q^{d_i}-1)}
$$
and 
\begin{equation}\label{eq-bd-gt}
\frac{|G|_{p'}}{(q-1)^r}=\prod_{1\leq i\leq r}{\frac{q^{d_i}-1}{q-1}}
\leq \prod_{1\leq i\leq r}{(q+1)^{d_i-1}}=(q+1)^{\sum{(d_i-1)}}=
(q+1)^{(d-r)/2},
\end{equation}
since $\sum{(d_i-1)}=N$  and $N=(d-r)/2$ (see
e.g.~\cite[2.4.1]{carter},~\cite[8.1.3]{springer}).  
\par
Inserting this in~(\ref{eq-ainf-almost}) we derive the first inequality
in~(\ref{eq-stronger}), and with~(\ref{eq-a1-almost}), we get
$$
A_1(G)\leq (q+1)^{(d+r)/2}\Bigl(1+\frac{2r|W|}{q-1}\Bigr),
$$
which is the second part of~(\ref{eq-stronger}).
\par
Now we explain why the extra factor involving the Weyl group can be
removed for products of groups of type $A$ and $C$.  Clearly it
suffices to work with $\G=GL(n)$ and $\G=CSp(2g)$. 
\par
For $\G=GL(n)$, with $d=n^2$ and $r=n$, Gow~\cite{gow} and
Klyachko~\cite{klyachko} have proved 
independently that  
$A_1(G)$ is equal to the number of symmetric matrices in $G$. The
bound
$$
A_1(G)\leq (q+1)^{(n^2+n)/2}
$$
follows immediately. 
\par
For $\G=CSp(2g)$, with $d=2g^2+g+1$ and $r=g+1$, the exact  analog of
Gow's theorem is due to Vinroot~\cite{vinroot}. Again,
Vinroot's result implies $A_1(G)\leq (q+1)^{(d+r)/2}$ in this case
(see~\cite[Cor 6.1]{vinroot}, and use the formulas for the order of unitary
and linear groups to check the final bound). 
\end{proof}

Here is the lemma used in the determination of $A_{\infty}(G)$ when
there is a character in general position of a split torus:

\begin{lemma}
\label{lm-geo-conj}
Let $\G/\Fp_{q}$ be a split connected reductive linear algebraic
group of dimension $d$ and
let $G=\G(\Fp_{q})$ be the finite group of its rational
points. Let $\Tt$ be a split torus in $\Gg$, $\theta$ a character of
$T=\Tt(\Fp_q)$. If $\Tt'$ is also a split torus for any pair $(\Tt',\theta')$
geometrically conjugate to $(\Tt,\theta)$, then
$R_{\Tt}^{\Gg}(\theta)$ is irreducible. 
\end{lemma}

\begin{proof}
If $R_{\Tt}^{\Gg}(\theta)$  is not irreducible, then by the inner
product formula for 
Deligne-Lusztig characters, there exists $w\in W$, $w\not=1$, such
that ${}^w\theta=\theta$ (see e.g.~\cite[Cor. 11.15]{digne-michel}). 
Let $\Tt'$ be a torus obtained from $\Tt$ by ``twisting by $w$'',
i.e., $\Tt'=g\Tt g^{-1}$ where $g\in\Gg$ is such that
$g^{-1}\frob(g)=w$ (see e.g.~\cite[3.3]{carter}). Let
$Y=\Hom(\Gg_m,\Tt)\simeq \Zz^r$ (resp. $Y'$) be 
the abelian group of cocharacters of $\Tt$ (resp. $\Tt'$); the
conjugation isomorphism $\Tt\ra\Tt'$ gives rise to a conjugation
isomorphism $Y\ra Y'$ (loc. cit.). Moreover, there is an action
of the Frobenius $\frob$ on $Y$ and a canonical
isomorphism $T\simeq Y/(\frob-1)Y$ (see
e.g.~\cite[Prop. 13.7]{digne-michel}), 
hence canonical 
isomorphisms of the character groups $\hat{T}$ and $\hat{T}'$ as
subgroups of the characters groups of $Y$ and $Y'$:
$$
\hat{T}\simeq \{\chi\,:\, Y\ra \Cc^{\times}\,\mid\, 
(\frob-1)Y\subset \ker\chi\},
\quad
\hat{T}'\simeq \{\chi\,:\, Y'\ra \Cc^{\times}\,\mid\, 
(\frob-1)Y'\subset \ker\chi\}.
$$
\par
Unraveling the definitions, a simple calculation shows that the
condition ${}^w\theta=\theta$ is precisely what is needed to prove
that the character $\chi$ of $Y$ associated to $\theta$, when
``transported'' to a character $\chi'$ of $Y'$ by the conjugation
isomorphism, still satisfies $\ker \chi'\supset (\frob-1)Y'$
(see in particular~\cite[Prop. 3.3.4]{carter}), so is associated with
a character $\theta'\in\hat{T}'$. 
\par
Using the characterization of geometric conjugacy
in~\cite[Prop. 13.8]{digne-michel}, it is
then clear  that $(\Tt,\theta)$ is
geometrically conjugate to $(\Tt',\theta')$, and since $w\not=1$, the
torus $\Tt'$ is not split. So by contraposition, the lemma is proved.
\end{proof}

\begin{example}\label{ex-a}
(1) Let $\ell$ be prime, $r\geq 1$ and let
$\G=GL(r)/\Fp_{\ell}$. Then 
$G=GL(r,\Fp_{\ell})$, $\G$ 
is a split connected reductive of rank $r$  and dimension $r^2$, with
connected center of dimension $1$. So from
Lemma~\ref{lm-incr} and Proposition~\ref{pr-bounds-p}, we get
$$
A_p(H)\leq (\ell+1)^{r(r-1)/2+r/p}
$$
for $p\in [1,+\infty]$ for any subgroup $H$ of $G$, and in particular
$$
A_{\infty}(H)\leq (\ell+1)^{r(r-1)/2}\quad\text{ and }\quad
A_1(H)\leq (\ell+1)^{r(r+1)/2}.
$$
\par
It would be interesting to know if there are
other values of $p$ besides $p=1$, $2$ and $+\infty$ (the latter when
$q$ is large enough) for which $A_p(GL(n,\Fp_{q}))$ can be computed
exactly. 
\par
(2) Let $\ell\not=2$ be prime, $g\geq 1$ and let
$\G=CSp(2g)/\Fp_{\ell}$. Then $G=CSp(2g,\Fp_{\ell})$ and $G$
is a split connected reductive group of rank $g+1$ and dimension
$2g^2+g+1$, with connected center. So from
Lemma~\ref{lm-incr} and Proposition~\ref{pr-bounds-p}, we get
$$
A_p(H)\leq (\ell+1)^{g^2+(g+1)/p}
$$
for $p\in [1,+\infty]$ for any subgroup $H$ of $G$, and in particular
$$
A_{\infty}(H)\leq (\ell+1)^{g^2}\quad\text{ and }\quad
A_1(H)\leq (\ell+1)^{g^2+g+1}.
$$
\end{example}

In the case of $G=SL(r,\Fp_q)$ or $G=Sp(2g,\Fp_q)$, which correspond
to $\G$ where the center is not connected, the bound for
$A_{\infty}(G)$ given by this example is still sharp if we see $G$ as
subgroup of 
$GL(r,\Fp_q)$ or $CSp(2g,\Fp_q)$, because both $d$ and
$r$ increase by $1$, so $d-r$ doesn't change. 
However, for $A_1(G)$, the exponent increases by one. 
Here is a slightly different argument
that almost recovers the  ``right'' bound.

\begin{lemma}\label{lm-sl-sp}
Let $\Gg=SL(n)$ or $Sp(2g)$ over $\Fp_q$, $d$ the dimension and $r$
the rank of $\Gg$, and $G=\Gg(\Fp_q)$. 
Then we have the following bounds
$$
A_p(G)\leq \kappa^{1/p}(q+1)^{(d-r)/2+r/p}
\Bigl(\frac{q+1}{q-1}\Bigr)^{1/p}.
$$
and
$$
A_p(G)\leq 
(q+1)^{(d-r)/2+r/p}\Bigl(\frac{q+1}{q-1}\Bigr)^{1/p}\Bigl(1
+\frac{2\kappa(r+1)|W|}{q-1}\Bigr)^{1/p}
$$
for any $p\in [1,+\infty]$, where $\kappa=n$ for $SL(n)$ and
$\kappa=2$ for $Sp(2g)$. 
\end{lemma}

The first bound is better for fixed $q$, whereas the second is almost
as sharp as the bound for $GL(n)$ or $CSp(2g)$ if $q$ is large.

\begin{proof}
As we observed before the statement, this holds for $p=+\infty$, so it
suffices to consider $p=1$ and then use the same interpolation
argument as for Proposition~\ref{pr-bounds-p}.
\par
Let $\Gg_1=GL(n)$ or $CSp(2g)$ for $\Gg=SL(n)$ or $Sp(2g)$
respectively, $G_1=\Gg_1(\Fp_q)$. We use the exact sequence 
$$
1\ra G\ra G_1\fleche{m} \Gamma=\Fp_q^{\times}\ra 1
$$
(compare
with Section~\ref{sec-spec-coset})
where $m$ is either the determinant or the multiplicator of a
symplectic similitude. Let $\rho$ be an irreducible representation of
$G$, and as in the proof of Lemma~\ref{lm-incr}, let $\pi(\rho)$ be
any irreducible representation of $G_1$ in the induced representation
to $G_1$. The point is that all ``twists'' $\pi(\rho)\otimes\psi$,
where $\psi$ is a character of $\Fp_q^{\times}$ lifted to $G_1$
through $m$, are isomorphic restricted to $G$, and hence each
$\pi(\rho)\otimes\psi$  contains $\rho$ when restricted to $G$, and
contains even all $\rho$ with the same $\pi(\rho)$. So if
$\pi\sim\pi'$, for  representations of $G_1$, denotes isomorphism when
restricted to $G$, we have 
$$
A_1(G)\leq \sum_{\{\pi\}/\sim}{\dim \pi}
$$
where the sum is over a set of representatives for this equivalence
relation. On the other hand, $\dim\pi=\dim\pi'$ for $\pi\sim \pi'$,
and for each $\pi$ there are $|\hat{\Gamma}/\hat{\Gamma}^{\pi}|$
distinct  representations equivalent to $\pi$, with notation as in
Lemma~\ref{lm-group-theory}.
Hence,
$$
A_1(G)\leq \frac{1}{q-1}
\sum_{\pi}{|\hat{\Gamma}^{\pi}|\dim\pi}.
$$
\par
From, e.g.,~\cite[Lemma 2.3]{sieve}, we know that
$\hat{\Gamma}^{\pi}$ has order at most $n$ (for $SL(n)$) or $2$ (for
$Sp(2g)$), which by applying Proposition~\ref{pr-bounds-p} yields the
first bound\footnote{\ This suffices for the
  applications in this paper.}, namely
$$
A_1(G)\leq \kappa\frac{(q+1)^{(d+r)/2}}{q-1},
\quad\quad\text{ with $\kappa=2$  or $n$}.
$$
\par
To obtain the refined bound, observe that in the
formula~(\ref{eq-depend}) for the dimension of an irreducible
representation $\rho$ of $G_1$, the first term is zero unless $\rho$
is a regular representation, and the second $t(\rho)$ is smaller by a
factor roughly $q$. If $\pi$ is regular, we have
$\hat{\Gamma}^{\pi}=1$ by Lemma~\ref{lm-no-twist} below. So it follows
that 
\begin{align*}
A_1(G)&\leq \frac{1}{q-1}\Bigl\{\sum_{\pi\text{ regular}}{\dim\pi}
+\kappa\sum_{\pi\text{ not regular}}{\dim\pi}
\Bigr\}\\
&\leq \frac{A_{\infty}(G_1)}{q-1}q^r(q-1)
+\kappa \sum_{\pi\text{ not regular}}{t(\rho)}
\end{align*}
(in the first term, $q^r(q-1)$ is the number of geometric
conjugacy classes for $G_1$, computed as in~(\ref{eq-conj-geo}), since
$r$ is the semi-simple rank of $G_1$). We have the analogue
of~(\ref{eq-sum-trho}): 
$$
\sum_{\pi\text{ not regular}}{t(\pi)}\leq
\frac{|G_1|_{p'}}{(q-1)^{r+1}}\frac{2(r+1)|W|}{q-1}(q+1)^{r+1}\leq
2(r+1)|W|\frac{(q+1)^{(d+r)/2+1}}{q-1},
$$
by~(\ref{eq-bd-gt}) (because 
$$
\Bigl|\langle\sum_{\pi\text{ not regular}}{\pi},
R_{\Tt}^{\G_1}(\theta)\rangle\Bigr|\leq \langle
R_{\Tt}^{\G_1}(\theta),R_{\Tt}^{\G_1}(\theta)\rangle,
$$
see~(\ref{eq-to-refinebis}), and the same argument leading to
(\ref{eq-to-refine})).
The bound
$$
A_1(G)\leq (q+1)^{(d+r)/2}\Bigl(1+\frac{2\kappa(r+1)|W|}{q-1}\Bigr)
$$
follows.
\end{proof}

\begin{lemma}\label{lm-no-twist}
Let $\Gg=GL(n)$ or $CSp(2g)$ over $\Fp_q$, $G=\Gg(\Fp_q)$.  For any
regular irreducible representation $\rho$ of $G$, we have
$\hat{\Gamma}^{\rho}=1$.
\end{lemma}

\begin{proof}
As above, let $m\,:\, \G\ra \G_m$ be the determinant or multiplicator
character. Let $\rho$ be a regular representation and $\psi$ a
character of $\Fp_q^{\times}$ such that $\rho\otimes\psi\simeq \rho$,
where $\psi$ is shorthand for $\psi\circ m$. We wish to show that
$\psi$ is trivial to conclude $\hat{\Gamma}^{\rho}=1$. For this
purpose, write  
$$
\rho=\sum_{(\Tt,\theta)\in\kappa}{
\frac{\eps_{\G}\eps_{\Tt}R_{\Tt}^{\G}(\theta)}
{\langle R_{\Tt}^{\G}(\theta),R_{\Tt}^{\G}(\theta)\rangle }
}
$$
for some unique geometric conjugacy class
$\kappa$. We have
$R_{\Tt}^{\G}(\theta)\otimes\psi=R_{\Tt}^{\G}(\theta(\psi|T))$ (see,
e.g.,~\cite[Prop. 12.6]{digne-michel}), so
$$
\rho\otimes\psi=\sum_{(\Tt,\theta)\in\kappa}{
\frac{\eps_{\G}\eps_{\Tt}R_{\Tt}^{\G}(\theta(\psi| T))}
{\langle R_{\Tt}^{\G}(\theta),R_{\Tt}^{\G}(\theta)\rangle }
}.
$$
\par
Since the distinct Deligne-Lusztig characters are orthogonal, the
assumption $\rho\simeq \rho\otimes\psi$ implies that for any fixed
$(\Tt,\theta)\in\kappa$, the pair $(\Tt,\theta(\psi|T))$ is also in
the geometric conjugacy class $\kappa$. 
Consider then the translation of this condition using the bijection
between  geometric conjugacy classes of 
pairs $(\Tt,\theta)$ and $\Fp_q$-rational conjugacy classes of
semi-simple elements in $\G^*$, the dual group of $\G$ (see,
e.g.,~\cite[Prop. 13.12]{digne-michel}). Denote by $s$ the conjugacy
class corresponding to $(\Tt,\theta)$. The pair $(\Tt,\psi|T)$
corresponds to a \emph{central} conjugacy class $s'$, because $\psi|T$
is the restriction of a global character of $\G$
(see the proof of~\cite[Prop. 13.30]{digne-michel}; alternately, use
the fact that both global characters and central conjugacy classes are
characterized by being invariant under the action of the Weyl
group\footnote{\ Think of $\Tt$ in $GL(n)$ being the diagonal
  matrices, with the Weyl group $\mathfrak{S}_n$ permuting the
  diagonal components.}),
and the definition of the correspondance shows that
$(\Tt,\theta\psi|T)$ corresponds to the conjugacy class $ss'$ (which
is well-defined because $s'$ is central). The assumption that
$(\Tt,\theta)$ and $(\Tt,\theta\psi|T)$ are geometrically conjugate
therefore means $ss'=s$, i.e, $s'=1$, and clearly this means $\psi=1$,
as desired.
\end{proof}

\begin{remark}
Here is a mnemonic device to remember the
bounds for $A_{\infty}(G)$ in (\ref{eq-stronger})\footnote{\ Which 
  explains why it seemed to the author to be a reasonable statement to
  look for...}: among the representations of $G$, we have the principal
series $R(\theta)$, parametrized by the characters of a maximal split
torus, of which there are about $q^r$, and those share a common
maximal dimension $A$. Hence
$$
q^rA^2\asymp \sum_{\theta} {\dim(R(\theta))^2}\asymp |G|\sim q^d,
$$
so $A$ is of order $q^{(d-r)/2}$. In other words, we expect that in
the formula $\sum{\dim(\rho)^2}=|G|$, the principal series contributes
a positive proportion. 
\par
The bound for $A_1(G)$ is also intuitive:
there are roughly $q^r$ conjugacy classes, and as many
representations, and for a ``positive proportion'' of them, the degree
of the representation is of the maximal size given by
$A_{\infty}(G)$. 
\end{remark}


\section{Probabilistic sieves}
\label{sec-proba-sieve}

The introduction of a general measure space $(X,\mu)$ as component of
the siftable set may appear yo be an instance of overenthusiastic French
abstraction. However, we believe that the generality involved may be
useful and that it suggests new problems in a probabilistic setting.
\par
To start with a simple example, let
$\sieve=(\Zz,\{\text{primes}\}, \Zz\ra \Zz/\ell\Zz)$ be the
classical sieve setting. Consider now a probability space $(X,\Sigma,\proba)$ 
(i.e., $\proba$ is a 
probability measure on $X$ with respect to a $\sigma$-algebra
$\Sigma$), and let $F=N\,:\, X\ra \Zz$ be an
integer-valued random variable. Then the triple
$(X,\proba,N)$ is a siftable set, and given any sieving sets
$(\Omega_{\ell})$ and prime sieve support   $\Lb$, it is tautological that the
measure, or rather probability, of the associated sifted set in $X$ is
equal to   
$$
\proba(N\in S(\Zz,\Omega;\Lb^*))=\proba(\{\omega\in X\,\mid\,
N(\omega)\mods{\ell}\notin \Omega_{\ell},\text{ for all }
\ell\in\Lb^*\}).
$$
\par
In other words, the sieve bounds in that context can give estimates
for the probability that the values of some integer-valued random
variable satisfy any condition that can be described by sieving sets. 
\par
If we are given natural integer-valued random variables, this
probabilistic setting gives a
precise meaning to such notions as ``the probability that an integer
is squarefree''. If the distribution law of $N$ is uniform on an
interval $1\leq n\leq T$, and we let $T\ra +\infty$, this is just the
usual ``natural 
density''. 

\begin{example}
Let $N_{\lambda}$ be a random variable with a Poisson distribution of
parameter $\lambda$, i.e., we have
$$
\proba(N_{\lambda}=k)=e^{-\lambda}\frac{\lambda^k}{k!},
\quad \text{ for } k\geq 0.
$$
Then one can easily show, e.g., that the probability that
$N_{\lambda}$ is squarefree (excluding $0$) tends to $\pi^2/6$ as
$\lambda$ goes to $\infty$.
\end{example}

The following setting seems to  have some interest as a way to get
insight into properties of ``random'' integers $n\in \Zz$.
\par
Consider a simple
random walk $S_n$, $n\geq 0$, on $\Zz$, i.e., a sequence of random
variables $S_n$ on $X$ such that $S_0=0$ and $S_{n+1}=S_n+X_{n+1}$
with $(X_n)_{n\geq 1}$ a sequence of independent random variables with
Bernoulli distribution $\proba(X_n=\pm 1)=\demi$ (or one could take
general Bernoulli distributions $\proba(X_n=1)=p$, 
$\proba(X_n=-1)=q$, for some $p$, $q\in ]0,1[$ with $p+q=1$).
These variables $(S_n)$ give a natural sequence of siftable sets
$(X,\proba,S_n)$. It turns out to be quite easy to estimate the
corresponding sieve constants; here the dependency on the random
variable component of the siftable set is the most important, so we
denote $\Delta(S_n,\Lb)$ the sieve constant.

\begin{proposition}\label{pr-proba-sieve-constant}
Let $(S_n)$ be a simple random walk on $\Zz$. With notation as above,
we have 
$$
\Delta(S_n,\Lb)\leq 1+\Bigl|
\cos\Bigl(\frac{2\pi}{L^2}\Bigr)\Bigr|^n\sum_{m\in\Lb}{m},
$$
for $n\geq 1$ and for any sieve support $\Lb$ consisting entirely of
\emph{odd} squarefree integers $m\leq L$.
\end{proposition}

It is natural to exclude even integers, simply because $S_n\mods{2}$
is not equidistributed: more precisely, we have $\proba(S_n\text{ is
  even})=0$ or $1$ depending on whether $n$ itself is even or odd. In
probabilistic terms, the random walk is not aperiodic. The simplest
way to avoid this problem would be to assume that the increments $X_n$
have distribution
$$
\proba(X_n=\pm 1)=\proba(X_n=0)=\frac{1}{3}
$$
(i.e., at each step the walker may decide to remain still).
The reader will have no trouble adapting the arguments below to this
case, without parity restrictions.

\begin{proof}
We will estimate the ``exponential sums'', which in the current
context, using probabilistic notation $\expect(Y)=\int_X{YdP}$ for the
integral, are simply 
$$
W(a,b)=\expect\Bigl(e\Bigl(\frac{a_1S_n}{m_1}\Bigr)
e\Bigl(-\frac{a_2S_n}{m_2}\Bigr)\Bigr)
$$
for $m_1$, $m_2\in\Lb$, $a_i\in (\Zz/m_i\Zz)^{\times}$. Using the
expression $S_n=X_1+\cdots +X_n$ for $n\geq 1$, independence, and the
distribution of the $X_i$, we obtain straightforwardly
$$
W(a,b)=\expect\Bigl(e\Bigl(\frac{(a_1m_2-a_2m_1)X_1}{m_1m_2}\Bigr)\Bigr)^n
=\Bigl(\cos2\pi\frac{a_1m_2-a_2m_1}{m_1m_2}\Bigr)^n.
$$
The condition that $m_i$ are odd, and that $(a_i,m_i)=1$, imply that
$|W(a,b)|=1$ if and only if $a_1=a_2$ and $m_1=m_2$, and otherwise
$$
|W(a,b)|\leq \Bigl|\cos \frac{2\pi}{m_1m_2}\Bigr|^n.
$$
Hence the sieve constant is bounded by
$$
\Delta(S_n,\Lb)\leq \max_{m_1,a_1}\Bigl\{1+\sum_{m_2}
\sums_{a_2\mods{m_2}}{\Bigl|\cos \frac{2\pi}{m_1m_2}\Bigr|^n}\Bigr\}
\leq 1+\Bigl|
\cos\Bigl(\frac{2\pi}{L^2}\Bigr)\Bigr|^n\sum_{m\in\Lb}{m}.
$$
\end{proof}

\begin{corollary}
With notation as above, we have:
\par
\emph{(1)} For any sieving sets $\Omega_{\ell}\subset \Zz/\ell\Zz$ for
$\ell$ odd, $\ell\leq L$, and $L\geq 3$, we have
$$
\proba(S_n\in S(\Zz,\Omega;L))\leq
\Bigl(1+L^2\exp\Bigl(-\frac{n\pi^2}{L^4}\Bigr)\Bigr)
H^{-1}
$$
where
$$
H=\sumb_{\stacksum{m\leq L}{m\text{ odd}}}{\prod_{\ell\mid m}
{\frac{|\Omega_{\ell}|}{\ell-|\Omega_{\ell}|}}}.
$$
\par
\emph{(2)} Let $\eps>0$ be given, $\eps\leq 1/4$. For any odd $q\geq
1$, any $a$ coprime with $q$, we have 
$$
\proba(S_n\text{ is prime and } \equiv a\mods{q})\ll \frac{1}{\varphi(q)}
\frac{1}{\log n}
$$
if $n\geq 2$, $q\leq n^{1/4-\eps}$, the implied constant depending
only on $\eps$.
\end{corollary}

Note that (2) is Theorem~\ref{th-proba} in the introduction. 

\begin{proof}
For (1), we take $\Lb$ to be the set of odd squarefree numbers $\leq
L$ (so $\Lb^*$ is the set of odd primes $\leq L$), and then since
$\cos(x)\leq 1-x^2/4$ for $0\leq x\leq 2\pi/9$, 
the proposition gives 
$$
\Delta\leq 1+L^2\Bigl(1-\frac{\pi^2}{L^4}\Bigr)^n\leq 1+L^2
\exp\Bigl(-\frac{n\pi^2}{L^4}\Bigr),
$$
and the result is a mere restatement of the large sieve inequality.
\par
For (2), we have to change the sieve setting a little bit. Consider
the sieve setting above $\sieve$, 
\emph{except} that for primes $\ell\mid q$, we take $\rho_{\ell}$ to
be reduction modulo $\ell^{\nu(\ell)}$, where $\nu(\ell)$ is the
$\ell$-valuation of $q$. Take the siftable set $(X,\proba,S_n)$, and
the sieve support 
$$
\Lb=\{mm'\,\mid\, mm'\text{ squarefree, } (m,2q)=1,\ m\leq L/q\text{ and
} m'\mid q\},
$$
with  $\Lb^*$ still the set of odd primes $\leq L$.
\par
Proceeding as in the proof of
Proposition~\ref{pr-proba-sieve-constant}, the sieve 
constant is bounded straightforwardly by
$$
\Delta\leq 1+\Bigl|\cos
\frac{2\pi}{L^2}\Bigr|^n\sum_{\stacksum{m\leq
    L/q}{(m,2q)=1}}{\sum_{m'\mid q}{mq}}
\leq 1+\tau(q)q^{-1}L^2\exp\Bigl(-\frac{n\pi^2}{L^4}\Bigr),
$$
where $\tau(q)$ is the number of divisors of $q$. 
\par
Finally, take
$$
\Omega_{\ell}=\begin{cases}
\{0\}&\text{ if } \ell\nmid q,\\
\Zz/\ell^{\nu(\ell)}\Zz-\{a\}&\text{ if } \ell\mid q.
\end{cases}
$$
If $S_n$ is a prime number congruent to $a$ mod $q$, then we have $S_n\in
S(\Zz,\Omega;\Lb^*)$, hence
$$
\proba(S_n\text{ is a prime } \equiv a\mods{q})\leq
\proba(S_n\in S(\Zz,\Omega;\Lb^*))\leq
\Delta H^{-1}
$$
where
$$
H=\sumb_{\stacksum{m\leq L/q}{(m,2q)=1}}{
  \sumb_{m'\mid q}{\frac{\varphi^*(m')}{\varphi(m)}}},\quad\quad\text{
  with }\quad\quad
\varphi^*(n)=\prod_{\ell^{\nu}\mid\mid n}{(\ell^{\nu}-1)}.
$$
Now the desired estimate follows on taking $q\leq n^{1/4-\eps}$ and
$L=qn^{\eps}$, using the classical lower bound (see
e.g.~\cite{bombieri},~\cite[(6.82)]{ant}) 
$$
\sumb_{\stacksum{m\leq L/q}{(m,2q)=1}}{\frac{1}{\varphi(m)}}\geq
\frac{\varphi(q)}{2q} \log L/q\gg \frac{\varphi(q)}{q}\log n
$$
(the implied constant depending only on
$\eps$) together with the cute identity
$$
\sumb_{m'\mid q}{\varphi^*(m')}=q
$$
which is trivially verified by multiplicativity.
\end{proof}

\begin{remark}
(1) It is important to keep in mind that, by the Central Limit Theorem,
$|S_n|$ is usually of order of magnitude $\sqrt{n}$ (precisely,
$S_n/\sqrt{n}$ converges weakly to the normal distribution with
variance $1$ as $n\ra +\infty$). So 
the estimate $\Delta\leq 1+L^2\exp(-n\pi^2/L^4)$, which gives a
non-trivial result in applications as long as, roughly speaking, $L\leq
n^{1/4}/(\log n)^{1/4}$, compares well with the classical large
sieve for integers $n\leq N$, where $\Delta\leq N-1+L^2$, which is
non-trivial  for $L\leq \sqrt{N}$. 
\par
(2) The second part is  an analogue of the Brun-Titchmarsh
inequality, namely (in its original form)
$$
\pi(x;q,a)\ll \frac{1}{\varphi(q)}\frac{x}{\log x}
$$
for $x\geq 2$, $(a,q)=1$ and $q\leq x^{1-\eps}$, the implied constant
depending only on $\eps>0$. However, from the previous remark we see
that it is weaker than could be expected, namely $q\leq n^{1/4-\eps}$
would have to be replaced by $q\leq n^{1/2-\eps}$. Here we have
exploited the flexibility of the sieve setting and sieve support.
For a different use of this flexibility, see Section~\ref{sec-ls-ff};
we want to point out here that the possibility of using a careful
non-obvious choice of $\Lb$ was first exploited by Zywina in his
preprint~\cite{zywina}.
\par
It would be quite interesting to know if the extension to $q\leq
n^{1/2-\eps}$ holds.
The point is that if we try to adapt  the  classical
method, which is to sieve for those $k$, $1\leq k\leq x/q$, such that
$qk+a$ is prime, we are led to some interesting and non-obvious (for
the author) probabilistic issues; indeed, if $S_n\equiv a\mods{q}$,
the (random) integer $k$ such that $S_n=kq+a$ can be described as
follows: we have $k=T_N$ where $N$ is a random variable
$$
N=|\{m\leq n\,\mid\, S_m\equiv a\mods{q}\}|
$$
and $(T_i)$ is a random walk with initial
distribution given by
$$
\proba(T_0=0)=1-\frac{a}{q},\quad \proba(T_0=-1)=\frac{a}{q},
$$
and independent identically distributed increments $V_i=T_i-T_{i-1}$
such that
$$
\proba(V_i=0)=1-\frac{1}{q},\quad\quad\quad
\proba(V_i=\pm 1)=\frac{1}{2q}.
$$
\par
So what is needed is to perform sieve on the siftable
set $(\{S_n\equiv a\mods{q}\},\proba,T_N)$. Since the length $N$ of the
auxiliary walk is 
random, this requires some care, and we hope to come back to
this. Note at least that if look at the same problem with
$(\Omega,\proba,T_i)$ for a fixed $i$, then we easily get by sieving
$$
\proba(qT_i+a\text{ is prime})\ll \frac{1}{\varphi(q)}\frac{1}{\log i}
$$
for all $q\leq i^{1/2-\eps}$, $\eps>0$, the implied constant depending
only on $\eps$.
\par
(3) Obviously, it would be very interesting to derive lower bounds or
asymptotic formulas for $\proba(S_n\text{ is prime})$ for instance, and for
other analogues of classical problems of analytic number theory. Note
that it is tempting to attack the problems with ``local'' versions of
the Central Limit Theorem and summation by parts to reduce to the
purely arithmetic deterministic case. Problems where such a reduction
is not feasible would of course be more interesting.
\end{remark}

In the next section, we will give another
example of probabilistic sieve, similar in spirit to the above,
although the basic setting will be rather deeper, and the results are
not accessible to a simple summation by parts.
\par
Finally, we remark that this probabilistic point of view should not be
mistaken with ``probabilistic models'' of integers (or primes), such
as Cramer's model: the values of the random variables we have
discussed are perfectly 
genuine integers.\footnote{\ To give a caricatural example, if it were
  possible to show 
  that, for some sequence of random variables $N_n$ distributed on
  disjoint subsets of integers, the probability $\proba(N_n\text{ and }
  N_n+2 \text{ are both primes})$ is always strictly positive, then
  the twin-prime conjecture would follow.}

\section{Sieving in arithmetic groups}\label{sec-sieve-arith}

We now start discussing examples of sieve settings which seem to be
either new, or have only been approached very recently. The first
example concerns sieving for elements in an arithmetic group
$G$. There are actually a number of different types of siftable sets
that one may consider here.
\par
Maybe the most obvious idea for analytic number theorists is to take the
group sieve setting defined by
$$
\sieve=(SL(n,\Zz),\{\text{primes}\},G\ra
  SL(n,\Fp_{\ell}))
$$
(where the last reduction map is known to be surjective for all
$\ell$), and look at the siftable set 
$X\subset G$ which is the set of those matrices with norm bounded by some 
quantity $T$, with $F_x=x$ and counting measure. In other words,
instead of sieving 
integers, we want to sieve integral unimodular matrices.
Of course, $SL(n,\Zz)$ may be replaced with
other arithmetic groups,   even possibly with infinite-index
subgroups.  
\par
Here the equidistribution approach leads to hyperbolic lattice point
problems (in the case $n=2$), and generalizations of those for
$n\geq 3$. The issue of uniformity with respect to $q$ when taking
``congruence towers''   $\Gamma\cap\Gamma(q)$, where $\Gamma(q)$ is
the principal congruence   subgroup, is the main issue, compared
with the results available in the literature (e.g., the work of Duke,
Rudnick and Sarnak~\cite{drs} gives individual equidistribution, with
methods that may be amenable to uniform treatments, whereas more
recent ergodic-theoretic methods by Eskin, Mozes, McMullen, Shah and
other, see e.g.~\cite{ems}, seem to be more problematic in this
respect). 
\par
A tentative and very natural application is the natural fact that
``almost all'' unimodular matrices with norm $\leq T$ have
irreducible characteristic polynomial, with an estimate for the
number of exceptional matrices (this question was also recently
formulated by Rivin~\cite[Conj. 8]{rivin}, where it
is observed that the qualitative form of this statement is likely to
follow from the results of Duke, Rudnick and Sarnak). The case of integral
matrices with arbitrary determinant can be treated very quickly as a
simple consequence of the higher-dimensional 
large sieve as in~\cite{gallagher} (in other words, embed
invertible matrices in the additive group $G=M(n,\Zz)\simeq
\Zz^{n^2}$, and  use abelian harmonic analysis).
\par
The setting of arithmetic groups suggests 
other types of siftable sets, which are of a more combinatorial flavor,
and the ``probabilistic'' theme of
Section~\ref{sec-proba-sieve} is also a natural fit.\footnote{\ A useful
survey on combinatorial and geometry group theory is given in the book
of de la Harpe~\cite{harpe}, and a survey of random walks on groups is
that of Saloff-Coste~\cite{saloff-coste}.}  
Theorem~\ref{th-sln} gives some first results of this kind.
\par
Let $G$ be a finitely generated
group. Assuming a 
symmetric set of generators $S$ to be fixed (i.e., with $S^{-1}=S$),
three siftable sets $(X,\mu,F)$ of great interest arise naturally: 
\par
\noindent -- the set $X$ of elements $g\in
G$ with word-length metric $\ell_S(g)$ at most $N$, for some integer 
$N\geq 1$, i.e, the  set of those elements  $g\in G$ that can be
written as  
$$
g=s_1\cdots s_k
$$
with $k\leq N$, $s_i\in S$ for $1\leq i\leq k$. Here we take $F_x=x$ 
for $x\in X$, and of course $\mu$ is the counting measure.
\par
\noindent -- the set $W$ of words of length $N$ in the alphabet $S$,
for some integer $N\geq 1$, with $F_w$ the ``value'' in $G$ of the 
word $w\in W$, i.e., the image of $w$ by the natural (surjective)
homomorphism $F(S)\ra G$ from the free group generated by $S$ to
$G$. Again $\mu$ is the counting measure.
\par
\noindent -- as in Section~\ref{sec-proba-sieve}, we may consider a
probabilistic siftable set $(\Omega,\proba,X_N)$, where $\Omega$ is some 
probability space, $\proba$ the associated probability measure, and
$X_N=\xi_1\cdots \xi_N$, where $(\xi_k)$ is a $N$-uple of 
$S$-valued random variables. The simplest case is when $(\xi_k)$ is an
independent vector, and the distribution of each $\xi_k$ is uniform:
$\proba(\xi_k=s)=1/|S|$. In other words, $X_N$ is then the $N$-th
step in the simple left-invariant random walk on $G$ given by
$S$. If $G=\Zz$ and $S=\{\pm 1\}$, we considered this in
Section~\ref{sec-proba-sieve}. 

\begin{remark}\label{rm-equiv}
Note that the last two examples are in fact equivalent: we have
$$
\proba(X_N\in A)=\frac{1}{|W|}|\{w\in W\,\mid\, F_w\in A\}|
$$
for any subset $A\subset G$. (Since $|W|=|S|^N$, this
explains why  the two statements of Theorem~\ref{th-sln} are
equivalent). Although this reduces one particular probabilistic case
to a ``counting'' sieve,  we may indeed wish to vary the
distribution of the factors $\xi_k$ of the random walk, and doing so
would not in general lead to such a reduction; even when possible,
this may not be desirable, because it would involve rather artificial
constructs. For instance, another natural type of random walk is the
random  walk given by factors $\xi_k$ where  
$$
\proba(\xi_k=s)=\frac{1}{|S|+1}\quad\text{ for $s\in S$},\quad\quad
\proba(\xi_k=1)=\frac{1}{|S|+1}.
$$
\par
This is also equivalent to replacing $S$ by $S\cup \{1\}$ (if $1\notin
S$ at least), but the set of words $w$ where each component $w_i$ may
be the identity is not vary natural.
\end{remark}

We now provide a concrete example by proving Theorems~\ref{th-sln}
and~\ref{th-non-conj},
indeed in a slightly more general case. Let  $\Gg$ be  either $SL(n)$
or $Sp(2g)$ for some $n\geq 3$ or $g\geq 2$, let 
$G=\Gg(\Zz)$, and let $S$ be a symmetric set of generators for $G$
(for 
instance, the elementary matrices with $\pm 1$ off the diagonal are
generators of $SL(n,\Zz)$, see Remark~\ref{rm-shalom}).\footnote{\ We
will also comment briefly on what happens 
for $G=SL(2,\Zz)$.} We consider either the group sieve setting
$$
\sieve=(G,\{\text{primes}\},G\ra G_{\ell}) 
$$
where $G_{\ell}=\Gg(\Zz/\ell\Zz)$ is $SL(n,\Fp_{\ell})$ or
$Sp(2g,\Fp_{\ell})$, 
and the maps are simply reduction modulo $\ell$, or the induced
conjugacy sieve setting.
It is well-known that the reduction maps are onto for all $\ell$ (see
e.g.~\cite[Lemma 1.38]{shimura} for the case of $SL(n)$).
\par
We will look here at the second type of siftable set
$\siftable=(W,\mu,F)$, i.e., $W$ 
is the set of words of length $N$ in $S$, and $F_w$ is the
``value'' of a word $w$ in $G$.  Equivalently,
we consider the simple left-invariant random walk $(X_n)$.
In that case, the qualitative form of
Theorem~\ref{th-sln} was proved by Rivin~\cite{rivin}, and the latest
version of Rivin's preprint also discusses quantitative forms of
equidistribution in $G_{\ell}$, using Property (T) in a manner
analogous to what we do.  
\par
We will obtain a bound for the large sieve constant by appealing
to~(\ref{eq-max-cosets}) and its analogue for the group sieve setting,
estimating the exponential sums 
$W(\pi,\tau)$ or $W(\varphi_{\pi,e,f},\varphi_{\tau,e',f'})$
of~(\ref{eq-matrix-coeffs}).\footnote{\ Of course the
  equidistribution approach may 
  also be used, but it is less efficient and not really quicker or
  simpler.} 
The crucial ingredient is the so-called ``Property $(\tau)$''.

\begin{proposition}\label{pr-tau}
Let $\Gamma$ be a finitely generated group, $I$ an arbitrary index set
and $\rho_i\,:\, \Gamma\ra G_i$ for $i\in I$ a family of surjective
homomorphisms onto finite groups, such that  $\Gamma$ has Property $(\tau)$
with respect to the family $(\ker\rho_{i})$ of finite index
subgroups of $\Gamma$. 
\par
Let  $S=S^{-1}$ be a symmetric finite generating set of $\Gamma$, and
for $N\geq 1$, let $W=W_N$ denote the set of words of length $N$ in
the alphabet $S$, and let $F_w$ denote the value of the word $w$ in $\Gamma$.
Assume that there exists a word $r$ in the alphabet $S$ of \emph{odd}
length $c$ such that $F_r=\mathrm{Id}\in \Gamma$.
\par
Then there exists $\alpha>0$ such that for any $i\in I$, any
representation $\pi\,:\, \Gamma\ra GL(V)$ that factors through $G_i$
and does not contain the trivial representation, any vectors $e$, $f$
in the space of $\pi$, we have
\begin{equation}\label{eq-estim-matcoeff}
\Bigl|\sum_{w\in W}{\langle \pi(F_w)e,f\rangle}\Bigr|
\leq \|e\|\|f\| |W|^{1-\alpha},
\end{equation}
for $N\geq 1$, where $\langle \cdot,\cdot\rangle$ is a
$\Gamma$-invariant inner product on $V$, and hence
\begin{equation}\label{eq-estim-trace}
\Bigl|\sum_{w\in W}{\Tr\pi(F_w)}\Bigr|\leq (\dim \pi)|W|^{1-\alpha}.
\end{equation}
\par
The constant $\alpha$  depends only on
$\Gamma$, $|S|$, the $(\tau)$-constant for $(S,\Gamma,\ker\rho_i)$ and
the length $c$ of the relation $r$. 
\end{proposition}  

We will recall briefly the definition of Property $(\tau)$ and the
associated $(\tau)$-constant in the course of the
proof; see e.g.~\cite[\S 4.3]{lubotzky} or~\cite{lubotzky-zuk} for
more complete surveys. This should also be compared
with~\cite[Th. 6.15]{saloff-coste}. 

\begin{proof}
Let $i\in I$ and let $\pi$ be a representation that factors through
$G_i$ and does not contain the trivial
representation. Clearly~(\ref{eq-estim-trace}) follows
from~(\ref{eq-estim-matcoeff}) 
since the trace of a matrix is the sum of the diagonal matrix
coefficients in an orthonormal basis.
\par
Let
$$
M=\frac{1}{|S|}\sum_{s\in S}{\pi(s)},\quad\quad
M'=\mathrm{Id}-M,
$$
which are both self-adjoint elements of the endomorphism ring
$\End(V)$, since $S=S^{-1}$. We then find by definition
$$
\frac{1}{|W|}\sum_{w\in W}{\langle \pi(F_w)e,f\rangle}
=\langle M^Ne,f\rangle.
$$
\par
Let $\rho\geq 0$ be the spectral radius of $M$, or equivalently the
largest of absolute values of the eigenvalues of $M$, which are real
since $M$ is self-adjoint. Then by Cauchy's inequality we have
$$
|\langle M^Ne,f\rangle|\leq \|e\|\|f\| \rho^N,
$$
so that it only remains to prove that there exists $\delta>0$, 
independent of $i$ and $\pi$, such that $\rho\leq 1-\delta$.
\par
Clearly $\rho=\max(\rho_+,\rho_-)$, where $\rho_+\in\Rr$ (resp. $\rho_-$) is
the largest eigenvalue and $\rho_-$ is the opposite of the smallest
eigenvalue (if it is negative) and $0$ otherwise. We bound each
$\rho_{\pm}$ separately, proving $\rho_{\pm}\leq 1-\delta_{\pm}$ with
$\delta_{\pm}$ independent of $i$ and $\pi$.
\par
For $\rho_+$, it is equivalent (by the variational characterization of
the smallest eigenvalue) to
prove that there exists $\delta_+>0$, independent of $i$ and $\pi$,
such that
$$
\frac{\langle M'(v),v\rangle}{\langle v,v\rangle}\geq \delta_+
$$
for any non-zero vector $v\in V$. But a simple and familiar
computation yields 
$$
\frac{1}{|S|}
\sum_{s\in S}{\|\pi(s)v-v\|^2}
=
2\langle M'(v),v\rangle
$$
and therefore tautologically we have
\begin{equation}\label{eq-def-tau}
\frac{\langle M'(v),v\rangle}{\langle v,v\rangle}\geq 
\frac{1}{2|S|}\inf_{\varpi}\inf_{v\not=0}
{
\max_{s\in S}\frac{\|\varpi(s)v-v\|^2}{\|v\|^2},
}
\end{equation}
where $\varpi$ ranges over \emph{all} unitary representations of
$\Gamma$ that factor through some $\ker \rho_i$ and do not contain the
trivial representation (and $\|\cdot\|$ on the right-hand side is the
unitary norm for each such representation). But it is precisely the
content of Property $(\tau)$ 
for $\Gamma$ with respect to $(\ker\rho_i)$ that this triple extremum
is $>0$ (see e.g.~\cite[Def. 4.3.1]{lubotzky}). 
\par
So we come to $\rho_-$. Here a suitable lower-bound follows from 
Theorem~6.6 of~\cite{saloff-coste} (due to Diaconis, Saloff-Coste,
Stroock), using the fact that any eigenvalue of
$M$ is also an eigenvalue of $M_{reg}$, where $M_{reg}$ is the
analogue of $M$ for the regular representation of $\Gamma$ on
$L^2(\Gamma/\ker\rho_i)$.
\par
For completeness, we prove what is needed here, adapting the arguments
to the case of a general  representation. It suffices
to prove that there exists $\delta_->0$ independent of $i$ and $\pi$
such that
\begin{equation}\label{eq-deltaminus}
\frac{\langle M''(v),v\rangle}{\langle v,v\rangle}\geq \delta_-
\end{equation}
for all non-zero $v\in V$, where now $M''=\mathrm{Id}+M$. We have
$$
2\langle M''(v),v\rangle=
\frac{1}{|S|}\sum_{s\in S}{\|\pi(s)v+v\|^2}.
$$
\par
Now let $r=s_1\cdots s_{c}$ be a word of odd length $c$ in the
alphabet $S$ such that $r$ is trivial in $\Gamma$; denote
$$
r_k=s_1\cdots s_k\text{ for } 1\leq k\leq r,\quad r_0=1.
$$
\par
For $v\in V$, we can write 
$$
v=\frac{1}{2}\Bigl((v+\pi(s_1)v)-(\pi(s_1)v+\pi(s_1s_2v))+\cdots
+(\pi(s_1\cdots s_{c-1})v+\pi(1)v)\Bigr)
$$
(the odd length is used here), hence by Cauchy's inequality we get
$$
\|v\|^2\leq
\frac{c}{4}\sum_{i=0}^{c-1}{\|\pi(r_i)v+\pi(r_{i}s_{i+1})v\|^2}
=\frac{c}{4}\sum_{i=0}^{c-1}{\|v+\pi(s_{i+1})v\|^2}
$$
(the representation is unitary). By positivity, since at worst all
$s_i$ are equal to the same generator in $S$, we get
\begin{equation}\label{eq-odd-length-bound}
\|v\|^2\leq \frac{c^2}{4}\sum_{s\in S}{\|\pi(s)v+v\|^2}=
\frac{c^2|S|}{2}\langle M''(v),v\rangle,
\end{equation}
which implies~(\ref{eq-deltaminus}) with
$\delta_-=\frac{2}{c^2|S|}>0$. 
\end{proof}

\begin{remark}
The odd-looking assumption on the existence of $r$ is indeed necessary
for such a general statement, because of periodicity issues.  Namely, if
(and in fact only if) all $S$-relations in $\Gamma$ are of even
length, the Cayley graph\footnote{\ Recall
  $C(\Gamma,S)$ has vertex set $\Gamma$ and as many edges from $g_1$
  to $g_2$ as there are elements $s\in S$ such that  $g_2=g_1 s$; this
  allows both loops and multiple edges, and those will occur if $1\in
  S$ or, in the Cayley graphs of quotients of $\Gamma$, if two
  generators have the same image.} 
$C(\Gamma,s)$ of 
$\Gamma$ with respect to $S$ is bipartite\footnote{\ I.e., the vertex 
  set $\Gamma$ is partitioned in two pieces $\Gamma_{\pm}$ and edges
  always go from one piece to another.}, and so are its finite quotients
$C(\Gamma/\ker\rho_i,S)$. In that case, it is well-known and easy to
see that $-1$ is an eigenvalue of 
$M_{reg}$ (the operator $M$ for the regular representation; take the
function such that $f(x)$ equals $\pm 1$ depending 
on whether the point is at even or odd distance from the origin) and
the argument above fails.  Alternately, this can be seen directly
with the exponential sums: the relations being of even length implies
that there is a 
well-defined surjective homomorphism $\eps\,:\,\Gamma\ra \Zz/2\Zz$ with
$\eps(s)=-1$ for $s\in S$. Viewing $\eps$ as a representation
$\Gamma\ra \{\pm 1\}\subset \Cc^{\times}$, we have
$$
\sum_{w\in W}{\eps(F_w)}=\begin{cases}
|W|&\text{ if $N$ is even}\\
-|W|&\text{ if $N$ is odd}.
\end{cases}
$$
We will describe an example of this for $\Gamma=SL(2,\Zz)$ below.
\par
The simplest way of ensuring that $r$ exists is to assume that
$1\in S$; geometrically, this means each vertex of the Cayley
graph
has a self-loop, and probabilistically, this means that one considers
a ``lazy'' random walk on the Cayley graph, with probability $1/|S|$ of
staying at the given element. 
\par
In fact, if we consider the effect of replacing $S$ by $S'=S\cup \{1\}$
(in the case where $1\notin S$), we have
$$
M_{S'}=\Bigl(1-\frac{1}{|S'|}\Bigr)M_S+\frac{1}{|S'|},
$$
with obvious notation, and so we obtain
$$
\rho_-\geq -1+\frac{2}{|S'|}
$$
directly (which is the same lower bound as the one we proved, in the
case $c=1$).
\end{remark}

With the  estimate of Proposition~\ref{pr-tau}, we can perform some
sieve. 

\begin{theorem}\label{th-sln-bis}
Let $\Gg=SL(n)$, $n\geq 3$, or $Sp(2g)$, $g\geq 2$,  be as before,
$G=\Gg(\Zz)$, and let  
$$
\sieve=(G, \{\text{primes}\},G\ra G_{\ell}=\Gg(\Zz/\ell\Zz))
$$
be the group sieve setting. Let $S=S^{-1}$ be a symmetric generating
set for  $G$, $(W,F)$ the siftable set of random 
products of length $N$ of elements of $S$. 
\par
\emph{(1)} For any sieve support
$\Lb$, the large sieve constant for the induced conjugacy sieve
satisfies 
\begin{equation}\label{eq-random-first}
\Delta(W,\Lb)\leq |W|+|W|^{1-\alpha}R(\Lb),
\end{equation}
where $\alpha>0$ is a constant depending only on $\Gg$ and $S$
and\ \footnote{\ With notation as in Section~\ref{sec-reprs}.}
$$
R(\Lb)=\max_{m\in \Lb}\Bigl\{A_{\infty}(G_m)\Bigr\}\times 
\sum_{n\in \Lb}{A_1(G_n)}.
$$
\par
\emph{(2)} There exists $\eta>0$ such that
\begin{equation}
  \label{eq-random-sec}
|\{w\in W\,\mid\, \det(F_w-T)\in \Zz[T]\text{ is reducible }\}|
\ll |W|^{1-\eta}
\end{equation}
where $\eta$ and the implied constant depend only on $\Gg$ and $S$.  
\par
\emph{(3)}  For any sieve support
$\Lb$, the large sieve constant for the group sieve
satisfies 
\begin{equation}\label{eq-random-third}
\Delta(W,\Lb)\leq |W|+|W|^{1-\alpha}\tilde{R}(\Lb),
\end{equation}
where $\alpha>0$ is as above and
$$
\tilde{R}(\Lb)=\max_{m\in \Lb}\Bigl\{\sqrt{A_{\infty}(G_m)}\Bigr\}\times 
\sum_{n\in \Lb}{A_{5/2}(G_n)^{5/2}}.
$$
\par
\emph{(4)}  There exists $\beta>0$ such that
\begin{equation}
  \label{eq-random-fourth}
|\{w\in W\,\mid\, \text{one entry of $F_w$  is a square }\}|
\ll |W|^{1-\beta}
\end{equation}
where $\beta$ and the implied constant depend only on $\Gg$ and $S$. 
\end{theorem}

It is clear that the fourth  part implies  Theorem~\ref{th-non-conj}.

\begin{lemma}\label{lm-4}
Let  $\Gg$ be as above.
\par
\emph{(1)} Property $(\tau)$ holds for the group $G=\Gg(\Zz)$ with
respect to the family of congruence subgroups $(\ker(G\ra
\Gg(\Zz/d\Zz)))_{d\geq 1}$.
\par
\emph{(2)} For any symmetric generating set $S=S^{-1}$, there exists
an $S$-relation of odd length.
\end{lemma}

\begin{proof}
(1) This is well-known; in fact, the group $G$ is a lattice in a
semisimple real Lie group with $\Rr$-rank $\geq 2$, and hence
it satisfies the stronger 
Property $(T)$ of Kazhdan, which means that in~(\ref{eq-def-tau}), the
infimum may be taken on \emph{all} unitary representations of $G$ not
containing the trivial representation and remains $>0$ (see,
e.g.,~\cite[Cor. 3.5]{harpe-valette},~\cite[Prop.3.2.3, Ex. 3.2.4, \S 
4.4]{lubotzky})).
\par
(2) If all $S$-relations are of even length, the homomorphism 
$$
F(S)\ra \{\pm 1\}
$$
defined by $s\mapsto -1$ induces a non-trivial homomorphism $G\ra
\{\pm 1\}$. However, there is no such homomorphism for the groups
under consideration (e.g., because its kernel $H$ will be a finite
index normal subgroup, hence by the Congruence 
Subgroup Property, due to Mennicke and Bass-Lazard-Serre in this case,
see~\cite[p. 64]{bass-milnor-serre} for references,
will factor through a principal congruence subgroup $\ker
(G\ra \Gg(\Zz/d\Zz))$ for some integer  $d\geq 1$, defining 
a non-trivial homomorphism\footnote{\ Here we use the fact that 
  $G\ra\Gg(\Zz/d\Zz)$ is surjective, see the first line of the next
  proof.} $\Gg(\Zz/d\Zz)\ra \{\pm
1\}$, which  is impossible since $\Gg(\Zz/d\Zz)$ is its own
commutator group).
\end{proof}

\begin{proof}[Proof of Theorem~\ref{th-sln-bis}]
(1) Let $m$, $n\in \Lb$, $\pi$, $\tau\in \Pi_m^*$, $\Pi_n^*$
respectively. Since the maps  $G\ra \Gg(\Zz/d\Zz)$ are onto for all
$d$ (e.g., because the family $(\rho_{d})$ is linearly
disjoint in the sense of Definition~\ref{def-chinese}, by Goursat's
lemma, as in~\cite[Prop. 5.1]{chavdarov}),  
we have in fact $G_{[m,n]}=\Gg(\Zz/[m,n]\Zz)$.  
\par
By Lemma~\ref{lm-ortho}, the representation
$[\pi,\bar{\tau}]$ of $G_{[m,n]}$ defined 
in~(\ref{eq-lcm-reps}) contains the identity
representation if and only if $(m,\pi)=(n,\tau)$, and then contains
it with multiplicity one. 
Let $[\pi,\bar{\tau}]_0$ denote the orthogonal of
the trivial component in the second case, and
$[\pi,\bar{\tau}]_0=[\pi,\bar{\tau}]$ otherwise.
\par
We can now appeal to
Proposition~\ref{pr-tau} applied to 
the representation 
$[\pi,\bar{\tau}]_0\circ\rho_{[m,n]}$  of $G$, 
using the family 
$(\rho_d\,:\, G\ra \Gg(\Zz/d\Zz))$ of 
congruence subgroups (since $\ker [\pi,\bar{\tau}]_0\circ
\rho_{[m,n]}\supset G_{[m,n]}$). The previous lemma ensures that all
required 
assumptions on $S$ and this family are valid, and
by~(\ref{eq-estim-trace}), the conclusion is 
the estimate
$$
\Bigl|W(\pi,\tau)-\delta(\pi,\tau)|W|\Bigr|
\leq (\dim \pi)(\dim\tau)|W|^{1-\alpha}
$$
for the exponential sum~(\ref{eq-sumexp}), where $\alpha$ depends only
on $\Gg$, $S$ and the relevant $(\tau)$ 
or $(T)$ constant.
\par
By Proposition~\ref{pr-w-to-sieve}, we obtain
$$
\Delta(W,\Lb)\leq |W|+|W|^{1-\alpha}\max_{m\in\Lb}A_{\infty}(G_m)
\sum_{n\in \Lb}{A_1(G_n)},
$$
as stated.
\par
(3) This is exactly similar, except that now we
use the basis of matrix coefficients for the group sieve setting, and
correspondingly we appeal to~(\ref{eq-estim-matcoeff}) and the fact
(see the final paragraphs of Section~\ref{sec-spec-group}) that the
sums $W(\varphi_{\pi,e,f},\varphi_{\tau,e',f'})$ are (up to the factor
$\sqrt{(\dim\pi)(\dim\tau)}$) of the type considered
in~(\ref{eq-estim-matcoeff}). 
\par
In the case where $[\pi,\tau]$ contains the trivial representation
(i.e., if $(m,\pi)=(n,\tau)$), we also use the fact that, when identified
with $\End(V_{\pi})$, the one-dimensional space of invariant vectors
in $\pi\otimes\bar{\pi}=[\pi,\pi]$ is spanned by homotheties  and the
orthogonal projection of a linear map $u\in \End(V_{\pi})$  is 
multiplication by $\Tr(u)/\sqrt{\dim\pi}$ (this is a corollary of the
orthogonality relations; note that $\|\mathrm{Id}\|^2=\dim\pi$). This
means that for a rank $1$ linear map of 
the form $u=e\otimes \bar{e}'$ (where $e$ is in the space of $\pi$,
and $e'$ in that of the 
contragredient), the projection is the multiplication by $\langle
e,e'\rangle/\sqrt{\dim\pi}$. Since the vectors are part of an
orthonormal basis, we get
\begin{align*}
\langle (\pi\otimes\bar{\pi})(e\otimes e'),f\otimes f'\rangle&=
\frac{\langle e,e'\rangle \langle f,f'\rangle}{\dim\pi}+
\langle [\pi,\pi]_0(e\otimes e'),f\otimes f'\rangle\\
&=
\frac{\delta((e,f),(e',f'))}{\dim\pi}
+\langle [\pi,\pi]_0(e\otimes e'),f\otimes f'\rangle.
\end{align*}
\par
Altogether, we obtain
$$
\Bigl|
W(\varphi_{\pi,e,f},\varphi_{\tau,e',f'})-\delta((\pi,e,f),(\tau,e',f'))
|W|\Bigr|
\leq \sqrt{(\dim \pi)(\dim\tau)}|W|^{1-\alpha}
$$
and hence
\begin{align*}
\Delta(W,\Lb)&\leq |W|+|W|^{1-\alpha}\max_{m,\pi,e,f}
\sqrt{\dim\pi}
\sum_{n\in\Lb}\sum_{\tau,e',f'}{\sqrt{\dim\tau}}\\
&\leq  |W|+|W|^{1-\alpha}\max_{m\in\Lb}{\sqrt{A_{\infty}(G_m)}}
\sum_{n\in\Lb}\sum_{\tau}{(\dim\tau)^{5/2}}\\
&= |W|+|W|^{1-\alpha}\max_{m\in\Lb}{\sqrt{A_{\infty}(G_m)}}
\sum_{n\in\Lb}{A_{5/2}(G_n)^{5/2}}.
\end{align*}
\par
(2) To obtain~(\ref{eq-random-sec}), we apply the large sieve
inequality for group sieves of Proposition~\ref{pr-conjug-sieve},
using~(\ref{eq-random-first}). This is completely standard; without
trying to get the sharpest 
result (see 
Section~\ref{sec-ls-ff} for more refined arguments in a similar
end-game), we select  the prime sieve  
support $\Lb^*=\{\ell\leq L\}$ for some $L\geq 2$, and take
$\Lb=\Lb^*$ (pedantically, the singletons of elements of
$\Lb^*$...). Letting $d=n^2-1$, 
$r=n-1$ (for $\Gg=SL(n)$) or $d=2g^2+g$, $r=g$ (for $\G=Sp(2g)$), we
have (using Lemma~\ref{lm-sl-sp} of Section~\ref{sec-reprs})
$$
R(\Lb)\ll L^{d+1},
$$
for $L\geq 2$, the implied constant depending only on $\Gg$.
\par
We take for sieving sets the conjugacy classes in the set
$\Omega_{\ell}\subset \Gg(\Fp_{\ell})$ of matrices with irreducible
characteristic polynomial in $\Fp_{\ell}[T]$. 
From (1) of Proposition~\ref{pr-local} in Appendix~B, we obtain
\begin{equation}\label{eq-lower-omega}
\frac{|\Omega_{\ell}|}{|G_{\ell}|}\gg 1
\end{equation}
for $\ell\geq 3$, where the implied constant depends on $\Gg$ (compare
with~\cite[\S 3]{chavdarov},~\cite[Lemma 7.2]{sieve}). 
\par
Since
those $w\in W$ for which $\det(F_w-T)$ is reducible are contained
in the sifted set $S(W;\Omega,\Lb^*)$, we have by
Proposition~\ref{pr-conjug-sieve}
$$
|\{w\in W\,\mid\, \det(F_w-T)\in \Zz[T]\text{ is reducible }\}|
\leq  \Delta H^{-1}\ll (|W|+|W|^{1-\alpha}L^{d+1})H^{-1},
$$
where $H\gg \pi(L)$ by~(\ref{eq-lower-omega}). Taking
$L=|W|^{\alpha/(d+1)}$, we get the bound stated.
\par
(4) Clearly, it suffices to prove the estimate for the number of $w\in
W$ for which the $(i,j)$-th component of $F_w$ is a square, where $i$
and $j$ are fixed integers from $1$ to $n$ or $2g$ in the $SL(n)$ and
$Sp(2g)$ cases respectively.  The principle is similar, using (2) to
estimate the large sieve constant for the sieve where
$\Lb^*=\{\ell\leq L\}$, $\Lb=\Lb^*$, with
$$
\Omega_{\ell}=\{g=(g_{\alpha,\beta})\in \Gg(\Fp_{\ell})\,\mid\,
g_{i,j}\text{ is not a square in $\Fp_{\ell}$}\}.
$$
\par
We get by Lemma~\ref{lm-sl-sp} the bound
$$
\tilde{R}(\Lb)\ll L^{1+(3d-r)/2}
$$
for $L\geq 2$, where the implied constant depends only on $\Gg$.
\par
Next by (2) of Proposition~\ref{pr-local}, we have
$$
\frac{|\Omega_{\ell}|}{|G_{\ell}|}\gg 1
$$
for $L\geq 3$ (for $L=2$, the left-hand side may vanish for
$SL(2)$), where the implied constant depends only on 
$\Gg$. (The proof in Appendix B uses the Riemann Hypothesis over
finite fields; the reader may find it interesting to see whether a
more elementary argument may be found).
\par
Hence the sieve bound is
$$
|\{w\in W\,\mid\, \text{ the $(i,j)$-th entry of $F_w$ is a square}\}
|\leq (|W|+\tilde{R}(L))H^{-1}
$$
with $H\gg \pi(L)$ for $L\geq 3$, the implied constant depending on
$\Gg$. We take $L=|W|^{\alpha/(1+(3d-r)/2)}$ if this is
$\geq 3$ and then obtain~(\ref{eq-random-fourth}). To deal with
those $N$ for which this $L$ is $<3$, we just enlarge the implied
constant in~(\ref{eq-random-fourth}).
\end{proof}

From part (2) of this theorem, we can easily deduce
Theorem~\ref{th-sln}.

\begin{corollary}\label{cor-almost-surely}
Let $\Gg=SL(n)$, $n\geq 2$, or $Sp(2g)$, $g\geq 1$, let
$G=\Gg(\Zz)$ and let  $S=S^{-1}$ be a symmetric generating set of
$G$. Let $(X_k)$  be
the associated simple left-invariant random walk on $G$. Then almost
surely there exist only finitely many $k$ such that $\det(X_k-T)$ is a
reducible polynomial.
\end{corollary}

Part of the point of this statement is that it
\emph{requires} some quantitative estimate for the probability that
$X_k$ has reducible characteristic polynomial.

\begin{proof}
For $n\geq 3$ (resp. $g\geq 2$), it suffices to apply the ``easy''
Borel-Cantelli lemma\footnote{\ If $A_n$ are events in a probability
  space such that the series $\sum{\proba(A_n)}$ converges, then almost
surely $\omega$ belongs to only finitely many $A_n$.}, since the 
estimate (2) above for $N=k$ is equivalent with
$$
\proba(\det(X_k-T)\text{ is reducible})\ll \exp(-\alpha k)
$$
with $\alpha=\eta \log |S|>0$, 
and this shows that the series
$$
\sum_{k\geq 0}{\proba(\det(X_k-T)\text{ is reducible})}
$$
converges. From the weaker bound~(\ref{eq-sldeux}) in
Remark~\ref{rm-sldeux} below, we see that this series remains
convergent for $SL(2,\Zz)=Sp(2,\Zz)$.
\end{proof}

The next corollary is a geometric application which answers a question
of Maher~\cite[Question 1.3]{maher}, and was suggested by Rivin's
paper~\cite{rivin}. See~\cite{ivanov} for a survey of the mapping
class group of surfaces,~\cite[Exp. 1, 9]{flp} for information on
pseudo-Anosov diffeomorphisms of surfaces.

\begin{corollary}\label{cor-anosov}
Let $G$ be the mapping class group of a closed orientable surface of
genus $g\geq 1$, let $S$ be a finite symmetric generating set of
$G$ and let $(X_k)$, $k\geq 1$, be the simple left-invariant random walk
on $G$. Then the set $X\subset G$ of non-pseudo-Anosov elements is
transient for this random walk.
\end{corollary}

\begin{proof}
We follow the arguments of Rivin. First of all, the mapping class
group $G$ may be defined as the group 
of diffeomorphisms of a fixed compact connected surface
$\Sigma_g$ of genus $g$ preserving the orientation, up to isotopy
(i.e., homotopy in the diffeomorphism group). The main point is that
the induced action on the integral homology $H_1(\Sigma_g,\Zz)$, which
preserves the intersection pairing, yields a surjective map
$$
\rho\,:\, G\ra Sp(2g,\Zz).
$$
Let $S$ be a generating set as above\footnote{The existence of such
  finite generating set is not obvious, of course, and is known as the
  Dehn-Lickorish Theorem; see e.g.~\cite[Th. 4.2.D]{ivanov}.}, and let
$S'=\rho(S)$, a finite symmetric generating set for $Sp(2g,\Zz)$.
The image $Y_k=\rho(X_k)$ of the random walk on $G$ is a random walk
on $Sp(2g,\Zz)$. 
\par
Note that  the steps $\xi_k$ are independent and
identically distributed, but is not
necessarily true that  each  $\xi_k$ is uniformly distributed on $S'$,
which means that we are not exactly in the setting of
Theorem~\ref{th-sln-bis}. 
However, we can easily prove the analogue of
Proposition~\ref{pr-tau} for any random walk on $Sp(2g,\Zz)$ defined
by identically distributed independent steps $\xi_k$ with the property
that $\proba(\xi_k=s')=p(s')>0$ for all $s'\in S'$, simply by
replacing the self-adjoint operator $M$ with
$$
M=\sum_{s\in S'}{p(s')\pi(s')},
$$
and using the identities
$$
\sum_{s'\in S'}{p(s')\|\pi(s')v\pm v\|^2}=2\langle
(\mathrm{Id}\pm M)v,v\rangle 
$$
to obtain the bounds
$$
\frac{\langle
(\mathrm{Id}- M)v,v\rangle }{\langle v,v\rangle}
\geq \frac{\min p(s')}{2}\inf_{\varpi}\inf_{v\not=0}
{
\max_{s\in S}\frac{\|\varpi(s)v-v\|^2}{\|v\|^2},
}
$$
and
$$
\|v\|^2\leq \frac{c^2}{4}\frac{1}{\min p(s')}\sum_{s'\in S'}
{
p(s')\|v+\pi(s')v\|^2
}
=
\frac{c^2}{2\min p(s')}\langle
(\mathrm{Id}+M)v,v
\rangle
$$
analogues of~(\ref{eq-def-tau}) and~(\ref{eq-odd-length-bound}). From
this, sieve bounds for the random walk $(Y_k)$ follow, comparable to
those for the simple random walk.
\par
Now, we need only use the fact (the ``homological criterion for
pseudo-Anosov diffeomorphism'') that it \emph{suffices} that the
following three conditions on the characteristic polynomial
$P=\det(T-\rho(X_k))$ hold for $X_k$ to be pseudo-Anosov:
\par
-- $P$ is irreducible;
\par
-- there is  no root of unity which is a zero of $P$;
\par
-- there is no $d\geq 2$ and polynomial $Q$ such that $P(X)=Q(X^d)$.
\par
Accordingly we have
$$
\proba(X_k\text{ is not pseudo-Anosov})\leq p_1+p_2+p_3
$$
where $p_1$, $p_2$, $p_3$ are the probabilities that
$\det(T-\rho(X_k))$ satisfy those three conditions. Assume first
$g\geq 2$. Then, according to (2) of Theorem~\ref{th-sln-bis} (adapted
to a non-simple random walk), there exists $\alpha>0$ such that
$$
p_1\ll \exp(-\alpha k) 
$$
for $k\geq 1$. To estimate $p_2$ and $p_3$ we can use simpler sieves
to obtain comparable bounds. For $p_2$, since $P$ is an integral
polynomial of degree $2g$ and hence may have only finitely many roots
of unity as zeros, we need only estimate the sifted sets for the
sifting sets 
$$
\Omega_{\ell}=\{g\in Sp(2g,\Fp_{\ell})\,\mid\, (\Phi_d\mods{\ell})\nmid
\det(T-g)\} 
$$
where $\Phi_d\in\Zz[X]$ is the $d$-th cyclotomic polynomial for some fixed
$d\geq 1$.  It is clear (by the same local density arguments of
Appendix~B that were already used) that $|\Omega_{\ell}|\gg
|Sp(2g,\Fp_{\ell})|$, and hence the sieve again yields
$p_2\ll \exp(-\alpha k)$ for $k\geq 1$.
\par
For $p_3$, we consider similarly
$$
\Omega'_{\ell}=\{g\in Sp(2g,\Fp_{\ell})\,\mid\, \det(T-g)\text{ is not
  of the form } Q(X^d)\} 
$$
for some fixed $d\geq 2$. We also trivially have $|\Omega'_{\ell}|\gg 
|Sp(2g,\Fp_{\ell})|$, and $p_3\ll \exp(-\alpha k)$.
\par 
Now we conclude that $\proba(X_k\text{ is not pseudo-Anosov})\ll
\exp(-\alpha k)$, and we can again apply the Borel-Cantelli lemma.
\end{proof}

\begin{remark}\label{rm-anosov}
Maher~\cite{maher} proved that the probability that $X_k$ is
pseudo-Anosov tends to 
$1$ as $k\ra +\infty$ using rather more of the geometry and structure
of the mapping class group, and asked about the possible
transience of the set of non-pseudo-Anosov elements. However, his
methods are also more general, and work for random walks on any subgroup
of $G$ that is not ``too small'' in some sense. It should be
emphasized that his condition encompasses groups which seem utterly
inapproachable by sieve as above, in particular the so-called Torelli
group which is the kernel of the homology action $\rho$ above. It may
seem surprising  that pseudo-Anosov should exist in this subgroup, but 
Maher's result shows that they remain ``generic''
(see~\cite[p. 250]{flp} for a construction which gives some examples,
and the observation that Nielsen had conjectured they did not exist).
It would be interesting to know if the random walk on the Torelli
group is still transient on the set of pseudo-Anosov elements.
\end{remark}

\begin{remark}\label{rm-shalom}
In the most classical sieves, estimating either the analogue of
$R(\Lb)$ or $H$ is not a significant part of the work, the latter
because once $\Omega_{\ell}$ is known, which is usually not a problem
there, it boils down to estimates for sums of multiplicative
functions, which are well understood.
\par
The results we have proved, and an
examination of Appendix~B show that 
when performing a sieve in some group 
setting, sharp estimates for
$R(\Lb)$ or for $H$ involve deeper tools. For the large sieve
constant, this 
involves the representation theory of the group in non-trivial ways.
For $H$, the issue of estimating $|\Omega_{\ell}|$ may quickly become
a difficult counting problem over finite fields. It is not hard to
envision situations where the 
full force of Deligne's work on exponential sums over finite fields
becomes really crucial, and not merely a convenience.
\par
Note that the use of the sharp upper bounds of
Proposition~\ref{pr-degrees} in the proof of Theorem~\ref{th-sln-bis}
is not necessary if one wishes merely to find a bound for the large
sieve constant of the type $|W|+|W|^{1-\alpha}L^{A}$ for some
$\alpha>0$ and $A\geq 0$: trivial bounds for $A_p(G)$ are
sufficient. 
\par
If no exact value of the $(\tau)$-constant for $G=\Gg(\Zz)$
and $S$ is known, the value of $\alpha$ coming
from Proposition~\ref{pr-tau}  is not explicit, so
knowing a specific value of $A$ is not particularly rewarding.
However, in some cases explicit Kazhdan, hence $(\tau)$, constants are
known for the groups we are considering. The question of such explicit
bounds was first 
raised by Serre, de la Harpe and Valette; the arguments above show
it is indeed a very natural question with concrete applications, such
as explicit sieve bounds.  
The  first results for arithmetic groups are due to Burger for
$SL(3,\Zz)$ (see~\cite[Appendice]{harpe-valette}).
\par 
To give an idea, we quote a result of Kassabov~\cite{kassabov},
improving an earlier one of Shalom~\cite[Cor. 1]{shalom}:
let $G=SL(n,\Zz)$ with $n\geq 3$, and let $S$ be the symmetric 
generating set (of $2(n^2-n)$ elements) of elementary matrices
$E_{i,j}(\pm 1)$ with $\pm
1$ in the $(i,j)$-th entry. Then, for any unitary representation $\pi$
of $G$ not containing the trivial representation, and  any non-zero
vector $v$ in the space of $\pi$, there exists $s\in S$ such that
$\|\pi(s)v-v\|\geq \eps_n \|v\|$, with $\eps_n=(42\sqrt{n}+920)^{-1}$.
\par
The standard commutator relation 
$$
E_{1,2}(1)[E_{1,3}(1),E_{3,2}(1)]^{-1}=1
$$
(which uses that $n\geq 3$...)
shows that there are relations  of odd length $\leq 5$ in terms of
$S$. Looking at the proof of Proposition~\ref{pr-tau}, we see that we
can take $\delta_+=\demi \eps_n^2 |S|^{-1}$ and
$\delta_-=\frac{1}{25(n^2-n)}>\delta_+$. This means that for this
particular generating set, $\alpha$ in
Theorem~\ref{th-sln-bis} can be taken to satisfy
$$
\alpha=-\frac{\log(1-\eps_n^2|S|^{-1})}{\log |S|}\geq 
\frac{1}{8(n^2-n)(21\sqrt{n}+460)^2\log( 2(n^2-n))}
$$
for $n\geq 3$. So we have
$$
|\{w\in W\,\mid\, \det(F_w-T)\text{ is reducible}\}|
\ll |W|^{1-\eta}
$$
for $N\geq 1$, the implied constant depending on $n$, with 
$\eta$ given by
$$
\eta=\frac{\alpha}{n^2}\geq \frac{1}{8n^3(n-1)(21\sqrt{n}+460)^2\log(
2(n^2-n))}\gg \frac{1}{n^5\log n}. 
$$
\par
Coming back to the probabilistic interpretation (which is more suited
to what follows), this means in particular that if $k$ is of order of
magnitude larger 
than $n^5\log n$, the probability that $\det(X_k-T)$ is irreducible
becomes close to $1$. It would be interesting to have
a more precise knowledge of this ``transition time''
$$
\tau_n=\min\{k\geq 1\,\mid\, \det(X_k-T)\text{ is irreducible}\},
$$
(which, of course, depends also on $S$).
\par
Note that, at
the very least, with this particular generating set, $\det(X_k-T)$ is
reducible for $k\leq t_n$ where
$t_n$ is the first time when all basis vectors have been moved at
least once. Since multiplying by $\xi_{k}$ means moving one of the $n$
basis vectors chosen uniformly, $t_n$ is the stopping time for the
``coupon collector problem''. Besides the obvious bound $t_n\geq
n$, it is well-known (see,
e.g.,~\cite[IX.3.d]{feller}) that
\begin{gather*}
\expect(t_n)=n(\log n+\gamma)+O(1),\text{ for $n\geq 1$},\quad
\variance(t_n)\sim \zeta(2)n^2\text{ as } n\ra
+\infty.
\end{gather*}
\par
The gap between upper and lower bounds for $\tau_n$ is quite large,
and numerical 
experiments strongly suggest that the lower bound is closer to the
truth (in fact, it suggests that $\expect(\tau)$ might be $\sim
c\expect(t_n)$ for some constant $c>1$ as $n\ra +\infty$). In terms of
possible improvements, it is interesting to note that the order of
magnitude of Kassabov's estimate of the Kazhdan
constant $\eps_n$ for this generating set $S$ is optimal, since Zuk has
pointed out that it must be $\geq \sqrt{2/n}$ 
at least (see~\cite[p. 149]{shalom}). 
\end{remark}

\begin{remark}\label{rm-sldeux}
If $G=SL(2,\Zz)$, although $G$ does not have Property $(T)$, it is
still true that Property $(\tau)$ holds for the congruence subgroups
$\ker(SL(2,\Zz)\ra SL(2,\Zz/d\Zz))$, by Selberg's $\lambda_1\geq 1/4$
theorem on the smallest eigenvalue of the hyperbolic laplacian acting
on congruence subgroups of $SL(2,\Zz)$. However, the second condition
of Lemma~\ref{lm-4} is not true. Indeed, there is a well-known
homomorphism
$$
SL(2,\Zz)\ra SL(2,\Fp_2)\simeq \mathfrak{S}_3\fleche{\eps} \{\pm 1\}
$$
(where the isomorphism in the middle is obtained by looking at the
action on the three lines in $\Fp_2^2$, and $\eps$ is the signature),
and (for instance) the generators 
$$
S=\Bigl\{
\begin{pmatrix}
1&\pm 1\\
0&1
\end{pmatrix},
\begin{pmatrix}
1& 0\\
\pm 1&1
\end{pmatrix}
\Bigr\}
$$
(which are the analogues of the generating set for $SL(n,\Zz)$
considered in the previous remark)
all map to transpositions in $\mathfrak{S}_3$. So $\eps(r)=-1$ for any
word of odd length in the alphabet $S$.
\par
Still, while this shows that Proposition~\ref{pr-tau} can not be
applied, it remains true that for an arbitrary symmetric set of generators
$S$ of $SL(2,\Zz)$ and for odd primes $\ell$, the Cayley graph of
$SL(2,\Fp_{\ell})$ with respect to $S$ is not bipartite (because any
homomorphism $SL(2,\Fp_{\ell})\ra \{\pm 1\}$ is still trivial for
$\ell\geq 3$). Hence
this Cayley graph contains some cycle of odd length, which is easily
checked to be $\leq 2d(\ell)+1$, where  $d(\ell)$ is the diameter of
the Cayley  graph.\footnote{\ Fix some vertex $x$ and find two
  vertices $y$ and   $z$ which are neighbors but satisfy $d(x,y)\equiv
  d(x,z)\mods{2}$ (those 
  exist, because otherwise the graph would be bipartite; note that
  $d(x,y)=d(x,z)$); then follow 
  a path 
  $\gamma_1$ of length   $d(x,y)\leq d(\ell)$ from $x$ to $y$, take the
  edge from $y$ to $z$, then follow a path of length $d(z,x)=d(x,y)$
  from $z$ to $x$ to obtain a loop of length $2d(x,y)+1\leq
  2d(\ell)+1$. The example of a cycle of odd length,  
  i.e., of the   Cayley graph of  
  $\Zz/m\Zz$ with $m$ odd with respect to $S=\{\pm 1\}$ shows that this
  is best possible for arbitrary graphs, and the Ramanujan graphs of
  Lubotzky-Philips-Sarnak give examples of expanding families where
  diameter and length of shortest loop (not necessarily of odd length)
  are of the same order of magnitude,
  see~\cite[Th. 3.3.1]{sarnak-cup}.}  Since  we have an expander
family (by Property $(\tau)$), there  is a bound 
$$
d(\ell)\ll \log \ell
$$
for $\ell\geq 3$, where the implied constant depends only on $S$
(since the $(\tau)$-constant, hence the expanding ratio, is fixed);
see, e.g.,~\cite[\S 6.4]{saloff-coste}. After a look at the
character table of $SL(2,\Fp_{\ell})$, it
is not difficult to check  that this leads to sieve constants such
that 
$$
\Delta\leq |W|+|W|\exp\Bigl(-\frac{cN}{\log L}\Bigr)
\Bigl\{\max_{n\in \Lb}{\psi(n)}\sum_{m\in\Lb}{\psi(m)^2}\Bigr\}
$$
where $L=\max{\Lb}$ and $c>0$ depends only on $S$ (see~(\ref{eq-psi})
for the definition of $\psi(m)$).
\par
For the problem of irreducibility
(which is not the most interesting question about quadratic
polynomials, perhaps...) taking $\Lb$ the odd primes $\leq L$, this
leads to 
\begin{equation}\label{eq-sldeux}
|\{w\in W\,\mid\, \det(F_w-T)\in \Zz[T]\text{ is  reducible}\}|
\ll |W|\exp(-c'\sqrt{N})
\end{equation}
where $c'>0$ depends only on $S$, and as observed already, this proves 
Corollary~\ref{cor-almost-surely} and Theorem~\ref{th-sln} for
$SL(2,\Zz)$.  
\end{remark}

\begin{remark}
The work of Bourgain, Gamburd and Sarnak (see~\cite{bgs} and Sarnak's
slides~\cite{sarnak-rademacher}) is based on another type of sieve
settings, which amounts to the following. First, we have a finitely
generated group $\Gamma$ which is a discrete subgroup of a matrix
group over $\Zz$, acting on an affine algebraic variety $V/\Zz$. Then
the sieve setting is $(\Gamma\cdot v,\{\text{primes}\},\rho_{\ell})$
where $\Gamma\cdot v$ is the orbit of a fixed element $v\in V(\Zz)$,
and $\rho_{\ell}$ is the reduction map to the finite orbit of the
reduction in $V(\Fp_{\ell})$ (with uniform density).  The siftable set
if a subset $Y$ of the orbit 
defined by the images of  elements of $\Gamma$ of bounded word-length
or bounded norm, with counting measure and identity map.
\end{remark}

\begin{remark}
If we consider an abstract finitely generated group $\Gamma$, and wish
to investigate by sieve methods some of its properties, the family of
reduction maps modulo primes makes no sense. We want to point
out a family $(\rho_{\ell})$ that may be of
use, inasmuch as it satisfies the linear disjointness condition
(Definition~\ref{def-chinese}). 
\par
Let
$\tilde{\Lambda}$ be the set of surjective homomorphisms
$$
\rho\,:\, G\ra \rho(G)=H
$$
where $H$ is a non-abelian finite simple group,
and let $\Lambda\subset
\tilde{\Lambda}$ be a set of representatives for the equivalence
relation $\rho_1\sim \rho_2$ if and only if there exists  an
isomorphism $\rho_1(G)\ra \rho_2(G)$ such that the 
triangle
\begin{equation}\label{eq-triangle}
\xymatrix{
&G \ar[dl]_-{{\rho_1}}\ar[dr]^-{{\rho_2}}&\\
\rho_1(G) \ar[rr]&& \rho_2(G)}
\end{equation}
commutes.

\begin{lemma}\label{lm-goursat}
The system $(\rho)_{\rho\in \Lambda}$ constructed in this manner is
linearly disjoint.
\end{lemma}

This is an easy adaptation of classical variants of the Goursat-Ribet
lemmas, and is left as an exercise (see e.g.~\cite[Lemma 3.3]{ribet}).
\par
To make an efficient sieve, it would be necessary, in practice, to
have some knowledge of $\Lambda$, such as the distribution of the
orders of the finite simple quotient groups of $G$.
This is of course in itself an interesting
problem (see, e.g., the book~\cite{lubotzky-segal}).
\end{remark}

\section{The elliptic sieve}
\label{sec-elliptic-amusing}

The next application is also  apparently new, although it
concerns a sieve which is a sort of ``twisted'' version of the
classical large sieve. Let
$E/\Qq$ be an elliptic curve given by a Weierstrass equation
$$
y^2+a_1xy+a_3y=x^3+a_2x^2+a_4x+a_6,\quad\text{ where } a_i\in\Zz.
$$
Assuming the rank of $E$ is positive, let $\Lambda_E$ be the set of
primes  $\ell$ of good reduction, and  for $\ell\in\Lambda_E$, 
let $\rho_{\ell}\,:\,E(\Qq)\ra 
E(\Fp_{\ell})$ be the 
reduction map. 
\par
The natural sets $X\subset E(\Qq)$ for sieving are the finite sets of
rational points  $x\in E(\Qq)$ with (canonical or na\"{\i}ve) height
$h(x)\leq T$ for some $T\geq 0$ (with again counting measure and
$F_x=x$ for  $x\in X$). There are interesting potential
applications of such sieves, because of the following interpretation: 
a rational point $x=(r,s)\in E(\Qq)$ (in affine coordinates, so
$x$ is non-zero in $E(\Qq)$) maps to a non-zero point $E(\Fp_{\ell})$
if and only 
if $\ell$  does not divide  the denominator of the affine coordinates
$r$ and $s$ of the 
point. This shows that integral or $S$-integral points (in the affine
model above) appear naturally as (subsets of) sifted sets.

\par
We use such ideas to prove Theorem~\ref{th-ell},  showing that
most rational points have denominators divisible by many (small)
primes. Recall that  $\omega_E(x)$ is the number of primes, without 
multiplicity, dividing the denominator of the coordinates of $x$, with
$\omega_E(0)=+\infty$. We also recall the statement:

\begin{theorem}\label{th-ell2}
Let $E/\Qq$ be an elliptic curve with rank $r\geq 1$. Then we have
\begin{equation}\label{eq-count-ell}
|\{x\in E(\Qq)\mid\, h(x)\leq T\}|\sim c_E T^{r/2}
\end{equation}
as $T\ra +\infty$, for some constant $c_E>0$, and moreover for any
fixed real number $\kappa$ with $0<\kappa<1$, we have
$$
|\{x\in E(\Qq)\,\mid\, h(x)\leq T\text{ and }
\omega_E(x)<\kappa \log\log T
\}|\ll T^{r/2}(\log\log T)^{-1},
$$
for $T\geq 3$, where the implied constant depends only on $E$ and
$\kappa$. 
\end{theorem}

\begin{proof}
Let $M\simeq \Zz^r$ be a subgroup of $E(\Qq)$ such that 
$$
E(\Qq)=M\oplus E(\Qq)_{tors},
$$
and let $(x_1,\ldots,x_r)$ be a fixed $\Zz$-basis of
$M$. Moreover, let $M'$ be the group generated by
$(x_2,\ldots, x_r)$. We will in fact perform sieving only on ``lines''
directed by $x_1$. 
\par
But first of all, since the canonical height is a positive definite
quadratic form on $E(\Qq)\otimes \Rr=M\otimes\Rr$, the asymptotic
formula~(\ref{eq-count-ell}) is clear: it amounts to nothing else but
counting integral points in  $M\otimes\Rr\simeq \Rr^r$ with norm
$\sqrt{h(x)}\leq \sqrt{T}$ (this being repeated as many times as there
are torsion cosets). 
\par
Moreover, we may (for convenience) measure the size of elements in
$E(\Qq)$ using the squared $L^{\infty}$-norm
$$
\|x\|_{\infty}^2=\max{|a_i|^2},\quad \text{ for }
x=\sum{a_ix_i}+t\text{ with } t\in E(\Qq)_{tors},
$$
i.e., we have $h(x)\asymp \|x\|_{\infty}^2$ for all $x\in M$, the
implied constants depending only on $E$.
\par
Now we claim the following:

\begin{lemma}\label{lm-ells}
For any fixed $\kappa\in ]0,1[$, any fixed $x'\in M'$, any fixed
torsion point $t\in E(\Qq)_{tors}$, we have
$$
|\{x\in (t+x')\oplus \Zz x_1\,\mid\, \|x\|_{\infty}^2\leq T\text{ and }
\omega_E(x)< \kappa\log\log T
\}|\ll \sqrt{T}(\log\log T)^{-1},
$$
for $T\geq 3$, the implied constant depending only on $E$, $\kappa$
and $x_1$, but not on $x'$ or $t$.
\end{lemma}

Taking this for granted, we conclude immediately that
$$
|\{x\in E(\Qq)\,\mid\, h(x)\leq T\text{ and }
\omega_E(x)<\kappa \log\log T
\}|\ll T^{r/2}(\log\log T)^{-1},
$$
by summing the
inequality of the lemma over all $x'\in M'$ with
$\|x'\|_{\infty}^2\leq T$ and over all $t\in E(\Qq)_{tors}$ (the
number of pairs $(t,x')$ is $\ll T^{(r-1)/2}$), 
the implied constant depending only on $E$ and the choice of basis of
$M$. 
\par
Next we come to the proof of this lemma.  Fix $x'\in M'$, $t\in
E(\Qq)_{tors}$. The left-hand side of the lemma being zero unless
$\|t+x'\|_{\infty}^2\leq T$, we assume that this is the case. We will
use the following group sieve setting:
\begin{gather*}
\sieve=(\Zz x_1, \Lambda_E,\Zz x_1\ra \rho_{\ell}(\Zz
x_1)\subset\rho_{\ell}(E(\Qq)))\\
X=\{mx_1\in G\,\mid\, \|t+x'+mx_1\|_{\infty}^2=m^2\leq T\},
\quad\quad F_x=x.
\end{gather*}
\par
For any prime $\ell\in\Lambda_E$, the finite group $G_{\ell}$ is a
quotient of  $\Zz x_1$ and is isomorphic to $\Zz/\nu(\ell)\Zz$ where
$\nu(\ell)$ is the order 
of the reduction of $x_1$ modulo $\ell$. (So this sieve is really an
ordinary-looking one for integers, except for the use of reductions
modulo $\nu(\ell)$ instead of reductions modulo primes).

\begin{lemma}\label{lm-silv}
Let $x_1$ be an infinite order point on $E(\Qq)$, and $\nu(\ell)$ the
order of $x_1$ modulo $\ell$. Then all but finitely many primes $p$
occur as the value of $\nu(\ell)$ for some $\ell$ of good reduction.
\end{lemma}

\begin{proof}
For a prime $p$, consider $p x_1\in E(\Qq)$. A prime $\ell$ of good
reduction divides the denominator of the coordinates of $px_1$ if
and only if $p\equiv 0\mods{\nu(\ell)}$, which means that $\nu(\ell)$
is either $1$ or $p$. So if $p$ is \emph{not} of the form $\nu(\ell)$,
it follows that $px_1$ is an $S$-integral point, where $S$ is the
union of the set of primes of bad reduction and the finite set of
primes where $x_1\equiv 0\mods{\ell}$.  By Siegel's finiteness theorem
(see, e.g.,~\cite[Th. IX.4.3]{silverman}), there are only finitely many
$S$-integral points in $E(\Qq)$, and therefore only finitely many $p$
for which $p$ is not of the form $\nu(\ell)$.
\end{proof}

(Note that this lemma is also a trivial  consequence of a
result of Silverman~\cite[Prop. 10]{silv-wieferich} according to
which all but finitely many 
\emph{integers} are of the form $\nu(\ell)$ for some $\ell$; in fact
Silverman's result depends on a stronger form of Siegel's theorem).

The lemma allows us to sieve  $X$ using as prime sieve support $\Lb^*$
the set of $\ell\in\Lambda_E$ is such that 
$\nu(\ell)$ is a prime number $p\leq L$ (where, in case the same prime
$p$ occurs as values of $\nu(\ell)$ for two or more primes, we keep
only one), and with $\Lb=\Lb^*$ (with the usual identification of
elements with singletons in $S(\Lambda)$).
\par
From the lemma, it follows that the inequality defining the
large sieve constant here, namely
\begin{equation}\label{eq-ell-ls}
\sum_{\ell\in \Lb}{\sums_{a\mods{\nu(\ell)}}{
\Bigl|
\sum_{|m|\leq \sqrt{T}}{\alpha(m)e\Bigl(\frac{am}{\nu(\ell)}\Bigr)
\Bigr|^2}}}\leq \Delta\sum_{|m|\leq \sqrt{T}}{|\alpha(m)|^2},
\end{equation}
for all $(\alpha(m))$, can be reformulated as
$$
\sums_{p\leq L}{\sums_{a\mods{p}}{
\Bigl|
\sum_{|m|\leq \sqrt{T}}{\alpha(m)e\Bigl(\frac{am}{p}\Bigr)
\Bigr|^2}}}\leq \Delta\sum_{|m|\leq \sqrt{T}}{|\alpha(m)|^2}.
$$
where $\sums$ in the sum over $p$ indicates that only those $p$ which
occur as $\nu(\ell)$ for some $\ell$ are taken into account.
We recognize the most standard large sieve inequality, and  by
positivity, it follows that
$$
\Delta\leq 2\sqrt{T}+L^2
$$
for $L\geq 2$.
We now apply Proposition~\ref{pr-weaker-ls}: we have
\begin{equation}
\label{eq-apply-wls}
\sum_{x\in X}{\Bigl(P(x,\Lb)-P(\Lb)\Bigr)^2}
\leq \Delta Q(\Lb)
\end{equation}
where $P(x,\Lb)$, $P(\Lb)$ and $Q(\Lb)$ are defined
in~(\ref{eq-pxl-pl}), for 
any given choice of sets $\Omega_{\ell}\subset G_{\ell}$ for
$\ell\in\Lambda_E$. 
\par
We let $\Omega_{\ell}=\{-\rho_{\ell}(t+x')\}$. By the remark before the
statement of Theorem~\ref{th-ell2}, we have $\rho_{\ell}(mx_1)\in
\Omega_{\ell}$ if and only if $\ell$ divides the denominator of the
coordinates of $t+x'+mx_1$, and therefore for $x=mx_1\in X$, we have
$$
P(mx_1,\Lb)\leq \omega_E(t+x'+mx_1).
$$
On the other hand, we have 
$$
P(\Lb)=\sum_{\ell\in \Lb}{\frac{1}{|G_{\ell}|}}=
\sum_{\ell\in\Lb}{\frac{1}{\nu(\ell)}}
=\sum_{p\leq L}{\frac{1}{p}}+O(1)=\log\log L+O(1)
$$
for any $L\geq 3$, because,  by Lemma~\ref{lm-silv}, the values
$\nu(\ell)\leq L$ range over all primes $\leq L$, with only finitely
many exceptions (independently of $L$).
\par
Hence there exists $L_0$ depending on $E$, $x_1$ and $\kappa$ only,
such that if $L\geq L_0$, we have
$$
P(\Lb)\geq \frac{1+\kappa}{2}\log\log L.
$$
\par
Putting together these two inequalities, we see that if we assume
$T\leq L^2$, say, and $L\geq L'_0$ for some other constant $L'_0$
(depending on $E$, $x_1$ and $\kappa$),
then for any 
$mx_1\in X$ such that $t+x'+mx_1$ satisfies
$\omega_E(t+x'+mx_1)<\kappa\log\log T$, we have
$$
\Bigl(P(x,\Lb)-P(\Lb)\Bigr)^2\gg (\log\log T)^2,
$$
the implied constant depending only on $E$, $x_1$ and $\kappa$. 
So it follows by positivity from~(\ref{eq-apply-wls}) and the
inequality $Q(\Lb)\leq P(\Lb)\ll \log\log T$ that
\begin{align*}
|\{x\in t+x'\oplus \Zz x_1\,\mid\, \|x\|_{\infty}^2\leq T
\text{ and }\omega_E(x)<\kappa \log\log T
\}|&\ll \Delta (\log\log T)^{-1}\\
&\ll (\sqrt{T}+L^2)(\log\log T)^{-1}
\end{align*}
for any $L\geq L'_0$. If $T^{1/2}\geq L'_0$, we take  $L=T^{1/2}$ and
prove the inequality of the lemma directly, and otherwise we need only
increase the resulting implied constant to make it valid for all
$T\geq 3$, since $L'_0$ depends only on $E$, $x_1$ and $\kappa$.
\end{proof}

It would be interesting to know whether there is some
``regular'' distribution for the function $\omega_E(x)$. Notice the
similarity between the above 
discussion and the Hardy-Ramanujan results concerning the normal
order of the number of prime divisors of an integer (see
e.g.~\cite[22.11]{hardy-wright}), but note that
since the denominators of rational points $x$ are typically of size
$\exp h(x)$, they should have around 
$$
\log\log\exp(h(x))=\log(h(x))
$$
prime divisors in order to be ``typical'' integers.  However, note also
that the prime divisors accounted for in the proof above are all $\leq
T^{1/2}\simeq \sqrt{h(x)}\simeq \sqrt{\log n}$; it \emph{is} typical
behavior for an integer $n\leq T$ to have roughly $\log\log\log T$
prime divisors of this size (much more precise results along those
lines are known, due in particular to Tur\'an, Erd\"os and Kac).
\par
Note also that, as mentioned during the discussion of
Proposition~\ref{pr-weaker-ls},
applying the (apparently stronger) form of the large sieve involving
squarefree numbers would only give a bound  for the number
of points which are $\Lb$-integral. Since (for any finite set $S$),
there are only finitely many $S$-integral points, and moreover this is
used in the proof of Lemma~\ref{lm-silv}, this would not be a
very interesting conclusion.
\par
We can relate this sieve, more precisely Lemma~\ref{lm-ells}, to
so-called \emph{elliptic divisibility sequences}, a notion introduced
by M. Ward and currently the subject of a number of investigations by
Silverman, T. Ward, Everest, and others (see 
e.g.~\cite{silv-eds},~\cite{ward},~\cite{einsiedler-everest-ward}).
This shows that the proposition above has very concrete
interpretations. 

\begin{proposition}\label{pr-eds}
Let  $(W_n)_{n\geq 0}$ be an \emph{unbounded} sequence of integers
such that 
\begin{gather*}
W_0=0,\quad W_1=1,\quad W_2W_3\not=0,\quad W_2\mid W_4\\
W_{m+n}W_{m-n}=W_{m+1}W_{m-1}W_n^2-W_{n+1}W_{n-1}W_m^2,\quad
\text{ for }m\geq n\geq 1,\\
\Delta= W_4W_2^{15}-W_3^3W_2^{12}+3W_4^2W_2^{10}-20W_4W_3^3W_2^7
\\\quad\quad\quad
+4W_4^3W_2^5+16W_3^6W_2^4+8W_4^2W_3^3W_2^2+W_4^4\not=0.
\end{gather*}
Then for any $\kappa$ such that $0<\kappa<1$, we have
$$
|\{n\leq N\,\mid\, \omega(W_n)< \kappa\log\log N\}|
\ll \frac{N}{\log\log N}
$$
for $N\geq 3$, where the implied constant depends only on $\kappa$ and
$(W_n)$. 
\end{proposition}

\begin{proof}
This depends on the relation between elliptic divisibility sequences
and pairs $(E,x_1)$ of an elliptic curve $E/\Qq$ and a point $x_1\in
E(\Qq)$. Precisely (see e.g.~\cite[\S 2]{einsiedler-everest-ward})
there exists such a pair $(E,x_1)$ with 
$x_1$ of infinite order such 
that if $(a_n)$, $(b_n)$, $(d_n)$ are the (unique) sequences of
integers with $d_n\geq 1$, $(a_n,d_n)=(b_n,d_n)=1$ and
$$
nx_1=\Bigl(\frac{a_n}{d_n^2},\frac{b_n}{d_n^3}\Bigr),
$$
then we have
$$
d_n\mid W_n\text{ for } n\geq 1
$$
(without the condition $\Delta=0$, this is still true provided
\emph{singular} elliptic curves are permitted; the  condition that
$(W_n)$  be unbounded implies that $x_1$ is of infinite order).
\par
Now the $d_n$ are precisely the denominators of the coordinates of the
points in $\Zz x_1$, and we have therefore
$$
\omega(W_n)\geq \omega(d_n)=\omega_E(nx_1).
$$
Hence Lemma~\ref{lm-ells} gives the desired result.
\end{proof}

The ``simplest'' example is the sequence $(W_n)$ given by
\begin{gather*}
W_0=0,\quad W_1=1,\quad W_2=1,\quad W_3=-1,\quad W_4=1,\\
W_n=\frac{W_{n-1}W_{n-3}+W_{n-2}^2}{W_{n-4}},\quad\text{ for }
n\geq 4
\end{gather*}
(sequence \texttt{A006769} in the Online Encyclopedia of Integer
Sequences), which corresponds to case of $E\,:\, y^2-y=x^3-x$ and
$x_0=(0,0)$.
\par
Finally, it will be noticed that the same reasoning and similar
results hold for elements of non-degenerate divisibility sequences
$(u_n)$ defined 
by linear recurrence relations of order $2$, e.g., $u_n=a^n-1$ where
$a\geq 2$ is an integer. (The analogue of Silverman's theorem here is
a result of Schinzel, and the rest is easy).


\section{Sieving  for Frobenius over finite fields}
\label{sec-ls-ff}

The final example of large sieve we discuss concerns the distribution of
geometric Frobenius conjugacy classes in finite monodromy groups,
refining the arguments and methods in~\cite{sieve}. It is a good
example of a coset sieve as in Section~\ref{sec-spec-coset}.
\par
The precise setting is as follows (see also~\cite{sieve}). Let $q$ be
a power of a prime $p$, 
let $U/\Fp_q$ be a smooth affine geometrically connected algebraic
variety of dimension $d\geq 1$ over $\Fp_q$. Put $\bar{U}=U\times
\bar{\Fp}_q$,  the extension of scalars to an algebraic closure of
$\Fp_q$.  
\par
Let $\bar{\eta}$ denote a geometric generic 
point of $U$. We consider the coset sieve with
\begin{equation}
\label{seq-setting-ff}
G=\pi_1(U,\bar{\eta}),\quad \ggeom=\pi_1(\bar{U},\bar{\eta}),\quad
G/\ggeom\simeq \Gal(\bar{\Fp}_q/\Fp_q)\simeq \hat{\Zz},
\end{equation}
so that we have the exact sequence
$$
1\ra \ggeom\ra G\fleche{d} \hat{\Zz}\ra 1,
$$
the last arrow being the ``degree''.
\par
We assume given  a family of representations
$$
\rho_{\ell}\,:\, \pi_1(U,\bar{\eta})\ra GL(r,k_{\ell})
$$
for $\ell$ in a subset $\Lambda$ of the set of prime numbers, where
$k_{\ell}$ is a finite field of characteristic $\ell$ and $r$ is
independent of $\ell$. By the
equivalence 
of categories between lisse sheaves of $k$-modules and continuous
actions of  $\pi_1(U,\bar{\eta})$ on finite dimensional $k$ vector
spaces, this corresponds equivalently to a system
$(\sheaf{F}_{\ell})$ of étale $k_{\ell}$-vector spaces. We then
put $G_{\ell}=\Imag(\rho_{\ell})$, the \emph{arithmetic monodromy
  group} of $\rho_{\ell}$, so that we have surjective maps
$G=\pi_1(U,\bar{\eta})\ra G_{\ell}$.
\par
The  siftable set we are interested in is given by  $X=U(\Fp_q)$, with
counting measure, with the map $x\mapsto F_x\in G_{\ell}^{\sharp}$
given by the geometric Frobenius conjugacy 
class at the rational point $x\in U(\Fp_q)$ (relative to the field
$\Fp_q$). Since, in the exact sequence above, we have
$d(F_x)=-1\in\hat{\Zz}$ for 
all $x\in U(\Fp_q)$, this 
corresponds to the sieve setting $(Y,\Lambda,(\rho_{\ell}))$ where $Y$,
as in Section~\ref{sec-spec-coset}, is the set of conjugacy classes in
$\pi_1(U,\bar{\eta})$ with degree $-1$.
\par 
Then, concerning the exponential sums of
Proposition~\ref{pr-coset-sieve}, we have two basic bounds.

\begin{proposition}\label{pr-ff}
Assume that the representations $(\rho_m)$ for $m\in S(\Lambda)$ are
such that, for all squarefree numbers $m$ divisible only by primes in
$\Lambda$, the map 
$$
\pi_1(\bar{U},\bar{\eta})\ra \ggeom_m=\prod_{\ell\mid m}{\ggeom_{\ell}}
$$
is onto. With notation as before and as in
Proposition~\ref{pr-coset-sieve}, we have:
\par
\emph{(1)} If $G_{\ell}$ is a group of order prime
to $p$ for all $\ell\in\Lambda$, then
$$
W(\pi,\tau)=\delta((m,\pi), (n,\tau))q^d+
O\Bigl(q^{d-1/2}|G_{[m,n]}|(\dim \pi)(\dim\tau)\Bigr)
$$
for $m$, $n\in S(\Lambda)$,
$\pi\in\Pi_m^*$, $\tau\in \Pi_n^*$, 
where the implied constant depends only on $\bar{U}$.
\par
\emph{(2)} If $d=1$ \emph{(}$U$ is a curve\emph{)} and if the sheaves
$\sheaf{F}_{\ell}$ are of the form
$\sheaf{F}_{\ell}=\tilde{\sheaf{F}}_{\ell}/\ell\tilde{\sheaf{F}}_{\ell}$
for some 
\emph{compatible family} of torsion-free $\Zz_{\ell}$-adic sheaves
$\tilde{\sheaf{F}}_{\ell}$, then
$$
W(\pi,\tau)=\delta((m,\pi), (n,\tau))q+
O\Bigl(q^{1/2}(\dim \pi)(\dim\tau)\Bigr)
$$
where the implied constant depends only on
the compactly-supported  
Euler-Poincaré characteristics of $\bar{U}$ and of the  compatible
system $(\tilde{\sheaf{F}}_{\ell})$ on $\bar{U}$.
\end{proposition}

Recall that a system $(\tilde{\sheaf{F}}_{\ell})$ of étale sheaves of
torsion-free $\Zz_{\ell}$-modules is \emph{compatible} if, for every $\ell$,
every extension field $\Fp_{q^r}$ of $\Fp_q$, any $v\in U(\Fp_{q^r})$,
the characteristic polynomial $\det(1-TF_v\mid
\tilde{\sheaf{F}}_{\ell})$ has integer coefficients and is independent
of $\ell$. Then the Euler-Poincaré characteristic
$\chi_c(\bar{U},\tilde{\sheaf{F}}_{\ell})$ is independent of $\ell$,
being the degree of the $L$-function
$$
\prod_{x\in |U|}{\det(1-TF_x\mid \tilde{\sheaf{F}}_{\ell})^{-1}}
$$
of the sheaf as a rational function ($|U|$ is the set of all closed
points of $U$). 

\begin{proof}
This is essentially Proposition~5.1 of~\cite{sieve}, in the case
$\ell=m$, $\ell'=n$ at least. We repeat the proof since it is quite
short.
\par
By~(\ref{eq-wpitau}) and the definition of $X$, we have 
$$
W(\pi,\tau)=\frac{1}{\sqrt{|\hat{\Gamma}_m^{\pi}|
|\hat{\Gamma}_n^{\tau}|}}
\sum_{u\in U(\Fp_q)}{\Tr([\pi,\bar{\tau}]\rho_{[m,n]}(F_u))}
$$
where the sum is the sum of local traces of Frobenius for a continuous
representation of $\pi_1(U,\bar{\eta})$. We can view
$[\pi,\bar{\tau}]$ as a representation acting on a
$\bar{\Qq}_{\ell}$-vector space for some prime $\ell\not=p$, and then
this expression may be interpreted as 
the sum of local traces of Frobenius at points in $U(\Fp_q)$ for some
lisse $\bar{\Qq}_\ell$-adic sheaf $\sheaf{W}(\pi,\tau)$ on $U$.
\par
By the Grothendieck-Lefschetz Trace Formula
(see, e.g.,~\cite{grothendiec},~\cite{sga4half},~\cite[VI.13]{milne}), we
have then
$$
W(\pi,\tau)=\frac{1}{\sqrt{|\hat{\Gamma}_m^{\pi}|
|\hat{\Gamma}_n^{\tau}|}}
\sum_{i=0}^{2d}{(-1)^i
\Tr (\frob\,\mid\, H^i_c(\barre{U},\sheaf{W}(\pi,\tau)))}
$$
where $\frob$ denotes the global geometric Frobenius on $\bar{U}$.
\par
Since the representation corresponding to $\sheaf{W}(\pi,\tau)$
factors through a finite group, this sheaf is pointwise pure of weight
$0$. Therefore, by Deligne's Weil II Theorem~\cite[p. 138]{weil2}, the
eigenvalues of the geometric Frobenius automorphism $\frob$ acting on 
$H^i_c(\barre{U},\sheaf{W}(\pi,\tau))$ are algebraic 
integers, all conjugates of which are of absolute value $\leq
q^{i/2}$. 
\par
This yields
$$
W(\pi,\tau)=\frac{1}{\sqrt{|\hat{\Gamma}_m^{\pi}|
|\hat{\Gamma}_n^{\tau}|}}
\Tr (\frob\,\mid\, H^{2d}_c(\barre{U},\sheaf{W}(\pi,\tau)))
+O\Bigl(\sigma'_c(\bar{U},\sheaf{W}(\pi,\tau))q^{d-1/2}\Bigr),
$$
with an absolute implied constant, where
$$
\sigma'_c(\bar{U},\sheaf{W}(\pi,\tau)))=
\sum_{i=0}^{2d-1}{\dim H^i_c(\barre{U},\sheaf{W}(\pi,\tau))}.
$$
For the ``main term'', we use the formula
$$
H^{2d}_c(\bar{U},\sheaf{W}(\pi,\tau))=
V_{\pi_1(\bar{U},\bar{\eta})}(-d)
$$
where $V=\sheaf{W}_{\barre{\eta}}(\pi,\tau)$ is the space on which the
representation which ``is''  the sheaf acts. But, by assumption, when
we factor the 
representation (restricted to the geometric fundamental group) as follows
$$
\pi_1(\bar{U},\bar{\eta})\fleche{\rho_{[m,n]}}
\ggeom_{[m,n]}\fleche{[\pi,\bar{\tau}]} 
GL((\dim\pi)(\dim\tau),\bar{\Qq}_{\ell}),
$$
the first map is \emph{surjective}. Hence we have
$$
V_{\pi_1(\bar{U},\bar{\eta})}(-d)
=W_{\ggeom_{[m,n]}}(-d)
$$
with $W$ denoting the space of $[\pi,\bar{\tau}]$. As we are dealing
with linear representations of finite groups in characteristic $0$,
this coinvariant space is the same as the space of invariants, and its
dimension is the
multiplicity of the trivial representation in $[\pi,\bar{\tau}]$
(acting on $W$, i.e., restricted to $\ggeom_{[m,n]}$). By
Lemma~\ref{lm-ortho}, we have therefore
$$
H^{2d}_c(\bar{U},\sheaf{W}(\pi,\tau))=0
$$
if $(m,\pi)\not=(n,\tau)$. Otherwise the dimension is
$|\hat{\Gamma}^{\pi}_m|$ and the Tate twist means the global
Frobenius acts on the invariant space by multiplication by $q^d$ (the
eigenvalue is exactly $q^d$, not a root of unity times $q^d$ because
in the latter case would correspond to a situation where
$[\pi,\bar{\tau}]$ vanishes identically on $Y_{[m,n]}$, which is
excluded by the choice of $\Pi_m^*$, $\Pi_n^*$ in
Proposition~\ref{pr-coset-sieve}). This
gives 
$$
W(\pi,\tau)=\delta((m,\pi),(n,\tau))q^d
+O\Bigl(\sigma'_c(\bar{U},\sheaf{W}(\pi,\tau))q^{d-1/2}\Bigr),
$$
with an absolute implied constant.
\par
To conclude, in Case~(1), we appeal to
Proposition~4.7 of~\cite{sieve}, which gives the desired estimate
directly. 
In Case~(2), we will apply Proposition~4.1 of~\cite{sieve}, but
however we argue a bit differently\footnote{\ The result there yields
$
\sigma'_c(\bar{U},\sheaf{W}(\pi,\tau))\leq
(\dim[\pi,\bar{\tau}])(1-\chi+\omega([m,n])w) 
$
and we do not want the term $\omega([m,n])$, which would lead to a
loss of $\log\log L$ below...}. Namely, we claim that
\begin{equation}\label{eq-claim}
\sigma'_c(\bar{U},\sheaf{W}(\pi,\tau))\leq
(\dim[\pi,\bar{\tau}])(1-\chi+w|S|),
\end{equation}
where  $\chi=\chi_c(\bar{U},\Qq_{\ell})$, and
$w$ is the sum of Swan conductors of 
$\tilde{\sheaf{W}}_{\ell}$ at the ``points at infinity'' $x\in
S\subset U(\bar{\Fp}_q)$, which is independent of $\ell$, being equal
to 
$$
\chi\rank
\tilde{\sheaf{W}}_{\ell}-
\chi_c(\bar{U},\tilde{\sheaf{W}}_{\ell})
$$
where both terms are independent of $\ell$. This provides the stated
estimate (2).
\par
To check~(\ref{eq-claim}), we  look at the proof of loc. cit. with the
current notation, and extract the bound 
$$
\sigma'_c(\bar{U},\sheaf{W}(\pi,\tau))\leq 
(\dim[\pi,\bar{\tau}])\Bigl(1-\chi+\sum_{x\in S}{
\max_{\ell\mid [m,n]}{\swan_x(\tilde{\sheaf{W}}_{\ell})}}\Bigr).
$$
\par
Then, by positivity of the Swan conductors, we note that
$$
\swan_x(\tilde{\sheaf{W}}_{\ell})
\leq \sum_{x\in S}{\swan_x(\tilde{\sheaf{W}}_{\ell})}=w
$$
for each $x$ and $\ell\mid[m,n]$ (we use here the compatibility of the
system), so that
$$
\max_{\ell\mid [m,n]}{\swan_x(\tilde{\sheaf{W}}_{\ell})}\leq w,
$$
and
$$
\sum_{x\in S}{\max_{\ell\mid [m,n]}
{\swan_x(\tilde{\sheaf{W}}_{\ell})}}
\leq w|S|
$$
which concludes the proof.
\end{proof}

To apply the bounds for the exponential sums to the estimation of the
large sieve constants, we need bounds for the quantities
\begin{equation}\label{eq-first-bd}
\max_{m,\pi}\Bigl\{ q^d+C(\dim\pi)\sumb_{n\leq L}{|G_{[m,n]}|
\sum_{\tau\in\Pi_n^*}{
(\dim \tau)}}\Bigr\}
\end{equation}
in the first case and
\begin{equation}\label{eq-second-bd}
\max_{m,\pi}\Bigl\{ q^d+C(\dim\pi)\sumb_{n\leq L}{
\sum_{\tau\in\Pi_n^*}{(\dim \tau)}}\Bigr\}
\end{equation}
in the second case.
\par
For this purpose, we make the following assumptions: for all
$\ell\in\Lambda$, and $\pi\in\Pi_{\ell}^*$, we have
\begin{equation}\label{eq-strong-bounds}
|G_{\ell}|\leq (\ell+1)^s,\quad \dim\pi\leq (\ell+1)^v,\quad
\sum_{\pi\in \Pi_{\ell}^*}{(\dim \pi)}\leq
(\ell+1)^t, 
\end{equation}
where $s$, $t$ and $v$ are non-negative integers. In the notation of
Section~\ref{sec-reprs}, the second and
third are implied by
$$
A_{\infty}(G_{\ell})\leq (\ell+1)^v,\quad A_1(G_{\ell})\leq (\ell+1)^t
$$
respectively.
\par 
Here are some examples; the first two are results proved in
Section~\ref{sec-reprs} (see Example~\ref{ex-a}).
\par
--  if $G_{\ell}$ is a subgroup of $GL(r,\Fp_{\ell})$,  we can take
$s=r^2$, $v=r(r-1)/2$, $t=r(r+1)/2$.
\par
--  if $G_{\ell}$ is a subgroup of symplectic similitudes for some
non-degenerate alternating  form of rank $2g$, we can take
$s=g(2g+1)+1$, $v=(s-(g+1))/2=g^2$, $t=g^2+g+1$.
\par
-- in particular, if $G_{\ell}\subset GL(2,\Fp_{\ell})$ and
$\ggeom=SL(2,\Fp_{\ell})$, we have
\begin{equation}
|G_{\ell}|\leq \ell^4,\quad \max(\dim\pi)=\ell+1,\quad
\sum_{\pi\in \Pi_{\ell}^*}{(\dim \pi)}\leq (\ell+1)^3.
\end{equation}
This particular case can be checked easily by looking at the character
table for $GL(2,\Fp_{\ell})$ and $SL(2,\Fp_{\ell})$. See also the
character tables of $GL(3,\Fp_{\ell})$ and $GL(4,\Fp_{\ell})$
in~\cite{steinberg} for those cases.

\begin{remark} 
In~\cite{sieve}, we used different assumptions on the size of
the monodromy groups and the degrees of their representations. The
crucial feature of~(\ref{eq-strong-bounds}) is that $A_1(G_{\ell})$
and $A_{\infty}(G_{\ell})$  are bounded by \emph{monic}
polynomials. Having polynomials with constant terms $>1$ would mean,
after multiplicativity is applied, that $A_{\infty}(G_m)$ and $A_1(G_m)$
would be bounded by polynomials times a divisor function; on average
over $m$, this would mean a loss of a power of logarithm, which in
large sieve situation (as above with irreducibility of zeta functions
of curves) is likely to overwhelm the saving coming from using
squarefree numbers in the sieve. In ``small sieve'' settings, the
loss from divisor functions is reduced to a power of 
$\log\log |X|$, which may remain reasonable, and may be sufficient
justification for using simpler but  weaker polynomial bounds.
\end{remark}

Let
\begin{equation}\label{eq-psi}
\psi(m)=\prod_{\ell\mid m}{(\ell+1)}.
\end{equation}
It follows by multiplicativity from~(\ref{eq-strong-bounds}) that we
have 
\begin{equation}\label{eq-strong-bounds2}
|G_{m}|\leq \psi(m)^s,\quad \dim\pi\leq \psi(m)^v,\quad
\sum_{\pi\in \Pi_{m}^*}{(\dim \pi)}\leq
\psi(m)^t, 
\end{equation}
for all squarefree $m$. 
\par
We wish to sieve with the prime sieve support
$\Lb^*=\{\ell\in\Lambda\,\mid\, \ell\leq L\}$ for some $L$.
The first idea for the sieve is to use the traditional sieve support
$\Lb_1$ which is the set of squarefree integers $m\leq  L$ divisible
only by primes in $\Lambda$.
However, since we have $\psi(m)\ll m\log\log m$, and this upper bound
is sharp (if $m$ has many small prime factors), the use of $\Lb_1$
leads to a loss of a power of a power of $\log \log L$ in the second
term in the estimation of~(\ref{eq-first-bd})
and~(\ref{eq-second-bd}). As described by  
Zywina~\cite{zywina}, this can be recovered using the trick of
sieving using only squarefree integers $m$ free of small prime
factors, in the sense that $\psi(m)\leq L+1$ instead of $m\leq L$
(which for primes $\ell$ remain equivalent with $\ell\leq L$). This
means we use the sieve support
$$
\Lb=\{m\in S(\Lambda)\,\mid\, m\text{ is squarefree and }
\psi(m)\leq L+1\}.
$$
\par
We quote both types of sieves:

\begin{corollary}\label{cor-ff}
With the above data and notation, let
$\Omega_{\ell}\subset G_{\ell}$, for all primes $\ell\in\Lambda$, be a
conjugacy-invariant subset of $G_{\ell}$ such that
$d(\Omega_{\ell})=-1$.
Then we have both
\begin{equation}\label{eq-first-est}
|\{u\in U(\Fp_q)\,\mid\, \rho_{\ell}(F_{u})\notin\Omega_{\ell}\text{
  for } \ell\leq L\}|
\leq (q^d+Cq^{d-1/2}(L+1)^A) H^{-1} 
\end{equation}
and
\begin{equation}\label{eq-second-est}
|\{u\in U(\Fp_q)\,\mid\, \rho_{\ell}(F_{u})\notin\Omega_{\ell}\text{
  for } \ell\leq L\}|
\leq (q^d+Cq^{d-1/2}L^A(\log\log L)^v) K^{-1},
\end{equation}
where
$$
H=\sumb_{\psi(m)\leq L+1}{\prod_{\ell\mid m}
{\frac{|\Omega_{\ell}|}{|\ggeom_{\ell}|-|\Omega_{\ell}|}}},\quad
K=\sumb_{m\leq L}{\prod_{\ell\mid m}
{\frac{|\Omega_{\ell}|}{|\ggeom_{\ell}|-|\Omega_{\ell}|}}},
$$
and
\par
\emph{(i)} If $p\nmid |G_{\ell}|$ for all $\ell\in\Lambda$, we can
take $A=v+2s+t+1$, and the constant $C$ depends only on $\bar{U}$.
\par
\emph{(ii)} If $d=1$ and the system $(\rho_{\ell})$ arises by reduction
of a compatible system of $\Zz_{\ell}$-adic sheaves on $U$, then we
can take $A=t+v+1$, 
and the constant $C$ depends only on the
Euler-Poincaré characteristic of $U$, the compactly-supported  
Euler-Poincaré characteristic of the compatible
system $(\tilde{\sheaf{W}}_{\ell})$ on $\bar{U}$, and on $s$, $t$, $v$
in the case of~\emph{(\ref{eq-second-est})}.
\end{corollary}

\begin{proof}
From Proposition~\ref{pr-w-to-sieve}, we must estimate
$$
\Delta=\max_{m,\pi}
{\sumb_{n}{\sum_{\tau\in\Pi^*_n}{|W(\pi,\tau)|}}},
$$
where $m$ and $n$ run over $\Lb$ and $\Lb_1$, respectively.
By Proposition~\ref{pr-ff}, this is bounded by the
quantities~(\ref{eq-first-bd}) 
and~(\ref{eq-second-bd}). Using~(\ref{eq-strong-bounds2}), the result
is now straightforward, using (in the case of~(\ref{eq-second-est}))
the simple estimate
$$
\sumb_{n\leq L}{\psi(n)^A}\ll L^{A+1}
$$
for $L\geq 1$, $A\geq 0$, the implied constant depending on $A$.
\end{proof}

This theorem can be used to get a slight improvement of the ``generic
irreducibility'' results for numerators of zeta functions of curves
of~\cite{sieve} (see Section 6 of that paper for some context and in
particular~Theorem 6.2): a small power of 
$\log q$ is gained in the upper 
bound, as in Gallagher's result~\cite[Th. C]{gallagher}. We only state
one special case, for a fixed genus (see the remark following the
statement for an explanation of this restriction).

\begin{theorem}\label{th-irred-log}
Let $\Fp_q$ be a finite field of characteristic $p\not=2$, let $f\in
\Fp_q[X]$ be a squarefree monic polynomial of degree $2g$, $g\geq
1$. For $t\in \Fp_q$ which is not a zero of $f$, let $P_t\in \Zz[T]$
be the numerator of the zeta function of the smooth projective model
of the hyperelliptic curve
\begin{equation}
  \label{eq-family}
C_t\,:\, y^2=f(x)(x-t),  
\end{equation}
and let $K_t$ be the splitting field of $P_t$ over $\Qq$, which has
degree $[K_t:\Qq]\leq 2^gg!$.
Then we have
$$
|\{t\in \Fp_q\,\mid\, f(t)\not=0\text{ and }
[K_t:\Qq]<2^gg!\}|
\ll q^{1-\gamma}(\log q)^{1-\delta}
$$
where $\gamma=(4g^2+2g+4)^{-1}$ and $\delta>0$, with $\delta\sim
1/(4g)$ as $g\ra+\infty$. The implied constant depends only on $g$.
\end{theorem}

\begin{proof}
Let $S(f)\subset \Fp_q$ be the set in question. 
Proceeding as in Section~8 of~\cite{sieve}, we set up a sieve using
the sheaves $\sheaf{F}_{\ell}=R^1\pi_!\Fp_{\ell}$ for $\ell>2$, $\pi$
denoting the projection from the family of curves~(\ref{eq-family}) to
the parameter space $U=\{t\,\mid\, f(t)\not=0\}\subset
\Aa^1$. Those sheaves are tame, obtained by reducing modulo $\ell$ a
compatible system, and the geometric monodromy of $\sheaf{F}_{\ell}$
is $Sp(2g,\Fp_{\ell})$ by a result of J.K. Yu (which also follows from
a recent more general result of C. Hall~\cite{hall}; see~\cite{k-hall} for a
write-up of this special case of Hall's
result). Using~(\ref{eq-first-bd}), and the proof of
Proposition~\ref{pr-ff} to bound explicitly the implied constant, we
obtain   
$$
|S(f)|\leq (q+4g\sqrt{q}L^A)H^{-1}
$$
where $A=2g^2+g+2$ (see Example~\ref{ex-a}) and
$$
H=\min_{1\leq i\leq 4}\Bigl\{\sumb_{\psi(m)\leq L}{\prod_{\ell\mid m}
{\frac{|\Omega_{i,\ell}|}{|\ggeom_{\ell}|-|\Omega_{i,\ell}|}}}\Bigr\},
$$
the sets $\Omega_{i,\ell}$ being defined as in~\cite[\S
7,8]{sieve}. For each of these we have
\begin{equation}\label{eq-bound-omega}
\frac{|\Omega_{i,\ell}|}{|\ggeom_{\ell}|-|\Omega_{i,\ell}|}=  
\frac{\delta_i}{1-\delta_i}+
O\Bigl(\frac{1}{\ell}\Bigr)
\end{equation}
for $\ell\geq 3$, for some $\delta_i\in ]0,1[$ which is
a ``density''  of conjugacy classes satisfying certain
conditions, either in the group of 
permutations on $g$ letters or the group of signed permutations of
$2g$ letters (this follows easily
from~\cite[\S 2]{gallagher} and Sections 7,~8 of~\cite{sieve}). The
implied constant depends only on $g$. Precisely, we have
$$
\delta_1\sim \frac{1}{2g},
\quad
\delta_2\sim \frac{1}{4g},
\quad
\delta_3\sim \frac{\log 2}{\log g}, 
\quad
\delta_4\sim \frac{1}{\sqrt{2\pi g}}
$$
as $g\ra +\infty$ (see~\cite[\S 8]{sieve}).
\par
Thus, we need lower bounds for sums of the type
$$
\sumb_{\psi(m)\leq L}{\beta(m)}=\sum_{\psi(m)\leq L}{\beta(m)\mu^2(m)}
$$
where $\beta$ is a multiplicative function, roughly constant at the
primes. This is a well-studied area of analytic number theory. We can
appeal for instance to Theorem 1 of~\cite{lau-wu}; in the notation of
loc. cit., we have $f(m)=\mu^2(m)\beta(m)$, $g(m)=\psi(m)$,
with $\kappa=\delta_i/(1-\delta_i)$, $\eta=1$, 
$\alpha=1$, $\theta=1$, $\alpha'=1$, $\theta'=0$, $t(p)=0$, $C_3=0$. We
obtain 
\begin{equation}
  \label{eq-w-l}
\sumb_{\psi(m)\leq L}{\beta(m)}\gg L(\log L)^{-1+\delta_i/(1-\delta_i)},  
\end{equation}
for $L\geq 3$, where the implied constant depends only on $g$.
\par
Taking $L=q^{1/2A}$, the upper bound for $S(f)$ then follows.
\end{proof}

\begin{remark}
In~\cite{sieve}, we obtained a result uniform in terms of $g$. Here it
is certainly possible to do the same, by checking the dependency of
the estimate~(\ref{eq-w-l}) on $g$. However, notice that the gain
compared to~\cite{sieve}\footnote{\ It seems that the exponent of $q$
  is better, but this reflects the use of the ``right'' bounds for
  degrees of representations of finite symplectic groups, and this
  exponent can be obtained with the method of~\cite{sieve} also.}  is
of size $(\log q)^{\delta}$ with 
$\delta\sim 1/(4g)$, and this becomes trivial as soon as $g$ is of
size $\log\log q$. This is a much smaller range than the (already
restricted) range where the estimate of~\cite{sieve} is non-trivial,
namely $g$ somewhat smaller than $\sqrt{\log q}$. 
\end{remark}

Now we prove Theorem~\ref{th-sq-jac} stated
in the introduction.

\begin{proof}[Proof of Theorem~\ref{th-sq-jac}]
We can certainly afford to be rather brief here. The sieve setting and
siftable set are the same as in Theorem~\ref{th-irred-log}. The number
of points of 
$C_t$ and $J_t$ are given by
$$
|C_t(\Fp_q)|=q+1-\Tr(\frob\,\mid\, H^1(\bar{C}_t,\Zz_{\ell})),
\quad\quad
|J_t(\Fp_q)|=|\det(1-\frob\,\mid\, H^1(\bar{C}_t,\Zz_{\ell}))|,
$$
(for any prime $\ell\nmid p$).
Thus, defining sieving sets
\begin{align*}
\Omega^J_{\ell}&=\{g\in CSp(2g,\Fp_{\ell})\,\mid\, g\text{ is
  $q$-symplectic and } \det(g-1)\text{ is a
  square in } \Fp_{\ell}\},\\
\Omega^C_{\ell}&=\{g\in CSp(2g,\Fp_{\ell})\,\mid\, g\text{ is
  $q$-symplectic and } q+1-\Tr(g)\text{ is a
  square in } \Fp_{\ell}\},
\end{align*}
(where $q$-symplectic similitudes are those with multiplicator $q$),
we have for any prime sieve support $\Lb^*$ the inclusion
$$
\{t\in  \Fp_q\,\mid\, f(t)\not=0\text{ and 
$|\mathsf{S}_t(\Fp_q)|$ is a square} \}\subset
S(U,\Omega^{\mathsf{S}};\Lb^*), 
$$
for $\mathsf{S}\in \{C,J\}$.  By (3) and (4), respectively, of
Proposition~\ref{pr-local} in Appendix~B, we have
$$
\frac{|\Omega^{\mathsf{S}}_{\ell}|}{|Sp(2g,\Fp_{\ell})|}
\geq \frac{1}{2}\Bigl(\frac{\ell}{\ell+1}\Bigr)^g.
$$
for $\ell\geq 3$. Thus if $\Lb$ is the set of odd prime $\leq L$, we
obtain 
$$
|\{t\in  \Fp_q\,\mid\, f(t)\not=0\text{ and 
$|\mathsf{S}_t(\Fp_q)|$ is a square} \}|
\leq (q+4g\sqrt{q}L^A)H^{-1}
$$
where $A=2g^2+g+2$, and 
$$
H=\sum_{3\leq \ell\leq L}{
\frac{|\Omega^{\mathsf{S}}_{\ell}|}{|Sp(2g,\Fp_{\ell})|}}
\geq \frac{1}{2}\sum_{3\leq \ell\leq L}{\Bigl(\frac{\ell}{\ell+1}\Bigr)^g}.
$$
\par
By the mean-value theorem we have
$$
\Bigl(\frac{\ell}{\ell+1}\Bigr)^g=1-\frac{g}{\ell+1}+O(g^2(\ell+1)^{-2})
$$
for $\ell\geq 3$, $g\geq 1$, with an absolute implied constant, and
thus by the Prime Number Theorem we have
$$
H\geq \frac{1}{2}\pi(L)+O(g\log\log L+g^2)
$$
with an absolute implied constant. For $L\gg g^2\log 2g$, this gives
$$
H\gg \frac{1}{2}\frac{L}{\log L}
$$
with an absolute implied constant, and
$$
|\{t\in  \Fp_q\,\mid\, f(t)\not=0\text{ and 
$|\mathsf{S}_t(\Fp_q)|$ is a square} \}|
\ll g^2(q+q^{1/2}L^A)L^{-1}(\log L),
$$
this time with no condition on $g$ and $L$ as this is trivial when
$L\leq Cg^2\log 2g$ (which explains the poorer dependency on $g$ than
follows from what we said). So choosing $L=q^{1/(2A)}$, we obtain the
uniform estimate 
$$
|\{t\in  \Fp_q\,\mid\, f(t)\not=0\text{ and 
$|\mathsf{S}_t(\Fp_q)|$ is a square} \}|
\ll q^{1-\gamma}(\log q)
$$
with $\gamma=1/(4g^2+2g+4)$ and where the implied constant is absolute.
\end{proof}

See also the end of Appendix~A for a lower bound sieve result on the
same families of curves.

\setcounter{section}{0}
\setcounter{theorem}{0}
\renewcommand{\thetheorem}{\thesection.\arabic{theorem}}
\setcounter{equation}{0}
\renewcommand{\theequation}{\thesection.\arabic{equation}}
\setcounter{figure}{0}
\setcounter{table}{0}

\appendix 
\section*{Appendix A: small sieves}
\renewcommand{\thesection}{A} 

If we are in a general sieving situation as described in
Section~\ref{sec-prelim}, we may in many cases be interested in a
lower bound, in 
addition to the upper bounds that the large sieve naturally
provides. For this purpose we can hope to appeal to the usual
principles of small sieves, at least when $\Lambda$ is the set of prime
numbers and for some specific sieve supports. We describe this for
completeness, with no   
claim to originality, and refer to books such as~\cite{HR}, the
forthcoming~\cite{if} or~\cite[\S 7]{ant} for more detailed coverage
of the principles of sieve theory.
\par
We assume that our sieve setting is of the type
$$
\sieve=(Y,\{\text{primes}\},(\rho_{\ell})),
$$
and our sieve support will be the set $\Lb$ of squarefree numbers
$d< L$ for some parameter $L$. We write $S(X;\Omega,L)$ for
the sifted set $S(X;\Omega,\Lb)$. The siftable set is $(X,\mu,F)$ as
before. 
\par
Let 
$$
\Omega_m=\prod_{\ell\mid m}{\Omega_{\ell}}
$$
for $m$ squarefree, and for an arbitrary integrable function $x\mapsto
\alpha(x)$, denote 
$$
S_d(X;\alpha)=\int_{\{\rho_{d}(F_x)\in\Omega_{d}\}}
{\alpha(x)d\mu(x)}.
$$
For $x\in X$, let $n(x)\geq 1$ be the integer defined
by
$$
n(x)=\prod_{\stacksum{\ell\leq L}{\rho_d(F_x)\in\Omega_{\ell}}}{\ell}
$$
so that for squarefree $d\in \Lb$, we have
$\rho_d(F_x)\in \Omega_d$  if and only if $d\mid n(x)$. 
\par
Then we have
\begin{align*}
\int_{S(X;\Omega,L)}{\alpha(x)d\mu(x)}&=
\int_{\{(n(x),P(L))=1\}}{\alpha(x)d\mu(x)}\\
&=
\sum_{(n,P(L))=1}{\Bigl(\int_{\{n(x)=n\}}\alpha(x)d\mu(x)\Bigr)}
=\sum_{(n,P(L))=1}{a_n}
\end{align*}
where $P(L)$ is the product of primes $\ell< L$ and
$$
a_n=\int_{\{n(x)=n\}}{\alpha(x)d\mu(x)}.
$$
\par
Note that
$$
\sum_{n\equiv 0\mods{d}}{a_n}=S_d(X;\alpha).
$$
\par
Let now $(\lambda_d^{\pm})$ be two sequences of real numbers supported
on $\Lb$ such that $\lambda_1^{\pm}=1$ and
$$
\sum_{d\mid n}{\lambda_d^-}\leq 0\leq \sum_{d\mid n}{\lambda_d^+}
$$
for $n\geq 2$. Then, if $\alpha(x)\geq 0$ for all $x$, we have
$$
\sum_{(n,P(L))=1}{a_n}\leq \sum_{n}{\Bigl(\sum_{d\mid
(n,P(L))}{\lambda_d^+}\Bigr)a_n}
= \sum_{d<L}{\lambda_d^+ \Bigl(\sum_{n\equiv 0\mods{d}}{a_n}\Bigr)}
= \sum_{d<L}{\lambda_d^+ S_d(X;\alpha)}
$$
and similarly
$$
\sum_{(n,P(L))=1}{a_n}\geq \sum_{d<L}{\lambda_d^- S_d(X;\alpha)}.
$$
\par
It is natural to introduce the approximations (compare~(\ref{eq-rd}))
\begin{equation}\label{eq-appro}
S_d(X;\alpha)=\nu_d(\Omega_d)H+
r_d(X;\alpha),
\end{equation}
(where $\nu_d$ is the a density as in Section~\ref{sec-prelim}), which
is really a definition of $r_d(X;\alpha)$,
where the ``expected main term'' is
$$
H=\int_{X}{\alpha(x)d\mu(x)}.
$$
Then, in effect, we have proved: 

\begin{proposition}
Assume $\alpha(x)\geq 0$ for all $x\in X$. 
Let $\lambda_d^{\pm}$ be arbitrary upper and lower-bound sieve
coefficients which vanish for $d\geq L$. We have then
$$
V^-(\Omega)H-R^-(X;L)\leq 
\int_{S(X;\Omega,L)}{\alpha(x)d\mu(x)}\leq V^+(\Omega)H+R^+(X;L)
$$
where
$$
V^{\pm}(\Omega)=\sum_{d<L}{\lambda_d^{\pm}\nu_d(\Omega_d)}
\quad\quad 
\text{ and }\quad\quad
R^{\pm}(X;L)=\sum_{d<L}{|\lambda_d^{\pm}r_d(X;\alpha)|}.
$$
\end{proposition}

In fact this is not quite what is needed for applications, because
$V^{\pm}(X)$ are not yet in a form that makes them easy to
evaluate. This next crucial step (called a ``fundamental lemma'')
depends on the choice of $\lambda^{\pm}_d$ (which is by no means
obvious) and on properties of
$\Omega_d$.  For instance, we have the following (see
e.g.~\cite[Cor. 6.2]{ant}; note this by no means the most general or
best result known). 

\begin{proposition}
Let $\kappa>0$ and $y>1$. There exist upper and lower-bound sieve
coefficients 
$(\lambda_d^{\pm})$, depending only on $\kappa$ and $y$, supported on
squarefree integers $<y$, bounded by one in absolute value, with the
following properties: 
 for all $s\geq 9\kappa+1$ and
$L^{9\kappa+1}<y$, we have 
\begin{align*}
\int_{S(X;\Omega,L)}{\alpha(x)d\mu(x)}<
\Bigl(1+ e^{9\kappa+1-s}K^{10}
\Bigr)
\prod_{\ell<L}{(1-\nu_{\ell}(\Omega_{\ell}))}H
+R^+(X;L^{s}),
\\
\int_{S(X;\Omega,L)}{\alpha(x)d\mu(x)}>
\Bigl(1- e^{9\kappa+1-s}K^{10}
\Bigr)
\prod_{\ell<L}{(1-\nu_{\ell}(\Omega_{\ell}))}H
+R^-(X;L^{s}),
\end{align*}
provided  the sieving sets
$(\Omega_{\ell})$ satisfy the condition 
\begin{equation}\label{eq-dimension}
\prod_{w\leq \ell<L}{(1-\nu_{\ell}(\Omega_{\ell}))^{-1}}
\leq K\Bigl(\frac{\log L}{\log w}\Bigr)^{\kappa},
\quad\quad\text{for all $w$ and $L$, $2\leq w<L<y$},
\end{equation}
for some $K\geq 0$.
\end{proposition}

In standard applications, $r_d(X;\alpha)$ should be ``small'',
as the remainder term in  some
equidistribution theorem. Note again that this can only be true if the
family $(\rho_{d})$ is linearly disjoint. 
If this remainder is 
well-controlled on average over $d<D$, for some $D$ (as large as
possible) we can apply the above for $L$ such that $L^s<D$ (with
$s\geq 9\kappa+1$). Note that when $s$ is large enough (i.e., $L$
small enough), the coefficient $1\pm e^{9\kappa+1-s}K^{10}$ will be
close to $1$, in particular it will be \emph{positive} in the lower
bound. 
\par
Note that the condition~(\ref{eq-dimension}) holds if 
$\nu_{\ell}(\Omega_{\ell})$ is of size $\kappa
\ell^{-1}$  on average.  This is the traditional
context of a ``small sieve'' of dimension $\kappa$; we see that in the
abstract framework, this means rather that the sieving sets
$\Omega_{\ell}$ are ``of codimension $1$'' in a certain sense. The
important case $\kappa=1$ (the classical ``linear sieve'')
corresponds intuitively to sieving sets defined 
by a single irreducible algebraic condition.
\par
Note also that the factor
$$
\prod_{\ell<L}{(1-\nu_{\ell}(\Omega_{\ell}))}
$$
is the natural one to expect intuitively if
$\nu_{\ell}(\Omega_{\ell})$ is interpreted as the probability of
$\rho_{\ell}(F_x)$ being in $\Omega_{\ell}$, and the various $\ell$
being independent.  To see the connection with the quantity $H^{-1}$ in the
large sieve bound~(\ref{eq-h}), note that if $\Lb$ is the full power
set of the prime sieve support $\Lb^*$, then multiplicativity gives
$$
H=\sum_{m\in\Lb}{\prod_{\ell\mid m}{\frac{\nu(\Omega_{\ell})}
{\nu(Y_{\ell}-\Omega_{\ell})}}}=
\prod_{\ell\in \Lb^*}{\Bigl(1+\frac{\nu(\Omega_{\ell})}
{\nu(Y_{\ell}-\Omega_{\ell})}\Bigr)}=
\prod_{\ell\in\Lb^*}{\frac{1}{1-\nu_{\ell}(\Omega_{\ell})}}
$$
(recall $\nu_{\ell}(Y_{\ell})=1$). So $H^{-1}$ has exactly the same
shape as the factor above. Of course, as in small sieves,  if $\Lb$ is
as large as the power set 
of $\Lb^*$, the large sieve constant will be much too big for the
large sieve inequality to be useful, and so ``truncation'' is needed.
\par
We conclude with a simple application, related to
Theorem~\ref{th-sq-jac} and Section~\ref{sec-ls-ff}. 

\begin{proposition}\label{th-lower-sieve}
Let $q$ be a power of a prime number $p\geq 5$, $g\geq 1$ an integer
and let $f\in \Fp_q[T]$
be a squarefree polynomial of degree $2g$. For $t$ not a zero of
$f$, let $C_t$ denote the smooth projective model of the
hyperelliptic curve
$y^2=f(x)(x-t)$,
and let $J_t$ denote its Jacobian variety.
There exists an absolute constant $\alpha\geq 0$ such that
\begin{align*}
|\{
u\in \Fp_q\,\mid\,f(t)\not=0\text{ and }
|C_t(\Fp_q)|\text{ has no odd prime factor $<q^{\gamma}$}
\}|\gg
\frac{q}{\log q}\\
|\{
u\in \Fp_q\,\mid\,f(t)\not=0\text{ and }
|J_t(\Fp_q)|\text{ has no odd prime factor $<q^{\gamma}$}
\}|\gg
\frac{q}{\log q}
\end{align*}
for any $\gamma$ such that
$$
\gamma^{-1}>\alpha (2g^2+g+1)(\log \log 3g),
$$
where the implied
constants depends only on 
$g$ and $\gamma$. 
\par
In particular, for any fixed $g$, there are
infinitely many points $t\in \bar{\Fp}_q$ such that
$|C_t(\Fp_{q^{\deg(t)}})|$ has at most $\alpha(2g^2+g+1)(\log\log 3g)+2$
prime factors, and similarly for $|J_t(\Fp_{q^{\deg(t)}})|$.
\end{proposition}

\begin{remark}
(1) It may well be that $|J_t(\Fp_q)|$ is even for all $t$, since if
$f$ has a root $x_0$ in $\Fp_q$, it will define a non-zero point of
order $2$ in $J_t(\Fp_q)$.
\par
(2) See e.g.~\cite{cojocaru} for results on almost prime values of
group orders of elliptic curves over $\Qq$ modulo primes; except for
CM curves, they are conditional on GRH.
\end{remark}

\begin{proof}
Obviously we use the same coset sieve setting and siftable set as in
Theorems~\ref{th-irred-log} and Theorem~\ref{th-sq-jac}, and consider
the sieving sets
\begin{align*}
\Omega^{J}_{\ell}&=
\{g\in CSp(2g,\Fp_{\ell})\,\mid\, g\text{ is
  $q$-symplectic and } \det(g-1)=0\in\Fp_{\ell}\}
\},\\
\Omega^{C}_{\ell}&=
\{g\in CSp(2g,\Fp_{\ell})\,\mid\, g\text{ is
  $q$-symplectic and } \Tr(g)=q+1\}
\},
\end{align*}
for $\ell\geq 3$, where $\mathsf{S}\in\{C,J\}$. By (5) and (6) of
Proposition~\ref{pr-local}, we have
$$
\nu_{\ell}(\Omega^{\mathsf{S}}_{\ell})=
\frac{|\Omega^{\mathsf{S}}_{\ell}|}{|Sp(2g,\Fp_{\ell})|}
\leq \min\Bigl(1,\frac{\ell^{g-1}}{(\ell-1)^g}\Bigr),
$$
from which~(\ref{eq-dimension}) can be checked to hold with $\kappa=1$
and $K\ll \log g$ (consider separately primes $\ell<g$ and $\ell\geq
g$).
\par
Coming to the error term $R^-(X;L)$, individual estimates
for $r_d(X;\alpha)$ with $\alpha(x)=1$ amount to estimates for the
error term in the Chebotarev density theorem. Using
Proposition~\ref{pr-ff} in the standard way we obtain
$$
r_d(X;\alpha)\ll gq^{1/2}|\Omega_{d}^{\mathsf{S}}|^{1/2}
\ll gq^{1/2}\Bigl(\psi(d)^{2g^2+g}d^{g-1}\varphi(d)^{-g}\Bigr)^{1/2},
$$
with absolute implied constants (see
also~\cite[Th. 1.3]{quad-twists}), and hence
$$
R^-(X;L^s)\ll gq^{1/2}L^{s(2g^2+g+1)/2}(\log\log L^s)^{g^2+g},
$$
for any $s\geq 1$, with an absolute implied constant.
\par
Let $s=\log 2+10\log K\ll \log\log 3g$, and let $\eps>0$ be arbitrarily
small. Then we can take 
$$
L=q^{(s(2g^2+g+1))^{-1}-\eps}
$$
in the lower bound sieve, which gives
\begin{multline*}
|\{t\in \Fp_q\,\mid\, f(t)\not=0\text{ and }
|\mathsf{S}_t(\Fp_q)|\text{ has no odd prime factor } <L
\}|
\\
\gg
q\prod_{\stacksum{\ell<L}{\nu_{\ell}(\Omega_{\ell}^{\mathsf{S}})<1}}
{(1-\nu_{\ell}(\Omega_{\ell}^{\mathsf{S}}))}
\gg
q\prod_{3g<\ell<L}
{\Bigl(1-\frac{\ell^{g-1}}{(\ell+1)^g}\Bigr)}
\end{multline*}
provided $L>3g$, say, the implied constant being absolute.
Putting all together, the theorem follows now easily.
\end{proof}

\setcounter{section}{0}
\setcounter{theorem}{0}
\renewcommand{\thetheorem}{\thesection.\arabic{theorem}}
\setcounter{equation}{0}
\renewcommand{\theequation}{\thesection.\arabic{equation}}
\setcounter{figure}{0}
\setcounter{table}{0}

\appendix 
\renewcommand{\thesection}{B} 
\section*{Appendix B: local density computations over finite fields}

In Sections~\ref{sec-sieve-arith},~\ref{sec-ls-ff}, and in the
previous Appendix, we have quoted various estimates for the
``density'' of certains subsets of matrix groups over finite fields,
which are required to prove lower (or upper) bounds for the saving
factor $H$ in certain applications of the large sieve inequalities.
We prove those statements here, relying mostly on the work of
Chavdarov~\cite{chavdarov} to link such densities with those of
polynomials of certain types which are much easier to compute.  In one
case, however, we use the Riemann Hypothesis over finite fields to
estimate a multiplicative exponential sum.

\begin{proposition}\label{pr-local}
Let $\ell\geq 3$ be a prime number.
\par
\emph{(1)} Let $G=SL(n,\Fp_{\ell})$ or $G=Sp(2g,\Fp_{\ell})$, with
$n\geq 2$ or $g\geq 1$. Then we have
$$
\frac{1}{|G|}
|
\{g\in G\,\mid\, \det(g-T)\in \Fp_{\ell}[X]\text{ is irreducible}\}
|
\gg 1
$$
where the implied constant depends only on $n$ or $g$.
\par
\emph{(2)} Let  $G=SL(n,\Fp_{\ell})$ or $G=Sp(2g,\Fp_{\ell})$, with
$n\geq 2$ or $g\geq 1$, let $i$, $j$ be integers with $1\leq i,j\leq
n$ or $1\leq i,j\leq 2g$ respectively. Then we have
$$
\frac{1}{|G|}
|
\{g=(g_{\alpha,\beta})\in G\,\mid\, g_{i,j}\in \Fp_{\ell}\text{
 is not a square}\}
|
\gg 1
$$
where the implied constant depends only on $n$ or $g$.
\par
\emph{(3)} Let $G=CSp(2g,\Fp_{\ell})$ with $g\geq 1$, and denote by
$m(g)\in\Fp_{\ell}^{\times}$ the multiplicator of a symplectic
similitude $g\in G$. Then for any $q\in \Fp_{\ell}^{\times}$, we have 
$$
\frac{1}{|Sp(2g,\Fp_{\ell}|}|
\{
g\in G\,\mid\, m(g)=q\text{ and } \det(g-1)\text{ is a
  square in } \Fp_{\ell}
\}
|\geq \frac{1}{2}\Bigl(\frac{\ell}{\ell+1}\Bigr)^g.
$$
\par
\emph{(4)} Let $G=CSp(2g,\Fp_{\ell})$ with $g\geq 1$. Then for any
$q\in \Fp_{\ell}^{\times}$, we have  
$$
\frac{1}{|Sp(2g,\Fp_{\ell}|}|
\{
g\in G\,\mid\, m(g)=q\text{ and } q+1-\Tr(g)\text{ is a
  square in } \Fp_{\ell}
\}
|\geq \frac{1}{2}\Bigl(\frac{\ell}{\ell+1}\Bigr)^g.
$$
\par
\emph{(5)} Let $G=CSp(2g,\Fp_{\ell})$ with $g\geq 1$. Then for any
$q\in \Fp_{\ell}^{\times}$, we have  
$$
\frac{1}{|Sp(2g,\Fp_{\ell}|}|
\{
g\in G\,\mid\, m(g)=q\text{ and } \det(g-1)=0
\}
|\leq \min\Bigl(1,\frac{\ell^{g-1}}{(\ell-1)^g}\Bigr).
$$
\par
\emph{(6)} Let $G=CSp(2g,\Fp_{\ell})$ with $g\geq 1$. Then for any
$q\in \Fp_{\ell}^{\times}$, we have  
$$
\frac{1}{|Sp(2g,\Fp_{\ell}|}|
\{
g\in G\,\mid\, m(g)=q\text{ and } q+1-\Tr(g)=0
\}
|\leq \min\Bigl(1,\frac{\ell^{g-1}}{(\ell-1)^g}\Bigr).
$$
\end{proposition}

\begin{proof}
(1) (Compare with~\cite[\S 3]{chavdarov},~\cite[Lemma 7.2]{sieve}).
Take the case $G=SL(n,\Fp_{\ell})$, for instance. We need to compute
$$
\frac{1}{|G|}\sum_{f\in\tilde{\Omega}_{\ell}}{|\{g\in
G\,\mid\,\det(g-T)=(-1)^nf\}|},
$$
where $f$ runs over the set $\tilde{\Omega}_{\ell}$ of irreducible
monic polynomials $f\in\Fp_{\ell}[T]$ of degree $n$ with
$f(0)=1$. For each $f$, we have 
$$
|\{g\in G\,\mid\,\det(g-T)=f\}|\gg
\frac{|G|}{\ell^{n-1}}
$$
by the argument in~\cite[Th. 3.5]{chavdarov} (note that the algebraic
group $SL(n)$ is 
connected and simply connected), and the number of $f$ is close to
$\frac{1}{n} \ell^{n-1}$ as $\ell\ra +\infty$, by identifying the set 
of $f$ with the set of Galois-orbits of elements of norm $1$ in
$\Fp_{\ell^n}$ which are of degree $n$ and not smaller. All this
implies the result for $G$, the  case of the
symplectic group being  similar.
\par
(2) By detecting squares using the Legendre character, we need to
compute 
$$
\frac{1}{2|G|}\sum_{\stacksum{
g\in G}{g_{i,j}\not=0}}{
\Bigl(1+\Bigl(\frac{g_{i,j}}{\ell}\Bigr)\Bigr)}
$$
where $(\frac{\cdot}{\ell})$ is the non-trivial quadratic character
of $\Fp_{\ell}^{\times}$. Let $\Gg$ be the algebraic group $SL(n)$ or
$Sp(2g)$ over $\Fp_{\ell}$, $d$ its dimension (either $n^2-1$ or
$2g^2-g$). Since  
$\Gg\cap\{g_{i,j}=0\}$
is obviously a proper closed subset of the geometrically connected
affine variety $\Gg$, the affine variety
$$
\Gg_{i,j}=\Gg-\Gg\cap\{g_{i,j}=0\}
$$
over $\Fp_{\ell}$ is geometrically connected of dimension $d$, and 
we have
$$
|\G_{i,j}(\Fp_{\ell})|=
|\{g\in \Gg(\Fp_{\ell})\,\mid\, g_{i,j}\not=0\}|\gg |\Gg(\Fp_{\ell})|,
$$
for $\ell\geq 3$. This means that it is enough to prove
$$
\sum_{g\in \Gg_{i,j}(\Fp_{\ell})}
{\Bigl(\frac{g_{i,j}}{\ell}\Bigr)}\ll
\ell^{d-1/2}
$$
for $\ell\geq 3$, the implied constant depending only on $\Gg$. Such a
bound follows (for instance) from the fact that this sum 
is a multiplicative character sum over the $\Fp_{\ell}$-rational points of
the geometrically connected affine algebraic variety $\Gg_{i,j}$ of
dimension $d$. 
\par
Instead of looking for an elementary proof (which may well exist), we
invoke the powerful $\ell$-adic cohomological formalism  (see
e.g.~\cite[11.11]{ant} for an introduction, and compare with the proof
of Proposition~\ref{pr-ff}). Using the (rank $1$) Lang-Kummer sheaf
$\mathcal{K}=\Lb_{(\tfrac{g_{i,j}}{\ell})}$, we have by the
Grothendieck-Lefschetz trace formula 
$$
\sum_{g\in \Gg_{i,j}(\Fp_{\ell})}
{\Bigl(\frac{g_{i,j}}{\ell}\Bigr)}
=\sum_{g\in \Gg_{i,j}(\Fp_{\ell})}
{\Tr(\frob_{g,\ell}\mid \mathcal{K})}=
\sum_{k=0}^{2d}{\Tr(\frob\mid H^k_c(\overline{\G}_{i,j},\mathcal{K}))}
$$ 
where $\frob_{g,\ell}$ (resp. $\frob$) is the local (resp. global)
geometric Frobenius for $g$ seen as defined over $\Fp_{\ell}$
(resp. acting on the cohomology of the base-changed variety to an
algebraic closure of $\Fp_{\ell}$). By Deligne's Riemann Hypothesis
(see, e.g.,~\cite[Th. 11.37]{ant}), we have
\begin{align*}
\sum_{g\in \Gg_{i,j}(\Fp_{\ell})}
{\Bigl(\frac{g_{i,j}}{\ell}\Bigr)}
&\ll q^{d}\dim H^{2d}_c(\overline{\G}_{i,j},\mathcal{K})+
q^{d-1/2}\sum_{k<2d}{\dim H^{k}_c(\overline{\G}_{i,j},\mathcal{K})}
\\
&
\ll q^{d}\dim H^{2d}_c(\overline{\G}_{i,j},\mathcal{K})+
q^{d-1/2}
\end{align*}
for $\ell\geq 3$, by results of Bombieri or Adolphson--Sperber that
show that the sum of dimensions of cohomology groups is bounded
independently of $\ell$ (see, e.g.,~\cite[Th. 11.39]{ant}).
\par
It therefore remains to prove that
$H^{2d}_c(\overline{\G}_{i,j},\mathcal{K})=0$. However, this space is
isomorphic (as vector space) to the space of coinvariants of the
geometric fundamental group of $\Gg_{i,j}$ acting on a one-dimensional
space through the character which ``is'' the Lang-Kummer sheaf
$\mathcal{K}$. This means that either the coinvariant space is zero,
and we are done, or otherwise the sheaf is geometrically trivial.  The
latter translates to the fact that the traces on $\mathcal{K}$ of the
local Frobenius 
$\frob_{g,\ell^{\nu}}$ of rational points
$g\in\Gg_{i,j}(\Fp_{\ell^{\nu}})$ over all extensions fields
$\Fp_{\ell^{\nu}}/\Fp_{\ell}$ depend only on $\nu$, i.e.,  the map
$$
g\mapsto \Bigl(\frac{N_{\Fp_{\ell^{\nu}}/\Fp_{\ell}}g_{i,j}}{\ell}
\Bigr)
$$
on $\Gg_{i,j}(\Fp_{\ell^{\nu}})$ depends only on  $\nu$. But 
this is clearly impossible for $SL(n)$ or $Sp(2g)$ with $n\geq 2$,
$g\geq 1$ (but not for $SL(1)$ or for $SL(2,\Fp_2)$...), because we
can explicitly write down  matrices even in $\Gg(\Fp_{\ell})$ both
with $g_{i,j}$ a non-zero square and $g_{i,j}$ not a square (taking
$\ell\geq 3$ for $SL(2,\Fp_{\ell})$).   
\par
(3) and (4): those are similar to (1). Namely, define first a 
 $q$-symplectic polynomial $f$ in $\Fp_{\ell}[X]$ to be one
of degree $2g$ such that\footnote{\ Unfortunately, this is not stated
  correctly   in~\cite{sieve}, although none of the results there are
  affected by 
  this slip...}
$$
f(0)=1,\quad\quad\text{ and }\quad\quad
(qT)^{2g}f(1/(qT))=f(T).
$$
We can express such a $q$-symplectic polynomial uniquely in the form
\begin{multline*}
f(T)=1+a_1(f)T+\cdots+a_{g-1}(f)T^{g-1}+
a_{g}(f)T^g+\\
qa_{g-1}(g)T^{g+1}+\cdots +
q^{g-1}a_1(f)T^{2g-1}+q^{g}T^{2g},
\end{multline*}
with $a_i(f)\in \Fp_{\ell}$, and this expression gives a bijection
$$
f\mapsto (a_1(f),\ldots,a_g(f))
$$
between the set of $q$-symplectic polynomials  and
$\Fp_{\ell}^g$. 
\par
Then we need to bound
\begin{equation}\label{eq-34-1}
\frac{1}{|Sp(2g,\Fp_{\ell})|}
\sum_{f\in \Omega^{\gamma}}
{
|\{g\in G\,\mid\, \det(1-Tg)=f\}|
}
\end{equation}
where we have put (in case (3) and (4) respectively)
\begin{align*}
\Omega^{(3)}&=\{f \in\Fp_{\ell}[T]\,\mid\, f\text{ is
  $q$-symplectic  and } f(1)\text{ is a
  square in } \Fp_{\ell}\},\\
\Omega^{(4)}&=\{f\in \Fp_{\ell}[T]\,\mid\, f\text{ is
  $q$-symplectic and } q+1-a_1(f)\text{ is a
  square in } \Fp_{\ell}\}.
\end{align*}
\par
Now it is  easy to check that we have
\begin{equation}\label{eq-34-2}
|\Omega^{\gamma}|=\frac{\ell^g+\ell^{g-1}}{2}\geq \frac{\ell^g}{2}
\end{equation}
for $\gamma=3$ or $4$ (recall $\ell$ is odd). Indeed, treating the
case $\gamma=3$ (the other is similar), we have
$$
|\Omega^{(3)}|=|\{f\,\mid\, f(1)=0\}|
+\frac{1}{2}\sum_{f(1)\not=0}{\Bigl(1+\Bigl(\frac{f(1)}{\ell}\Bigr)\Bigr)}.
$$
\par
The first term is $\ell^{g-1}$ since $f\mapsto f(1)$ is a non-zero
linear functional 
on $\Fp_{\ell}^g$. The first part of the second sum is
$(\ell^g-\ell^{g-1})/2$, and the last is 
$$
\sum_{(a_2,\ldots,a_{g})}
\sum_{a_g\not=-\tilde{f}(1)}{\Bigl(\frac{a_g+\tilde{f}(1)}{\ell}\Bigr)}
$$
where $\tilde{f}(1)$ is defined by $f(1)=a_g+\tilde{f}(1)$ (note that
$\tilde{f}(1)$ depends only on $(a_2,\ldots,a_g)$). Because of the
summation over the free variable $a_g$, this expression vanishes. 
\par
Now appealing to Lemma~7.2 of~\cite{sieve} (itself derived from the
work of Chavdarov), we obtain
\begin{equation}\label{eq-34-3}
\frac{1}{|Sp(2g,\Fp_{\ell})|}
|\{g\in G\,\mid\, \det(1-Tg)=f\}|\geq 
\frac{1}{(\ell+1)^g}
\end{equation}
for all $q$-symplectic polynomials $f$, and hence the stated bound
follows by~(\ref{eq-34-1}),~(\ref{eq-34-2}),~(\ref{eq-34-3}).
\par
(5) and (6): this is again similar to (3) and (4), where we now deal
with 
$$
\frac{1}{|Sp(2g,\Fp_{\ell})|}
\sum_{f\in \Omega^{\gamma}}{
|\{g\in G\,\mid\, \det(1-Tg)=f\}|
}
$$
with now
\begin{align*}
\Omega^{(5)}&=
\{f \in\Fp_{\ell}[T]\,\mid\, f\text{ is
  $q$-symplectic  and } f(1)=0\},\\
\Omega^{(6)}&=\{f\in \Fp_{\ell}[T]\,\mid\, f\text{ is
  $q$-symplectic and } q+1=a_1(f)\}.
\end{align*}
\par
We have in both cases $|\Omega^{\gamma}|=\ell^{g-1}$
since the condition is a linear one on the coefficients. By the proof
of  Lemma~7.2 of~\cite{sieve} we also have
$$
|\{g\in G\,\mid\, \det(1-Tg)=f\}|
\leq \frac{1}{(\ell-1)^g}
$$
for all $f$, and therefore 
$$
\frac{1}{|Sp(2g,\Fp_{\ell})|}
\sum_{f\in \Omega^{\gamma}}{
|\{g\in G\,\mid\, \det(1-Tg)=f\}|
}\leq \frac{\ell^{g-1}}{(\ell-1)^g}.
$$
Since the quantity to estimate is also at most $1$ for trivial
reasons, we have the desired result.
\end{proof}

\end{document}